\documentclass[a4paper,10pt]{article}
\everymath{\displaystyle}
\usepackage{amssymb,amsmath,amsthm}
\usepackage{graphicx}

\def\N{\mathbb N}
\def\R{\mathbb R}
\def\T{\mathbb T}
\def\Zi{{\mathbb Z}\setminus {\cal I}}
\usepackage{amsfonts}
\usepackage{hyperref}
\usepackage{latexsym}
\usepackage{amsthm}
\usepackage{amsmath}
\usepackage{amssymb}
\usepackage{amsopn}
\usepackage{geometry}
\usepackage{amscd}
\usepackage{mathrsfs}
\usepackage{dsfont}
\usepackage[english]{babel}
\usepackage{caption}
\font\teneufm=eufm10
\font\seveneufm=eufm7
\font\fiveeufm=eufm5
\newfam\eufmfam
\textfont\eufmfam=\teneufm
\scriptfont\eufmfam=\seveneufm
\scriptscriptfont\eufmfam=\fiveeufm

\newcommand{\normastretta}{|\!\!|}

\newcommand{\Om}{\Omega}

\newenvironment{pf}{\noindent{\sc Proof}.\enspace}{\rule{2mm}{2mm}\smallskip}
\newenvironment{pfn}{\noindent{\sc Proof} \enspace}{\rule{2mm}{2mm}\smallskip}
\newtheorem{theorem}{Theorem}[section]
\newtheorem{proposition}{Proposition}[section]
\newtheorem{lemma}{Lemma}[section]
\newtheorem{corollary}{Corollary}[section]

\newtheorem{remark}{Remark}[section]
\newtheorem{remarks}{Remark}[section]
\newtheorem{definition}{Definition}[section]
\newcommand{\be}{\begin{equation}}
\newcommand{\ee}{\end{equation}}
\newcommand{\teta}{\theta}
\newcommand{\om}{\omega}
\newcommand{\e}{\varepsilon}

\renewcommand{\a }{\alpha }
\renewcommand{\b }{\beta }
\newcommand{\s }{\sigma }
\newcommand{\ii }{{\rm i} }
\renewcommand{\d }{\delta }
\newcommand{\D }{\Delta}
\def\Ta{\mathcal T}
\newcommand{\g }{\gamma}

\renewcommand{\l }{\lambda }
\renewcommand{\L }{\Lambda }
\newcommand{\m }{\mu }

\renewcommand{\t }{\tau }
\renewcommand{\o }{\omega }
\renewcommand{\O }{\Omega }
\newcommand{\C}{\mathbb{C}}
\newcommand{\Z}{\mathbb{Z}}

\newcommand{\mm}{{\rm m}}
\newcommand{\fr}{{\hat f}}
\newcommand{\ft}{{\tilde f}}
\newcommand{\Fr}{{\hat F}}
\newcommand{\Ft}{{\tilde F}}

\newcommand{\pluto}{\mathtt j}
\newcommand{\norma}{\|}
\newcommand{\Hsr}{\mathcal{H}_{s,r}}
\newcommand{\Hsrn}{\mathcal{H}_{s,r}^{\rm null}}
\newcommand{\Hrn}{\mathcal{H}_{r}^{\rm null}}
\newcommand{\rro}{{\rho}}

\newcommand{\cc}{{1}}
\newcommand{\CC}{{6}}
\newcommand{\dueCC}{{12}}
\newcommand{\duecc}{{2}}
\newcommand{\trecc}{{3}}
\newcommand{\quattrocc}{{4}}

\newcommand{\CCterzi}{{2}}
\newcommand{\CCquarti}{{3/2}}
\newcommand{\CCNmezzi}{{3N}}
\newcommand{\CCNquarti}{{3N/2}}
\newcommand{\quattrofrattoCC}{{\frac23}}

\newcommand{\treCC}{{18}}

\newcommand{\ccN}{{N}}

\usepackage{color}

\title{\bf KAM theory for the \\
Hamiltonian derivative wave equation}

\author{Massimiliano Berti, Luca Biasco, Michela Procesi}

\date{}

\begin{document}

\maketitle

\noindent
{\bf Abstract:} We prove an infinite dimensional KAM theorem
 which implies  the existence of
Cantor families of small-amplitude,
 reducible, elliptic, analytic, invariant tori
of  Hamiltonian derivative wave equations.
\\[1mm]
{\sc 2000AMS subject classification}: 37K55, 35L05.




\section{Introduction}

In the last years many progresses have been done concerning  KAM theory for nonlinear Hamiltonian PDEs.
The first  existence
results  
were given
by Kuksin \cite{Ku} and Wayne \cite{W1} for semilinear
wave (NLW) and Schr\"odinger equations (NLS) in one space dimension ($1d$) under Dirichlet boundary conditions,
see  \cite{Po2}-\cite{Po3} and \cite{KP} for further developments.
The approach of these papers consists  in generating iteratively
a sequence of symplectic changes of variables 
 which bring the Hamiltonian into a  constant coefficients ($=$reducible) normal form  with an elliptic  ($=$linearly stable)
invariant torus at the origin.
Such a  torus is filled by quasi-periodic solutions with zero Lyapunov exponents.
This procedure requires
 to solve, at each step, constant-coefficients linear ``homological equations"
by imposing the ``second order Melnikov" non-resonance conditions.
Unfortunately these (infinitely many)
conditions are violated
already for periodic boundary conditions.

In this case, existence of quasi-periodic solutions for semilinear $ 1d $-NLW and NLS equations,
was first proved by Bourgain \cite{B1} by extending the Newton approach
introduced by Craig-Wayne \cite{CW} for periodic solutions.
Its main advantage 
is to require only the 
``first order Melnikov"
non-resonance conditions
(the minimal assumptions)  for solving the
homological equations.
Actually, developing this perspective,
Bourgain  
 was able to
prove in \cite{B3}, \cite{B5}  also
the existence of
quasi-periodic solutions  for NLW and NLS (with Fourier multipliers) in higher space dimensions,
see also the recent extensions  in \cite{BBo}, \cite{W2}.
The main drawback
of this approach is that the
homological equations are  linear PDEs with non-constant coefficients.
Translated in the KAM language
 this implies
a non-reducible normal form around the torus and then a lack of  informations about the
stability of the quasi-periodic solutions.

Later on, existence of reducible elliptic  tori was proved by  Chierchia-You  \cite{CY} for semilinear $ 1d $-NLW,
and, more recently,  by   Eliasson-Kuksin \cite{EK} for
NLS  (with Fourier multipliers) in any space dimension,
see also Procesi-Xu \cite{PX}, Geng-Xu-You \cite{GYX}.

\smallskip

An important problem concerns the study of PDEs where the
nonlinearity involves derivatives. A comprehension of this
situation is of major importance since most of the models coming
from Physics are of this kind.

\smallskip

In this direction KAM theory 
has been extended to deal with KdV  equations
 by Kuksin \cite{K2}-\cite{k1},
Kappeler-P\"oschel \cite{KaP},  and,
for the $1d$-derivative NLS (DNLS) and Benjiamin-Ono equations,  by Liu-Yuan \cite{LY}.
The key idea of these results is again to provide only a non-reducible normal form around the torus.
However, in this cases, the  homological
equations with  non-constant coefficients are only  {\it scalar} (not an infinite system as in the Craig-Wayne-Bourgain
approach).
We remark that the KAM proof 
 is more delicate for DNLS and Benjiamin-Ono, because these equations
 are less ``dispersive" than KdV, i.e.
  the eigenvalues of the principal part of the differential operator 
  grow only quadratically  at infinity, and not cubically as  for KdV.
As a consequence of this difficulty,  
the quasi-periodic solutions 
in \cite{K2}, \cite{KaP} are analytic, in  \cite{LY}, only $ C^\infty $.
 Actually, for the applicability of these KAM schemes,
the more dispersive the equation is,
the more derivatives in the nonlinearity  can be supported. 
The limit case 
of the derivative nonlinear wave equation (DNLW) -which is not dispersive at all-  is excluded  by these approaches.

\smallskip

In the paper \cite{B1} (which proves the existence
of quasi-periodic solutions for semilinear $ 1 d $-NLS and NLW),
Bourgain  claims, in the last remark, that his analysis works also
for the 
Hamiltonian ``derivation" wave equation
$$
y_{tt} - y_{xx} + g(x)y  = \Big(  - \frac{d^2}{d x^2} \Big)^{1/2} F(x,y) \, ,
$$
see also \cite{B2}, page 81. Unfortunately no details are  given.
However, Bourgain  \cite{B2} provided
a detailed proof
of the existence of periodic solutions for the non-Hamiltonian equation
$$
y_{tt} - y_{xx} + \mm y + y_t^2  = 0 \, , \quad \mm \neq 0 \, .
$$
These kind of problems have been then reconsidered
by Craig in \cite{C}  for  more general Hamiltonian derivative wave equations like
$$
y_{tt} - y_{xx} + g(x) y + f(x, D^\b y)=0 \, , \quad x \in \T \, ,
$$
where $ g(x) \geq 0 $ and $ D $ is the first order pseudo-differential operator
$ D:=\sqrt{-\partial_{xx}+ g(x)} $.
The perturbative analysis of Craig-Wayne \cite{CW} for the search of periodic solutions
works when
$ \b < 1 $.
The main reason is that the wave equation vector field gains one derivative and then
the 
nonlinear term  $ f(D^\b u) $  has a strictly weaker effect on the dynamics for  $ \b < 1 $.
The case $ \b = 1 $ is left as an open problem. Actually,  in this case,
the small divisors problem for periodic solutions has the same level of difficulty of quasi-periodic solutions
with $ 2 $ frequencies.

\smallskip

The goal of this paper is to extend KAM theory to deal with the
Hamiltonian  derivative wave equation
\begin{equation}\label{pseudoNLW}
y_{tt} - y_{xx} + \mm y + f(D y)=0 \, , \quad \mm >  0 \, , \quad D :=\sqrt{-\partial_{xx}+\mm} \, , \quad x \in \T \, ,
\end{equation}
with real analytic nonlinearities (see Remark \ref{carmelo})
\be\label{funzi}
f(s) = a s^3 + \sum_{k\geq 5} f_k s^k \, , \quad a \neq 0 \, .
\ee
We  write equation \eqref{pseudoNLW}
as
the  infinite dimensional  Hamiltonian system
$$
u_t=-\ii \partial_{\bar u}H\,,\qquad
\bar u_t=\ii \partial_{u}H\, ,
$$
with Hamiltonian
\begin{equation}\label{bilbo}
H (u,\bar u) : =\int_\T
\bar u D u+ F\Big(\frac{u+\bar u}{\sqrt 2}\Big)
\,dx\,, \quad F(s) := \int_0^s f \, ,
\end{equation}
in the 
complex unknown
$$
u:=\frac{1}{\sqrt 2} (Dy+\ii y_t)\,,\qquad
\bar u:=\frac{1}{\sqrt 2} (Dy-\ii y_t)\,,\qquad
\ii:=\sqrt{-1} \, .
$$
Setting
$ u = \sum_{j\in \Z} u_j e^{\ii jx} $ (similarly for $ \bar u$),
we obtain the Hamiltonian in infinitely many coordinates
\begin{equation}\label{tuc}
H = 
\sum_{j\in\Z} \l_j u_j \bar u_j +
\int_\T
F\Big(\frac{1}{\sqrt 2}\sum_{j\in\Z}(u_je^{\ii j x}+\bar u_je^{-\ii j x})\Big)
\,dx   
\end{equation}
where
\be\label{faramir}
 \l_j:=\sqrt{j^2+\mm}
\ee are the eigenvalues of the diagonal operator $D$. Note that
the nonlinearity in \eqref{pseudoNLW} is $ x $-independent
implying,  for \eqref{bilbo}, the conservation of the momentum $ -
\ii \int_{\T} \bar u \partial_x u \, dx $. This symmetry allows to
simplify somehow  the KAM proof (a similar idea was used by
Geng-You \cite{GY}).

\smallskip

For every choice of   the \textsl{tangential sites}
${\cal I} := \{\pluto_1,\dots, \pluto_n \} \subset \Z $, $n\geq 2 $,
the integrable Hamiltonian $\sum_{j\in\Z}\l_j u_j \bar u_j$
has the invariant tori
$ \{ u_j\bar u_j=\xi_j,\ {\rm for}\ j\in\mathcal{I}\,,\ u_j=\bar u_j=0
\ {\rm for}\ j\not\in\mathcal{I}  \}$
parametrized by the actions $ \xi = (\xi_j)_{j\in\mathcal{I}}\in\R^n $.
The  next KAM result states the existence of nearby
invariant tori 
for the complete Hamiltonian $H$ in \eqref{tuc}.


\begin{theorem}\label{thm:DNLW}
The equation \eqref{pseudoNLW}-\eqref{funzi}
admits
Cantor families of small-amplitude, analytic,  quasi-periodic solutions 
  with zero Lyapunov exponents
and whose linearized equation is reducible to constant coefficients.
Such Cantor families have 
asymptotically full measure at the origin
 in  the set of parameters. 
\end{theorem}

\smallskip

The proof of Theorem \ref{thm:DNLW} is based on the abstract 
KAM Theorem \ref{thm:IBKAM},
which provides a reducible normal form (see \eqref{Hnew})
around the elliptic invariant torus, and on the measure estimates
Theorem \ref{thm:measure}.
The key point in proving Theorem \ref{thm:measure} is the asymptotic bound \eqref{freco}
on the perturbed normal frequencies $ \Omega^\infty (\xi ) $
after the  KAM iteration. This allows to prove
that the
second order Melnikov non-resonance conditions \eqref{Cantorinf}
are fulfilled for an asymptotically full measure set of parameters
(see \eqref{consolatrixafflictorum}). The estimate \eqref{freco}, in turn,
 is achieved by exploiting  the {\it quasi-T\"oplitz}  property
 of the perturbation.
This notion 
has been introduced by Procesi-Xu \cite{PX} in the context of NLS in higher space dimensions
and it is similar, in spirit, to the T\"oplitz-Lipschitz property in 
Eliasson-Kuksin \cite{EK}.
The precise
formulation of quasi-T\"oplitz functions,
adapted to the DNLW setting,  is given in  Definition \ref{topbis_aa}
below.

\smallskip

Let us roughly explain
the  main ideas and techniques for proving
Theorems \ref{thm:IBKAM}, \ref{thm:measure}. These theorems concern, as usual,
a parameter dependent family of analytic Hamiltonians of the form
\begin{equation}\label{HKAM}
H = \o (\xi)\cdot y + \O(\xi)\cdot z \bar z + P(x,y,z,\bar z;\xi)
\end{equation}
where $ (x, y)  \in \T^n \times \R^n $,  $ z, \bar z  $ are infinitely many variables,
$ \o(\xi)\in\R^n $, $\O(\xi)\in\R^\infty $ and $\xi\in\R^n$.
The frequencies $ \Om_j (\xi)$ are close to the unperturbed frequencies $ \l_j $ in \eqref{faramir}.

As well known,
the main difficulty of the KAM iteration
which provides a reducible KAM normal form like \eqref{Hnew}
is to fulfill, at each iterative step, the 
second order Melnikov non-resonance conditions.
Actually, following the formulation of the KAM theorem given in \cite{BB10}, it is sufficient
to verify 
\begin{equation}\label{cervicale}
| \o^{\infty} (\xi) \cdot k +\O_i^{\infty}(\xi) - \O_j^{\infty}(\xi)|
\geq \frac{\g}{1+|k|^\tau} \, , \quad  \g > 0 \, , 
\end{equation}
only for the ``final" frequencies $ \om^\infty (\xi) $ and
$ \Om^\infty (\xi) $, see \eqref{Cantorinf}, and
not along the inductive iteration.


The application of the usual KAM theory (see e.g. \cite{Ku},
\cite{Po2}-\cite{Po3}), to the DNLW equation
provides  only the asymptotic decay estimate
\be\label{DNLW}
\O_j^\infty(\xi)= j + O(1)  \quad {\rm for} \quad j \to + \infty \, . 
\ee
Such a bound is not enough: the set of parameters $ \xi $ satisfying \eqref{cervicale}
could be empty. 
Note that for the semilinear NLW equation (see e.g. \cite{Po2}) the frequencies
decay asymptotically faster, namely like $ \O_j^\infty(\xi) = j + O (1 / j ) $.

\smallskip

The key idea for verifying the second order Melnikov non-resonance conditions \eqref{cervicale} for DNLW
is to  prove the higher order asymptotic decay estimate (see \eqref{freco}, \eqref{carbonara})
\begin{equation}\label{delfi}
   \O_j^\infty(\xi) =
    j + a_+(\xi) +\frac{\mm}{2j}+  O (\frac{\g^{2/3}}{j})
    \quad {\rm for} \quad j \geq O(\g^{-1/3})
   \end{equation}
where $ a_+(\xi) $ is a constant independent of $ j $
(an analogous expansion holds  for $ j \to - \infty $
with a possibly different limit constant  $ a_- (\xi) $).
In this way infinitely many conditions in \eqref{cervicale} are
verified by imposing only
 first order Melnikov conditions like
$ |\o^\infty (\xi) \cdot k + h | \geq  2 \g^{2/3} / |k|^\tau $,
$ h \in \Z $.
Indeed,
for $ i > j >   O(|k|^\t \g^{-1/3})  $, we get
\begin{eqnarray*}
   |\o^\infty (\xi) \cdot k +\O_i^\infty(\xi)-\O_j^\infty(\xi)|
   &=&
|\o^\infty (\xi) \cdot k + i - j +\frac{\mm(i-j)}{2ij} +O(\g^{2/3}/j)|
\\
&\geq&   2 \g^{2/3}  |k|^{-\tau} - O(|k|/j^2) -O(\g^{2/3}/j)  \geq
\g^{2/3}  |k|^{-\tau}
\end{eqnarray*}
noting that $ i - j $ is integer and $|i-j|=O(|k|)$
(otherwise no small divisors occur).
 We refer to section \ref{sec:meas} for the precise arguments,
 see in particular
Lemma \ref{lem:inclu}.

\smallskip

The asymptotic decay \eqref{freco}  for the perturbed frequencies $ \Omega^\infty (\xi) $
 is achieved thanks to
 the ``quasi-T\"oplitz"  property of the perturbation
 (Definition \ref{topbis_aa}).
 Let us roughly explain  this notion.
The new
normal frequencies after each KAM step
are $ \Om_j^+ = \Om_j + P_{j}^0  $ where  the corrections $ P_{j}^0 $
are the coefficients of  the quadratic form
$$
P^0 z \bar z := \sum_j P_{j}^0 z_j \bar z_j \, , \quad  \
P_{j}^0 
:= \int_{\T^n} (\partial^2_{z_j\bar z_j}  P)(x,0, 0, 0; \xi) \, d x \, .
$$
We say that a quadratic form $ P^0 $ is quasi-T\"oplitz if it has the form
$$
P^0 = T + R
$$
where $T$ is a T\"oplitz matrix (i.e. constant on the diagonals)
and $R$ is a ``small" remainder
satisfying $ R_{jj} = O(1/ j)$ (see Lemma \ref{frodo}).
Then  \eqref{delfi} follows with $ a := T_{jj} $ which is independent of $ j $.

\smallskip
Since the quadratic perturbation $P^0$ along the KAM iteration
 does not depend only on the quadratic perturbation at
 the previous steps, we need to extend the notion of quasi-T\"oplitz to
 general (non-quadratic) analytic functions.

The preservation of the quasi-T\"oplitz property
of the perturbations $ P $ at each KAM step (with just slightly modified parameters)
holds in view of   the following
key facts:
\begin{enumerate}
\item
the Poisson bracket of two quasi-T\"oplitz functions is quasi-T\"oplitz (Proposition \ref{festa}),
\item the hamiltonian flow generated by
a quasi-T\"oplitz function preserves the quasi-T\"oplitz property
 (Proposition \ref{main}),
 \item
the solution of the homological equation with a quasi-T\"oplitz  perturbation
is quasi-T\"oplitz (Proposition \ref{FP}).
\end{enumerate}
We note that, in \cite{EK}, the analogous properties 1  (and therefore 2)
for T\"oplitz-Lipschitz functions is proved only when one of them 
is quadratic.

The definition of
quasi-T\"oplitz functions 
heavily relies on properties of
projections.
However,  for an analytic function in infinitely many variables, 
such projections may not be well defined unless the
Taylor-Fourier series
(see \eqref{hatf}) is \textsl{absolutely} convergent.
 For such reason, instead of the $\sup$-norm,
  we use the majorant norm (see \eqref{normadue}, \eqref{normadueA}),
  for which
 the bounds
 \eqref{proiezA} and  \eqref{caligola} on projections hold
(see also Remark \ref{failsup}).

We underline that the majorant norm
of a vector field introduced  in \eqref{normadueA} is very different
from the weighted norm introduced by P\"oschel in
\cite{Po}-Appendix C, which works {\it only} in finite
dimension, see comments in \cite{Po} after Lemma C.2 and
Remark \ref{differencePoschel}.
As far as we know this majorant norm
of  vector fields is new. In Section \ref{sec:2}
we show its properties, in particular the key
estimate of
  the majorant norm of the commutator of
two vector fields (see Lemma \ref{settimiosevero}).



\smallskip

Before concluding this introduction we also mention
the recent KAM theorem of
Greb\'ert-Thomann \cite{GT} for the
quantum harmonic oscillator with semilinear nonlinearity.
Also here the eigenvalues grow to infinity only linearly. 
We quote the
 normal form results 
of Delort-Szeftel \cite{DS}, Delort \cite{De},
for quasi-linear wave equations, 
where only finitely many steps of normal form can be performed.
Finally we mention also the recent work by G\'{e}rard-Grellier
\cite{GG} on Birkhoff normal form for a degenerate ``half-wave''
 equation.

\smallskip

The paper is organized as follows:
\begin{itemize}
\item
In {\sc section \ref{sec:2}} we define the {\it majorant} norm of formal  power series of scalar functions
(Definition \ref{majsc}) and vector fields (Definition \ref{MNV}) and
we  investigate the relations with the notion of analiticity,
see Lemmata \ref{fundamA}, \ref{mare}, \ref{ttata}, \ref{XMvectA} and  Corollary \ref{cor11}.
Then we prove Lemma \ref{settimiosevero} on commutators.
\item
In {\sc section \ref{sec:3}}
we define the  T\"oplitz (Definition \ref{matteo_aa}) and Quasi-T\"oplitz functions (Definition \ref{topbis_aa}).
Then we prove that this class of functions is closed under Poisson brackets
(Proposition \ref{festa}) and composition with the Hamiltonian flow (Proposition \ref{main}).

\item In {\sc section \ref{sec:4}} we state the abstract KAM Theorem \ref{thm:IBKAM}.
The first part of  Theorem \ref{thm:IBKAM} follows by the KAM Theorem 5.1  in \cite{BB10}.
The main novelty is part II, in particular the asymptotic estimate  \eqref{freco} of the normal frequencies.

\item
In {\sc section \ref{sec:5}}  we prove the abstract KAM Theorem \ref{thm:IBKAM}.

We first perform (as in Theorem 5.1  in \cite{BB10}) a first normal form step, which makes
Theorem \ref{thm:IBKAM} suitable for the direct application to the wave equation.

In Proposition \ref{FP} we prove that the solution of the homological equation with a quasi-T\"oplitz  perturbation
is quasi-T\"oplitz.  Then the main results of the KAM step concerns
the asymptotic estimates of the perturbed frequencies (section \ref{sec:normalf})
and the T\"oplitz estimates of the new perturbation (section \ref{sec:P+}).

\item   In {\sc section \ref{sec:meas}}  we prove Theorem \ref{thm:measure}:  the second order Melnikov non-resonance
conditions are fulfilled for a set of parameters with large measure,
see \eqref{consolatrixafflictorum}.
We  use the conservation of momentum to avoid the presence of double eigenvalues.
\item In {\sc section \ref{sec:DNLW}} we  finally
apply the abstract KAM Theorem \ref{thm:IBKAM} to the DNLW equation
\eqref{pseudoNLW}-\eqref{funzi}, proving Theorem \ref{thm:DNLW}.
We first verify that the Hamiltonian  \eqref{tuc} is quasi-T\"oplitz (Lemma \ref{finzioni}), as well as the
Birkhoff normal form Hamiltonian \eqref{Birkhoff} of Proposition \ref{BNF}.
The main technical difficulties concern
 the proof in Lemma \ref{isaia}
that the generating function
\eqref{brodino} of the Birkhoff symplectic transformation is also quasi-T\"oplitz
 (and the small divisors Lemma \ref{poeschel}).
In section \ref{sec:aa} we prove that the perturbation, obtained
after the introduction of the action-angle variables,
is still quasi-T\"oplitz (Proposition \ref{qtop0}). Finally  in section \ref{sec:pr} we prove Theorem
\ref{thm:DNLW} applying Theorems \ref{thm:IBKAM} and \ref{thm:measure}.
\end{itemize}

\vskip10pt
\noindent{\bf  Acknowledgments :}  We thank Benoit Gr\'ebert for pointing out a technical mistake in the previous version.

\section{Functional setting}\label{sec:2}\setcounter{equation}{0}

Given a finite subset $ \cal I \subset \Z $  (possibly empty),  $ a \geq 0, p > 1/2 $,
we define the Hilbert space
$$
\ell^{a,p}_{\cal I} := \Big\{  z=\{ z_j\}_{j\in\Z\setminus {\cal I}} \, , \ z_j \in \C \ \ : \
\| z \|_{a,p}^2 := \sum_{j\in\Z\setminus {\cal I}} |z_j|^2 e^{2a|j|} \langle j
\rangle^{2p} < \infty  \Big\} \, .
$$
When $ \cal I = \emptyset $ we denote $ \ell^{a,p} := \ell^{a,p}_{\cal I}  $. We consider the direct product 
\begin{equation}\label{E}
E:= \C^n \times \C^n \times \ell^{a, p}_{\cal I} \times \ell^{a, p}_{\cal I} 
\end{equation}
where $ n $ is the cardinality of $ \cal I $.
We endow  the space $ E $  with the $(s,r)$-weighted norm 
\begin{equation}\label{normaEsr}
v =  (x,y,z,\bar z) \in E \, , \quad
\|v\|_E:=\|v\|_{E,s,r}= \frac{|x|_\infty}{s} + \frac{|y|_1}{r^2}
 +\frac{\|z\|_{a,p}}{r}+\frac{\|\bar z\|_{a,p}}{r}
\end{equation}
where, $ 0 < s, r < 1 $, and 
$ |x|_\infty := \max_{h =1, \ldots, n} |x_h| $,
$ |y|_1 := \sum_{h=1}^n |y_h| $.
Note that, for all $s'\leq s$, $ r' \leq r $,
\begin{equation}\label{empire}
\| v\|_{E,s',r'}\leq \max \{ s / s', (r / r' )^2 \}\| v\|_{E,s,r}\,.
\end{equation}
We shall also use the notations
$$
z_j^+ = z_j \, , \quad z_j^- = \bar z_j \, .
$$
We identify a vector $ v \in E $ 
with the sequence $ \{ v^{(j)} \}_{ j \in {\cal J} } $ with indices in 
\be\label{J} {\cal J}:=
 \left\{ j=(j_1,j_2),\  j_1 \in \{1,2,3,4\}, \
 j_2 \in
\begin{cases}
\{1, \ldots,  n\} \quad \ \; \quad \  {\rm if} \   j_1 = 1, 2 \cr  \Zi  \qquad
\quad \qquad{\rm if} \ j_1 = 3,4
 \end{cases} \right\}
\ee
and components
$$
v^{(1,j_2)}:=x_{j_2}\,, \  v^{(2,j_2)}:=y_{j_2}\, \ (1\leq j_2\leq
n), \ v^{(3,j_2)}:=z_{j_2}\,,  \  v^{(4,j_2)}:=\bar z_{j_2} \
(j_2\in\Zi ) \, ,
$$
more compactly 
$$
v^{(1,\cdot)}:=x \, , \
v^{(2,\cdot)}:=y, \, , \
v^{(3,\cdot)}:=z, \, , \
v^{(4,\cdot)}:=\bar z \, .
$$
We denote by $\{ e_j\}_{j\in{\cal J}}$
the orthogonal basis of the Hilbert space $E$,
where
$e_j$ is the sequence with all zeros, except the $j_2$-th
entry of its $j_1$-th components, which is $ 1 $.
Then every $v\in E$ writes 
$ v=\sum_{j\in {\cal J}} v^{(j)} e_j $, $  v^{(j)}\in\C $.
We also define the toroidal domain
\be\label{Dsr}
D (s,r) := \T^n_s \times D(r) := \T^n_s \times B_{r^2} \times B_r
\times  B_r \subset E
\ee
where $ D(r):= B_{r^2} \times B_r \times B_r $,
\be\label{seconda}
 \T^n_s := \Big\{ x \in \C^n \, : \, \max_{h=1, \ldots, n} |{\rm Im} \, x_h | < s  
 \Big\} \, , \ \
 B_{r^2}  := \Big\{ y \in \C^n \, : \, |y |_1 < r^2 \Big\}
\ee
and $ B_{r} \subset \ell^{a,p}_{\cal I} $ is the open ball of radius $ r $
centered at zero. We think $\T^n$ as the $ n $-dimensional torus
$ \T^n := 2\pi \R^n/\Z^n $, namely $ f:D(s,r) \to \C $
means that $f$ is $ 2\pi$-periodic in each $ x_h $-variable, $h=1,\ldots,n$.
\begin{remark}
If $ n = 0 $  then $ D(s,r) 
 \equiv B_r\times B_r \subset
\ell^{a,p}\times \ell^{a,p}$.
\end{remark}

\subsection{Majorant norm}

\subsubsection{Scalar functions}

We consider {\sl formal}  power series with infinitely many variables
\begin{equation}\label{formalpowerA}
f(v)=f(x,y,z, \bar z) = \sum_{(k,i,\a,\b) \in {\mathbb I}} f_{k,
i, \alpha, \b} \, e^{\ii k \cdot x} y^i  z^\alpha \bar z^{\b}
\end{equation}
with coefficients $ f_{k, i, \alpha, \b}\in \C $ and multi-indices in
\be\label{I}
\mathbb{I}:=\Z^n \times \N^n \times \N^{(\Zi )} \times \N^{(\Zi )}
\ee
where
\be\label{finNZ}
\N^{(\Z \setminus {\cal I})}  := \Big\{ \a := (\a_j)_{j \in \Z \setminus {\cal I}} \in \N^\Z \ \ {\rm with} \ \
|\a| := \sum_{j \in \Z \setminus {\cal I}} \a_j < + \infty \Big\} \, .
\ee
In \eqref{formalpowerA} we use the standard multi-indices notation
$ z^\a {\bar z}^\b := \Pi_{j \in \Z \setminus {\cal I}} \, z_j^{\a_j}  {\bar z}_j^{\b_j} $.
We denote 
the monomials
\be\label{commodo}
{\mathfrak m}_{k,i,\a,\b}(v) =
{\mathfrak m}_{k,i,\a,\b}(x,y,z,\bar z):= e^{\ii k\cdot x} y^i z^\a {\bar z}^{\b}\, .
\ee
\begin{remark}
If $n=0$ the set $\mathbb I$ reduces to $ \N^{\Z } \times \N^{\Z }$ and the formal series to
$ f(z, \bar z) = \sum_{(\a,\b) \in {\mathbb I}} f_{\alpha, \b} \,  z^\alpha \bar z^{\b} $.
\end{remark}
We  define the ``majorant" of $ f $ as
\be\label{formalpower1A}
\big( Mf\big) (v)=\big( Mf\big)(x,y,z, \bar z) :=
\sum_{(k,i,\a,\b) \in {\mathbb I}} |f_{k, i, \alpha, \b}| e^{\ii k\cdot x } y^i z^\alpha {\bar z}^{\b} \,.
\ee
We now discuss the convergence of  formal series.
\begin{definition}\label{augusto}
A series
$$
\sum_{(k, i, \alpha, \b )\in\mathbb{I}  } c_{k, i, \alpha, \b} \, , \quad c_{k, i, \alpha, \b} \in \C  \, ,
$$
is {\sl absolutely convergent} if the function $ \mathbb{I}
  \ni (k,i,\a,\b) $ $\mapsto c_{k,i,\a,\b} \in \C $
is in $ L^1( \mathbb{I}, \mu) $ where $ \mu $ is the counting
measure of $ \mathbb{I}$. Then we set
$$
\sum_{(k, i, \alpha, \b )\in\mathbb{I}  } c_{k, i, \alpha, \b} := \int_{\mathbb{I}} c_{k,i,\a,\b} \, d \mu \, .
$$
\end{definition}
By the properties of the Lebesgue integral, given any sequence $ \{ I_l \}_{l\geq 0} $  of finite
subsets $ I_l \subset \mathbb{I} $ with $ I_l \subset I_{l+1} $ and  $ \cup_{l\geq 0} I_l =\mathbb{I} $,
the absolutely convergent series 
$$
\sum_{k, i, \alpha, \b } c_{k, i, \alpha, \b}:= \sum_{(k, i,
\alpha, \b )\in\mathbb{I}  } c_{k, i, \alpha, \b}
=\lim_{l\to\infty}\sum_{(k,i,\a,\b)\in I_l} c_{k,i,\a,\b}\,.
$$
\begin{definition} {\bf (Majorant-norm: scalar functions)}\label{majsc}
The majorant-norm of a formal  power series
\eqref{formalpowerA} is
\be\label{normadue}
\norma f \norma_{s,r} := \sup_{(y,z, \bar z) \in D(r)}
\sum_{k,i,\a, \b} |f_{k,i,\a, \b}| e^{|k|s} |y^i| |z^\a| |{\bar z}^{\b}|
 \ee
where  $ |k| := |k|_1:= |k_1| + \ldots + |k_n| $.
\end{definition}
By \eqref{formalpowerA} and  \eqref{normadue} we clearly have
$ \norma f \norma_{s,r} = \norma M f \norma_{s,r} $.

For every subset of indices $ I \subset \mathbb{I} $, we define
the projection
\begin{equation}\label{toblerone}
(\Pi_I f)(x,y,z, \bar z) := \sum_{ (k,i,\a, \b) \in I} f_{k,i,\a,
\b} e^{\ii k \cdot x} y^i z^\a {\bar z}^{\b}
\end{equation}
of the formal power series $ f $  in \eqref{formalpowerA}.
Clearly
\begin{equation}\label{proiezA}
\norma \Pi_I f \norma_{s,r} \leq \norma f \norma_{s,r}
\end{equation}
and, for any $ I, I' \subset \mathbb{I} $, it results
\begin{equation}\label{olinto}
\Pi_I \Pi_{I'}=\Pi_{I\cap I'}=\Pi_{I'}\Pi_I\,.
\end{equation}
Property \eqref{proiezA} is
one of the main advantages of the majorant-norm  
with respect to the
usual sup-norm   
\begin{equation}\label{laodicea}
|f|_{s,r}:=\sup_{v \in D(s,r)}|f(v)|\,.
\end{equation}
We now define useful projectors on the time Fourier indices. 
\begin{definition}
Given $\varsigma=(\varsigma_1,\ldots,\varsigma_n) \in \{ +,-\}^n$
we define
\begin{equation}\label{ficulle}
f_\varsigma := \Pi_\varsigma f := \Pi_{\Z^n_\varsigma \times  \N^n \times
\N^{(\Zi)} \times \N^{(\Zi)}} f = \sum_{k\in \Z^n_\varsigma
,i,\a,\b} f_{k, i, \alpha, \b} \, e^{\ii k \cdot x} y^i  z^\alpha
\bar z^{\b}
\end{equation}
where
\begin{equation}\label{romero}
\Z^n_\varsigma:= \Big\{ k \in \Z^n  \quad {\rm with} \quad
\begin{cases}
k_h \geq 0 \quad \mbox{\rm if} \quad \varsigma_h=+ \cr k_h < 0
\quad  \mbox{\rm if} \quad  \varsigma_h=-
\end{cases}
 \forall\,1\leq h\leq n\ \Big\}.
\end{equation}
\end{definition}
Then any formal series $ f $ can be decomposed as
\begin{equation}\label{falco}
f = \sum_{\varsigma \in \{ +,-\}^n}\Pi_\varsigma f
\end{equation}
and \eqref{proiezA} implies $\norma \Pi_\varsigma f \norma_{s,r} \leq \norma  f
\norma_{s,r}. $

\smallskip

We now investigate
the relations between formal power series
with finite majorant norm and analytic functions.
We recall that a function $ f : D(s,r) \to \C $ is
\begin{itemize}
\item
{\sc analytic}, if  $ f \in  C^1(D(s,r), \C ) $, namely
the 
Fr\'echet differential  $ D(s,r) \ni v \mapsto df (v) \in {\cal L}(E,\C) $ is continuous,
\item  {\sc weakly analytic},
if $ \forall v \in D(s,r) $, $ v' \in E \setminus \{ 0 \} $,  there exists $ \e > 0 $  such that  the function
\begin{equation*}
\{ \xi\in \C\,, \ |\xi|<\e\  \} \ \mapsto \ f(v + \xi v') \in \C
\end{equation*}
is analytic  in the usual sense of one complex variable.
\end{itemize}
A well known result (see e.g. Theorem 1, page 133 of \cite{PT}) states that a  function $ f $ is
\be\label{awalb}
{\rm analytic} \quad \Longleftrightarrow \quad   {\rm weakly \ analytic \ and \  locally \ bounded} \, .
\ee
\begin{lemma}\label{fundamA}
Suppose that the  formal power  series  \eqref{formalpowerA} is
absolutely convergent for all $ v \in D(s,r)$.
Then $ f(v) $ and $Mf(v)$, defined in
\eqref{formalpowerA} and \eqref{formalpower1A},
are well defined and
weakly analytic in $ D(s,r)$.

\noindent
If, moreover,  the sup-norm $|f|_{s,r}<\infty,$ resp. $|Mf|_{s,r}<\infty,$ then
$f$, resp. $Mf$, is analytic  in $ D(s,r)$.
\end{lemma}

\begin{pf}
Since the series \eqref{formalpowerA} is
absolutely convergent the functions $ f $, $ M f $,
and, for all $ \varsigma \in \{ + , - \}^n $,
$ f_\varsigma  :=\Pi_\varsigma f $, $M f_\varsigma $ (see \eqref{ficulle})
are well defined (also the series in \eqref{ficulle} is absolutely convergent).

We now prove that each $ M f_\varsigma$ is weakly analytic, namely
$ \forall v \in D(s,r) $, $ v' \in E \setminus \{ 0 \},$
\begin{equation}\label{tiberio}
M f_\varsigma (v + \xi v')= \sum_{ k\in \Z^n_\varsigma,i,\a,\b}
|f_{k, i , \a, \b}| {\mathfrak m}_{k,i,\a,\b}(v+\xi v')
\end{equation}
is analytic in  $ \{ |\xi|< \e \} $, for $ \e $ small enough (recall the notation \eqref{commodo}).
Since each $ \xi \mapsto {\mathfrak m}_{k,i,\a,\b} (v + \xi v') $ is  entire,
the analyticity of $ M f_\varsigma (v + \xi v') $ follows  once we prove that
the series \eqref{tiberio}  is totally convergent, namely
 \begin{equation}\label{lecalle}
\sum_{ k\in \Z^n_\varsigma,i,\a,\b} |f_{k, i , \a, \b}|
\sup_{|\xi|<\e} |{\mathfrak m}_{k,i,\a,\b}(v + \xi v') | < + \infty \, .
\end{equation}
Let us prove \eqref{lecalle}. We claim that, for $ \e $ small enough, there is
$ v^{\varsigma} \in D(s,r) $ such that
\be\label{marcoaurelio}
\sup_{|\xi|<\e} \big| {\mathfrak m}_{k,i,\a,\b} (v + \xi v') \big| \leq
{\mathfrak m}_{k,i,\a,\b} (v^\varsigma) \, , \quad \forall k\in
\Z^n_\varsigma,i,\a,\b \, .
\ee
Therefore  \eqref{lecalle} follows by
\begin{eqnarray*}
\sum_{ k\in \Z^n_\varsigma,i,\a,\b} |f_{k, i , \a, \b}|
\sup_{|\xi|<\e} |{\mathfrak m}_{k,i,\a,\b}(v + \xi v') | & \leq & \sum_{ k\in
\Z^n_\varsigma,i,\a,\b} |f_{k, i , \a, \b}|
{\mathfrak m}_{k,i,\a,\b}(v^\varsigma) \\
&= & M f_\varsigma (v^\varsigma) < + \infty \, .
\end{eqnarray*}
Let us construct $v^\varsigma\in D(s,r)$ satisfying
\eqref{marcoaurelio}. Since $v=(x,y,z,\bar z)\in D(s,r)$ we
have $x \in \T^n_s $ and, since $ \T^n_s $ is open, there is
$ 0 < s' < s $ such that $ |{\rm Im} (x_h)| < s' $, $ \forall 1\leq
h\leq n $. Hence, for $\e$ small enough,
\begin{equation}\label{massimo}
\sup_{|\xi|<\e} \big|  {\rm Im} (x + \xi x')_h \big|\leq s' <s \,,
\quad \forall\, 1\leq h\leq n\,.
\end{equation}
The vector $ v^\varsigma := ( x^\varsigma, y^\varsigma, z^\varsigma, \bar z^\varsigma)$
with  components
\begin{eqnarray}
&& x^\varsigma_h:= -\ii \varsigma_h s'\,,\qquad \qquad
y^\varsigma_h := |y_h|+\e |y'_h|\,,\qquad 1\leq h\leq n\,,
\nonumber
\\
&& z^\varsigma_h := |z_h|+\e |z'_h|\,, \qquad \bar z^\varsigma_h
:= |\bar z_h|+\e |\bar z'_h|\,,\qquad h\in \Z\, , \label{decimo}
\end{eqnarray}
belongs to  $ D(s,r) $ because $ |{\rm Im} \, x^\varsigma_h| = s'<s $,
$ \forall\, 1\leq h\leq n $, and also
   $ (y^\varsigma, z^\varsigma, \bar z^\varsigma) \in D(r) $ for $ \e $ small enough, because
$(y,z,\bar z)\in D(r) $ and $ D(r) $ is open.
Moreover,   $ \forall k\in \Z^n_\varsigma $, by \eqref{massimo},
\eqref{romero} and \eqref{decimo}, 
\be\label{pal} \sup_{|\xi|<\e} \big| e^{\ii k\cdot (x+\xi x')}
\big| \leq e^{|k|s'}  = e^{\ii k\cdot x^\varsigma} \, .
\ee
By  \eqref{commodo},    \eqref{decimo}, \eqref{pal},  we get \eqref{marcoaurelio}.
Hence each $ M f_\varsigma $ is weakly analytic and,
by the decomposition \eqref{falco}, also $ f $ and $ Mf $ are weakly analytic.
The final statement follows by \eqref{awalb}.
\end{pf}

\begin{corollary}\label{cor11}
 If $ \norma  f \norma_{s,r} < + \infty $ then
$ f $ and $ Mf $ are  analytic  and
\be\label{supfm}
|f|_{s,r}, |Mf|_{s,r} \ \leq\  \norma  f \norma_{s,r}\, .
\ee
\end{corollary}

\begin{pf}
For all $ v=( x, y, z, \bar z) \in \mathbb{T}^n_s \times D(r) $,
we have $ |e^{\ii k\cdot x}|\leq e^{|k|s} $ and
\begin{eqnarray*}
|f(v)|\,,\ |Mf(v)|
&\leq&
\sum_{k,i, \alpha, \b} |f_{k, i, \alpha, \b}| e^{|k| s} |y^i| |z^\alpha| |{\bar z}^{\b}|
\stackrel{\eqref{normadue}} \leq \norma  f \norma_{s,r}  < + \infty
\end{eqnarray*}
by assumption. Lemma \ref{fundamA} implies that $ f $, $ Mf $ are analytic.
\end{pf}

Now, we associate to any analytic function $ f:
D(s,r) \to \C $ the formal Taylor-Fourier power series
\be\label{hatf}
\mathtt f (v) := \sum_{(k,i,\a,\b) \in {\mathbb I}} f_{k,
i, \alpha, \b} \, e^{\ii k \cdot x} y^i  z^\alpha \bar z^{\b}
\ee
(as \eqref{formalpowerA}) with  Taylor-Fourier coefficients
\begin{equation}\label{subtuumpraesidium}
f_{k,i,\a, \b} := \frac{1}{(2\pi)^n} \int_{\mathbb{T}^n} e^{-\ii
k\cdot x}  \frac{1}{i!\a!\b!} (\partial_y^i \partial_z^\a
\partial_{\bar z}^\b f)(x,0,0,0) \, dx
\end{equation}
where $\partial_y^i \partial_z^\a \partial_{\bar z}^\b f$ are the
partial derivatives\footnote{For a multi-index $ \a = \sum_{1\leq j\leq k} e_{i_j}$, $ |\a | = k $, the partial derivative is
\begin{equation}\label{ulmobis}
\partial^{\a}_z f (x,y,z,\bar z) :=
 \frac{\partial^{k}}{ \partial \t_1 \ldots \partial \t_k}_{|\t=0}
 f ( x,y,z+\t_1 e_{i_1} + \ldots + \t_k e_{i_k},\bar z )
 \, .
\end{equation}}. 

\smallskip

What is the relation between $ f $ and its formal Taylor-Fourier series  $ \mathtt f $ ?

\begin{lemma}\label{mare}
Let $f:D(s,r)\to\C$ be analytic.
 If its associated
 Taylor-Fourier power series  \eqref{hatf}-
\eqref{subtuumpraesidium} is absolutely convergent
 in $D(s,r)$, and the sup-norm
\be\label{suphat}
\Big| \sum_{k,i,\a,\b} f_{k,i,\a, \b} \, e^{\ii k \cdot x}  y^i
z^\a {\bar z}^{\b}\Big|_{s,r}<\infty \, ,
\ee
then $ f = \mathtt f $,  $ \forall\, v\in D(s,r) $.
\end{lemma}
\begin{pf}
Since the Taylor-Fourier series   \eqref{hatf}-
\eqref{subtuumpraesidium} is absolutely convergent and \eqref{suphat} holds,
 by Lemma \ref{fundamA}
 the function $  \mathtt  f:D(s,r)\to\C $ is analytic.
The functions $ f =  \mathtt  f $ are equal
if  the Taylor-Fourier coefficients 
\begin{equation}\label{tabba}
  f_{k,i,\a, \b} =  \mathtt  f_{k,i,\a, \b}\,, \quad \forall \, k,i,\a, \b \, ,
\end{equation}
where the coefficients $  \mathtt  f_{k,i,\a, \b} $ are defined from $  \mathtt  f $ as in \eqref{subtuumpraesidium}.
Let us prove \eqref{tabba}. Indeed, for example,
\begin{eqnarray}
 \mathtt  f_{0,0,e_h,0}
&=&
\frac{1}{(2\pi)^n} \int_{\mathbb{T}^n}\frac{d}{d\xi}_{|\xi=0}
\sum_{k\in\Z^n,\, m\in \N}
f_{k,0,m e_h,0} e^{\ii k \cdot x} \xi^m \label{totcon} \\
&=&
\sum_{k\in\Z^n,\, m\in \N}
\frac{1}{(2\pi)^n} \int_{\mathbb{T}^n}\frac{d}{d\xi}_{|\xi=0}
f_{k,0,m e_h,0} e^{\ii k \cdot x} \xi^m
=f_{0,0,e_h,0}\,, \nonumber
\end{eqnarray}
using that the above series totally converge for $ r' < r $, namely
\begin{eqnarray*}
\sum_{k\in\Z^n,\, m\in \N}
\sup_{x\in \R,\, |\xi|\leq r'} |f_{k,0,m e_h,0} e^{\ii k \cdot x} \xi^m|
&\leq&
 \sum_{k\in\Z^n,\, m\in \N}
|f_{k,0,m e_h,0}| (r')^m
\\
&\leq&
\sum_{k,i,\a,\b} |f_{k,i,\a,\b}{\mathfrak m}_{k,i,\a,\b}(0,0,r' e_h,0)|
<\infty
\end{eqnarray*}
recall 
\eqref{commodo}.
For the others $ k, i, \a, \b $ in \eqref{tabba} is analogous.
\end{pf}

The above arguments also show the unicity of the Taylor-Fourier
expansion. 

\begin{lemma}\label{ttata}
If an analytic function $f:D(s,r)\to \C$
equals  an absolutely convergent formal series, i.e.
$ f(v) = \sum_{k,i,\a,\b} \tilde f_{k,i,\a, \b} e^{\ii k\cdot x} y^i z^\a \bar z^\b $,
then its Taylor-Fourier coefficients \eqref{subtuumpraesidium} are
$  f_{k,i,\a, \b}  = \tilde f_{k,i,\a, \b} $.
\end{lemma}

The majorant norm of $ f $
is equivalent to the sup-norm of its majorant $  Mf $.

\begin{lemma} \label{abscon}
\be\label{equiv12} |M f|_{s,r} \leq \norma f \norma_{s,r}  \leq
2^n |M f|_{s,r} \, . \ee
\end{lemma}
\begin{pf}
 The first inequality in \eqref{equiv12}
 is \eqref{supfm}.
The second one follows by
\begin{equation}\label{scipione}
\norma \Pi_\varsigma f \norma_{s,r} \leq |M f|_{s,r}\, , \quad
\forall\, \varsigma \in\{ +,-\}^n\,,
\end{equation}
where $\Pi_\varsigma f$ is defined in
 \eqref{ficulle}.
Let us prove \eqref{scipione}.
Let
$$
 D^+(r):= \Big\{  (y,z,\bar z)\in D(r)  \ :  \ y_h\geq 0\,, \, \forall\, 1\leq h\leq n\,,\  z_l,\bar z_l\geq 0\, , \forall\, l\in\Z \setminus
 {\cal I} \Big\}\, .
 $$
For any $ 0 \leq \s < s $, we have
\begin{eqnarray*}
 |M f|_{s,r} & =  & \sup_{(x,y,z, \bar z) \in D(s,r)}
 \Big| \sum_{k, i, \a, \b} |f_{k, i, \a, \b}| e^{\ii k \cdot x } y^i z^\a {\bar z}^{\b} \Big|  \\
 & \geq &
\sup_{x_1 = -\ii\varsigma_1 \s, \ldots, x_n = -\ii \varsigma_n \s,
 (y,z, \bar z) \in D^+(r)}
\Big| \sum_{k, i, \a, \b} |f_{k, i, \a, \b}| e^{\ii k \cdot x } y^i z^\a {\bar z}^{\b} \Big| \\
& \stackrel{\eqref{romero}}\geq & \sup_{ (y,z, \bar z) \in D^+(r)}
\sum_{k\in \Z^n_\varsigma, i, \a, \b} |f_{k, i, \a, \b}| e^{|k| \s} |y^i| |z^\a| |{\bar z}^{\b}|  \\
& = &  \sup_{ (y,z, \bar z) \in D(r)} \sum_{k\in \Z^n_\varsigma,
i, \a, \b} |f_{k, i, \a, \b}| e^{|k| \s} |y^i| |z^\a| |{\bar
z}^{\b}| =  \norma \Pi_\varsigma f \norma_{\s ,r} \, .
\end{eqnarray*}
Then \eqref{scipione} follows  since for every function
$g$ we have
$
\sup_{0\leq\s<s}\norma g\norma_{\s,r}=\norma g \norma_{s,r}\,.$
\end{pf}

\begin{definition}\label{orderscalarA} {\bf (Order relation: scalar  functions)}
Given  formal power series
$$
f = \sum_{k,i,\a,\b} f_{k,i,\a, \b} \, e^{\ii k \cdot x}  y^i
z^\a {\bar z}^{\b} \, , \ \ g = \sum_{k,i,\a,\b} g_{k,i,\a, \b} \,
e^{\ii k \cdot x} y^i z^\a {\bar z}^{\b} \, ,
$$
with $ g_{k,i, \a, \b} \in \R^+ $, we say that
\be\label{precedA}
f \prec  g \quad {\rm if} \quad | f_{k,i,\a,\b} | \leq
g_{k,i,\a,\b}\, ,\ \forall k,i,\a, \b \,.
\ee
\end{definition}
Note that,  by the definition  \eqref{formalpower1A} of majorant series,
\be\label{precemaj}
f \prec  g \quad \, \Longleftrightarrow \quad \,  f \prec  M f \prec g  \, .
\ee
Moreover, if $ \norma g \norma_{s,r} < + \infty $, then
$ f \prec  g $ $ \Longrightarrow $ $  \norma f \norma_{s,r} \leq
\norma g \norma_{s,r} $.

For any $\varsigma\in \{ +,-\}^n$
define $ q_\varsigma := (q_\varsigma^{(j)})_{j \in {\cal J}} $ as 
\begin{equation}\label{q}
q_\varsigma^{(j)} :=
\begin{cases}
-\varsigma_h \, \ii\qquad \mbox{if}\quad j=(1,h)\,, \ \ 1\leq
h\leq n\,, \cr 1\qquad\ \ \ \ \, \mbox{otherwise} \, .
\end{cases}
\end{equation}

\begin{lemma}\label{niennaA}
Assume  $ \norma f \norma_{s,r}, \norma g \norma_{s,r} < + \infty
$. Then
\begin{equation}\label{sumproA}
\quad f+ g\prec Mf+ Mg \, , \qquad f \cdot g\prec Mf \cdot Mg
\end{equation}
and \be\label{McompA} M \big(\partial_j (\Pi_\varsigma f)\big)  =
q_\varsigma^{(j)} \partial_j \big(M(\Pi_\varsigma f)\big)\,, \quad
j\in {\cal J} \,, \ee where $\partial_j$ is short for
$\partial_{v^{(j)}}$ and $ q_\varsigma^{(j)} $ are defined in \eqref{q}.
\end{lemma}
\begin{pf}
Since the series which define $ f $ and $ g $ are absolutely
convergent, the bounds \eqref{sumproA} follow by summing and
multiplying the series term by term. Next \eqref{McompA} follows by
differentiating the series term by term.
\end{pf}

An immediate consequence of  \eqref{sumproA} is
\be\label{fgrA}
\norma f+g\norma_{s,r} \leq \norma f\norma_{s,r} +
\norma g\norma_{s,r} \,, \quad \norma f \, g \norma_{s,r} \leq
\norma f \norma_{s,r} \norma g \norma_{s,r} \, .
\ee
The next lemma extends property \eqref{sumproA} for  infinite series.

\begin{lemma}\label{pollo}
Assume that  $ f^{(j)} $, $ g^{(j)}$  are formal power series satisfying
\begin{enumerate}
\item $f^{(j)}\prec g^{(j)}$,  $ \forall j\in {\cal J} $,
\item    $\norma g^{(j)}\norma_{s,r}<\infty$,  $ \forall j\in {\cal J} $,
\item $\sum_{j\in {\cal J}} |g^{(j)}(v)|<\infty$, $\forall\, v\in D(s,r) $, \label{num2}
\item $ g(v):=\sum_{j\in {\cal J}} g^{(j)}(v)$
is bounded in $D(s,r)$, namely $|g|_{s,r}<\infty.$
\end{enumerate}
Then the function $ g : D(s,r) \to \C $ is analytic, its Taylor-Fourier coefficients (defined
as in \eqref{subtuumpraesidium}) are
\begin{equation}\label{dabi}
g_{k,i,\a,\b} = \sum_{j\in {\cal J}} g^{(j)}_{k,i,\a,\b}\geq 0 \,,\qquad
\forall\, (k,i,\a,\b)\in\mathbb{I} \, ,
\end{equation}
and $\norma g\norma_{s,r}<\infty.$
Moreover
\begin{enumerate}
\item
$\sum_{j\in {\cal J}} |f^{(j)}(v)|<\infty$, $\forall\, v\in D(s,r),$
\item
$f(v):=\sum_{j\in {\cal J}} f^{(j)}(v)$
is analytic in $D(s,r)$,
\item  $ f \prec g $ and  $\|f \|_{s,r}\leq \|g\|_{s,r} < \infty $.
 \end{enumerate}
\end{lemma}

\begin{pf}
For each monomial $ {\mathfrak m}_{k,i,\a,\b} (v) $ (see  \eqref{commodo}) and
$ v = (x,y,z,\bar z)\in D(s,r) $, we have 
\be\label{mposit}
|{\mathfrak m}_{k,i,\a,\b}(v)|={\mathfrak m}_{k,i,\a,\b}(v_+)\,, 
\ee
where $v_+:=(\ii\,  {\rm Im}\, x, |y|,|z|,|\bar z|) \in D(s,r) $ with
$|y|:=(|y_1|,\ldots,|y_n|)$ and $|z|$, $|\bar z|$ are similarly defined.

Since $ \norma g^{(j)}\norma_{s,r}<\infty $ (and $ f^{(j)}\prec g^{(j)}$) the series
\be\label{espgj}
g^{(j)}(v) := \sum_{k,i,\a,\b} g^{(j)}_{k,i,\a,\b}  {\mathfrak m}_{k,i,\a,\b} (v)  \, , \qquad
g_{k,i,\a,\b}^{(j)} \geq 0
\ee
is absolutely convergent. 
For all  $ v \in D(s,r) $ we prove
that 
\begin{eqnarray}
\sum_{j\in{\cal J}}\sum_{k,i,\a,\b} |g^{(j)}_{k,i,\a,\b}  {\mathfrak m}_{k,i,\a,\b} (v)|
& \stackrel{\eqref{espgj}, \eqref{mposit}}  = &
\sum_{j\in{\cal J}}\sum_{k,i,\a,\b} g^{(j)}_{k,i,\a,\b}  {\mathfrak m}_{k,i,\a,\b} (v_+) \nonumber \\
& \stackrel{\eqref{espgj}} = & \sum_{j\in{\cal J}} g^{(j)}(v_+) = g(v_+)<\infty   \label{pela}
\end{eqnarray}
by assumption \ref{num2}. Therefore, by  Fubini's theorem, we exchange the order of the series
\be\label{agribus}
g(v)=
\sum_{j\in{\cal J}}\sum_{k,i,\a,\b} g^{(j)}_{k,i,\a,\b}  {\mathfrak m}_{k,i,\a,\b} (v)
=
\sum_{k,i,\a,\b} \Big(\sum_{j\in{\cal J}} g^{(j)}_{k,i,\a,\b}\Big)  {\mathfrak m}_{k,i,\a,\b} (v)
\ee
proving that  
$g$ is equal to an absolutely convergent  series.
Lemma \ref{fundamA} and the assumption 
$|g|_{s,r}<\infty$ imply that $ g $ is analytic in $D(s,r)$.
Moreover \eqref{agribus} and Lemma \ref{ttata} imply 
\eqref{dabi}. The $ g_{k,i,\a,\b} \geq 0 $ because $ g_{k,i,\a,\b}^{(j)} \geq 0 $, see \eqref{espgj}.
Therefore $ Mg = g $, and, by  \eqref{equiv12} and the assumption $|g|_{s,r}<\infty$, we
deduce  $ \|g\|_{s,r}<\infty.$

Concerning $ f $ we have
\begin{eqnarray*}
\sum_{j\in{\cal J}} |f^{(j)} (v)|
\leq
\sum_{j\in{\cal J}}\sum_{k,i,\a,\b} \Big| f^{(j)}_{k,i,\a,\b}  {\mathfrak m}_{k,i,\a,\b} (v)\Big|
\leq
\sum_{j\in{\cal J}}\sum_{k,i,\a,\b} g^{(j)}_{k,i,\a,\b}  |{\mathfrak m}_{k,i,\a,\b} (v)|
\stackrel{\eqref{pela}}
<\infty
\end{eqnarray*}
and, arguing as for $  g $, its Taylor-Fourier coefficients are
$ f_{k,i,\a,\b}=\sum_{j\in{\cal J}} f_{k,i,\a,\b}^{(j)} $, $ \forall (k,i,\a,\b)\in\mathbb{I} $.
Then
$$
|f_{k,i,\a,\b}|
\leq
\sum_{j\in{\cal J}} |f_{k,i,\a,\b}^{(j)}|
\leq \sum_{j\in{\cal J}} g_{k,i,\a,\b}^{(j)}
\stackrel{\eqref{dabi}} = g_{k,i,\a,\b}\,.
$$
Hence $f\prec g$ and $\|f \|_{s,r}\leq \|g\|_{s,r}<\infty$.
Finally $f$ is analytic by Lemma \ref{fundamA}.
\end{pf}

\begin{lemma}\label{baddi}  Let $ \norma f\norma_{s,r} < \infty $. Then,  $ \forall 0 < s' < s $, $ 0 < r' < r $,
we have $ \norma\partial_j f\norma_{s',r'} < \infty $.
\end{lemma}

\begin{pf}
It is enough to prove the lemma
for each $ f_\varsigma=\Pi_\varsigma f$
defined in \eqref{ficulle}.  By  $ \norma f\norma_{s,r}<\infty $ and Corollary \ref{cor11} 
the functions $ f_ \varsigma $, $ M f_ \varsigma $ are analytic and
$$
\norma\partial_j f_\varsigma\norma_{s',r'}
\stackrel{\eqref{equiv12}}\leq 2^n |M(\partial_j f_\varsigma)|_{s',r'}
\stackrel{\eqref{McompA}}=2^n
|\partial_j(M f_\varsigma)|_{s',r'}
\leq c|M f_\varsigma|_{s,r}
\stackrel{\eqref{equiv12}}\leq c\|f_\varsigma\|_{s,r}
$$
for a suitable $ c := c(n,s,s',r,r') $, having used  the Cauchy estimate (in one variable).
\end{pf}

We conclude this subsection with a simple result on representation of differentials.
\begin{lemma}\label{compod}
Let $ f : D(s,r)\to \C$ be Fr\'echet differentiable at  $ v_ 0$.
Then
\be\label{maritozzo1}
df(v_0)[v]=\sum_{j\in{\cal J}} \partial_j f(v_0)v^{(j)} \, , \quad \forall v = \sum_{j\in {\cal J}} v^{(j)} e_j \in E \, ,
\ee
and
\be\label{maritozzo}
\sum_{j\in{\cal J}} |\partial_j f(v_0)v^{(j)}|
\leq \|df (v_0)\|_{{\cal L}(E,\C)} \| v\|_E \, .
\ee
\end{lemma}
\begin{pf}
\eqref{maritozzo1} follows
by the continuity of the differential $df (v_0) \in {\cal L}(E, \C) $.
Next, consider a vector $\tilde v=( \tilde v^{(j)})_{j\in{\cal J}} \in E $ such that  $ |\tilde v_j|=|v_j| $ and
$$
 \tilde v^{(j)}  (\partial_j f)(v_0) = |(\partial_j f)(v_0) v^{(j)}| \, \, , \quad \forall\, j\in{\cal J} .
$$
Hence
$ df(v_0) [\tilde v] = $
$ \sum_{j\in{\cal J}}  \tilde v^{(j)} (\partial_j f)(v_0)  = $ $ \sum_{j\in{\cal J}}  |(\partial_j f)(v_0) v^{(j)}| $
which gives \eqref{maritozzo} because $\|\tilde v\|_E=\|v\|_E. $
\end{pf}




\subsubsection{Vector fields}

We now consider a  {\sl formal} vector field
\begin{equation}\label{nerva}
X(v) :=  \Big( X^{(j)}(v) \Big)_{j \in {\cal J}}
\end{equation}
where
each component $ X^{(j)} $
is a  formal  power series 
\begin{equation}\label{cesare}
X^{(j)}(v)=X^{(j)} (x,y,z, \bar z) =\sum_{k,i, \a,\b}
X^{(j)}_{k,i,\a,\b}
 \, e^{\ii k\cdot x} y^i z^\a \bar z^\b
\end{equation}
as in \eqref{formalpowerA}. We  define its ``majorant" vector field componentwise, namely
\begin{equation}\label{claudio}
MX(v) := \Big( (MX)^{(j)}(v) \Big)_{j \in {\cal J}}
:= \Big( M X^{(j)}(v) \Big)_{j \in {\cal J}} \, . 
\end{equation}
We  consider vector fields $X:D(s,r)\subset E \to E $, see \eqref{E}.
\begin{definition}
The vector field $X$ is absolutely convergent at $ v $ 
if every component $X^{(j)}(v) $, $ j \in  {\cal J} $, is absolutely convergent
(see  Definition \ref{augusto}) and
$$
\Big\| \big( X^{(j)}(v)\big)_{j\in {\cal J}} \Big\|_E < +
\infty\,.
$$
\end{definition}
The properties of the space $E$ in \eqref{E}  (as target space),   that we will use are:
\begin{enumerate}
\item $ E $ is a separable Hilbert space times a finite
dimensional space,
\item
the ``monotonicity property" of the norm
\begin{equation}\label{sgamirro}
 v_0 ,  v_1 \in E \ \ \ {\rm with} \ \ \  |v^{(j)}_0|\leq | v^{(j)}_1 |\, , \
\forall\, j\in {\cal J} \quad \Longrightarrow \quad \| v_0
\|_E\leq \| v_1 \|_E\,.
\end{equation}
\end{enumerate}
For $X:D(s,r)\to E$ we define the sup-norm
\begin{equation}\label{tiatira}
|X|_{s,r}:=\sup_{v \in D(s,r)}\|X(v) \|_{E,s,r}\,.
\end{equation}
\begin{definition} {\bf (Majorant-norm: vector field)}\label{MNV}
The majorant norm of a  formal vector field $ X $ as in \eqref{nerva} is
\begin{eqnarray}\label{normadueA}
\norma X \norma_{s,r} &:=&
  \sup_{(y,z, \bar z) \in D(r)}
\Big\| \Big(
 \sum_{k,i,\a, \b}  | X_{k,i,\a, \b}^{(j)}| e^{|k|s}
 |y^i| |z^\a| |{\bar z}^{\b}|
\Big)_{j\in {\cal J}} \Big\|_{E,s,r} \nonumber
\\
&=&  \sup_{(y,z, \bar z) \in D(r)} \Big\|
 \sum_{k,i,\a, \b}  | X_{k,i,\a, \b}| e^{|k|s}
 |y^i| |z^\a| |{\bar z}^{\b}| \Big\|_{E,s,r}
\end{eqnarray}
where
$$
  X_{k,i, \a, \b} :=
 \big(  X_{k,i, \a, \b}^{(j)}
 \big)_{j\in {\cal J}}\qquad
 {\rm and}
\qquad
 | X_{k,i, \a, \b}| :=
 \big(  |X_{k,i, \a, \b}^{(j)}|
 \big)_{j\in {\cal J}}\,.
$$
\end{definition}

\begin{remark}\label{differencePoschel}
The stronger norm 
(see \cite{Po2}) 
$$
\normastretta  X \normastretta_{s,r}  :=
\Big\| \Big( \sup_{(y,z, \bar z) \in D(r)}
 \sum_{k,i,\a, \b}   | X_{k,i,\a, \b}^{(j)}| e^{|k|s}
 |y^i| |z^\a| |{\bar z}^{\b}|
\Big)_{j\in {\cal J}} \Big\|_{E,s,r}
$$
is not suited for infinite dimensional systems: for $ X = Id $ we have $ \normastretta X \normastretta_{s,r} = + \infty $.
\end{remark}

\noindent
By   \eqref{normadueA} and  \eqref{claudio} we get
$ \norma X \norma_{s,r}=\norma M X \norma_{s,r} $.
For a subset of indices $ I \subset \mathbb{I} $ we define the projection
$$
(\Pi_I X)(x,y,z, \bar z) := \sum_{ (k,i,\a, \b) \in I} X_{k,i,\a,
\b} \, e^{\ii k \cdot x} y^i z^\a {\bar z}^{\b} \, .
$$
\begin{lemma} {\bf (Projection)}\label{lem:pro}
$ \forall I\subset \mathbb{I} $,
\begin{equation}\label{caligola}
 \norma \Pi_I X \norma_{s,r} \leq \norma X \norma_{s,r} \,.
\end{equation}
\end{lemma}

\begin{pf}
See \eqref{normadueA}.
\end{pf}

\begin{remark}\label{failsup}
The estimate \eqref{caligola} may fail for the sup-norm $ | \ |_{s,r} $ and suitable $ I $.
\end{remark}

Let us define  the  ``ultraviolet" reps. infrared projections 
\be\label{pisciati}
(\Pi_{|k| \geq K} X)(x,y,z, \bar z) := \sum_{ |k| \geq K,i,\a, \b}
X_{k,i,\a, \b} \, e^{\ii k \cdot x} y^i z^\a {\bar z}^{\b} \,,\quad \Pi_{|k| < K} := Id-\Pi_{|k| \geq K}  .
\ee
\begin{lemma}  {\bf (Smoothing)}  $ \forall\, 0 < s' < s $,
\begin{equation}\label{smoothl}
\norma\Pi_{|k| \geq K}X \norma_{s',r}\leq \frac{s}{s'} \, e^{-K(s-s')}\norma X
\norma_{s,r} \, .
\end{equation}
\end{lemma}
\begin{pf}
Recall
\eqref{normadueA} and use $ e^{|k| s'} \leq e^{|k| s} e^{- K (s-s')} $, $ \forall |k| \geq K $.
\end{pf}

We decompose  each formal vector field  
\begin{equation}\label{falco2}
X = \sum_{\varsigma \in \{ +,-\}^n}\Pi_\varsigma X
\end{equation}
applying \eqref{falco} componentwise
\begin{equation}\label{ficulle2}
 X_\varsigma := \Pi_\varsigma X :=
\Big( \Pi_\varsigma  X^{(j)} \Big)_{j \in {\cal J}}  
\end{equation}
recall \eqref{ficulle}. Clearly \eqref{caligola} implies
\be\label{ProXs}
\norma X_\varsigma \norma_{s,r} \leq \norma  X \norma_{s,r} \, .
\ee
In the next lemma we prove that,
if  $ X $ has finite majorant norm, then  it is analytic. 
\begin{lemma}\label{XMvectA}
Assume
\begin{equation}\label{forte}
\norma  X \norma_{s,r} < + \infty\,.
\end{equation}
Then the series in \eqref{nerva}-\eqref{cesare}, resp.
\eqref{claudio},   absolutely converge to the {\rm analytic}
vector field $X(v)$, resp. $MX(v)$, for every $v\in D(s,r).$
Moreover the sup-norm defined in \eqref{tiatira} satisfies
\be\label{XMXA}
|X|_{s,r},\ |MX|_{s,r} \ \leq\  \norma X \norma_{s,r}\,.
\ee
\end{lemma}
\begin{pf}
By \eqref{forte} and Definition \ref{MNV}, for each $ j \in {\cal
J} $,  we have
$$
\sup_{(y,z, \bar z) \in D(r)}  \sum_{k,i,\a,\b}
|X^{(j)}_{k,i,\a,\b}| e^{|k|s} |y^i| |z^\a| |\bar z^\b|  < +
\infty
$$
and Lemma \ref{fundamA} (and Corollary \ref{cor11}) implies that each coordinate function
$ X^{(j)} $, $ (MX)^{(j)} : D(s,r) \to \C $ is analytic.
Moreover \eqref{XMXA} follows applying  \eqref{supfm} componentwise.
By  \eqref{forte} the maps
$$
X\,,\,MX : D(s,r) \to E
$$
are bounded. Since $ E $ is a separable Hilbert space (times a
finite dimensional space), Theorem 3-Appendix A in \cite{PT},
implies that $ X$, $ MX : D(s,r) \to E $ are analytic.
\end{pf}

Viceversa, we associate to an analytic vector field $ X: D(s,r) \to E $
a formal Taylor-Fourier  vector field \eqref{nerva}-\eqref{cesare}
developing  each component $ X^{(j)} $ as in
 \eqref{hatf}-\eqref{subtuumpraesidium}.

\begin{definition}\label{ordvetA}
{\bf (Order relation: vector fields)} Given formal vector
fields $X$, $Y$, we say that
$$
X \prec Y
$$
if each coordinate  $  X^{(j)} \prec Y^{(j)}  $, $ j \in {\cal J}
$, according to Definition \ref{orderscalarA}.
\end{definition}

If $ \norma Y \norma_{s,r} < + \infty $ and
\be\label{orderXYA}
X \prec  Y \quad \Longrightarrow \quad  \norma X\norma_{s,r} \leq
\norma Y \norma_{s,r} \, .
\ee
Applying Lemma \ref{niennaA} component-wise we get
\begin{lemma}\label{niennaB}
If $\norma X \norma_{s,r}, \norma Y \norma_{s,r}<\infty$ then
$ X + Y \prec M X +  MY $ and  $ \norma X+Y\norma_{s,r} \leq \norma X\norma_{s,r} + \norma
Y\norma_{s,r} $.
\end{lemma}

\begin{lemma} \label{absconA}
\be\label{equiv12A} |M X|_{s,r} \leq \norma X \norma_{s,r}  \leq
2^n |M X|_{s,r} \, . \ee
\end{lemma}

\begin{pf}
As for Lemma \ref{abscon} with 
$f\rightsquigarrow X$, $|\sum_{k,i,\a,\b}|\rightsquigarrow
\|\sum_{k,i,\a,\b}\|_E$ and using \eqref{sgamirro}.
\end{pf}

We define the 
space of analytic  vector fields
$$
   {\cal V}_{s,r} :=
   {\cal V}_{s,r,E}:=
   \Big\{ X:  D(s,r) \to E
    \ \mbox{ with norm}\
     \norma X\norma_{s,r} < + \infty   \Big\}\,.
$$
By Lemma \ref{XMvectA} if $ X \in {\cal V}_{s,r} $ then $ X $ is analytic, namely the Fr\'echet differential
$ D(s,r) \ni v \mapsto dX (v) \in {\cal L}(E,E) $ is  continuous.
The next lemma bounds its operator norm from $   (E,s,r) := (E, \| \ \|_{E,s,r}) $ to $(E,s',r')$, see \eqref{normaEsr}.
\begin{lemma}\label{Cauchyx} {\bf (Cauchy estimate)}
Let $  X \in {\cal V}_{s,r} $. Then, for $ s/2 \leq s' < s $, $ r/2 \leq
r' < r $,
\begin{equation}\label{nitore}
\sup_{v \in D(s',r')}\| dX(v) \|_{{\cal L}((E,s,r),(E,s',r'))}
\leq 4\d^{-1}
|X|_{s,r}   
\end{equation}
where the sup-norm $ |X|_{s,r}   $ is defined in \eqref{tiatira} and 
\begin{equation}\label{diffusivumsui}
  \d :=  \min\Big\{ 1- \frac{s'}{s},   1-  \frac{r'}{r}  \Big\}\,.
\ee
\end{lemma}
\begin{pf}
In  the Appendix.
\end{pf}

The commutator of  two vector fields $ X, Y :D(s,r)\to E $ is 
\be\label{comm}
[X,Y](v):=d X(v) [Y(v)]- d Y(v)[X(v)]\,,\quad \forall\, v\in
D(s,r)\,.
\ee
The next lemma is the fundamental result of this section.
\begin{lemma}\label{settimiosevero} {\bf (Commutator)}
Let $ X, Y \in {\cal V}_{s,r} $. Then, for $ r/2 \leq r' < r $, $ s/2
\leq s' < s $,
\begin{equation}
\norma [ X, Y] \norma_{s',r'}
\leq
2^{2n+3} \delta^{-1} \norma
X\norma_{s,r} \norma Y\norma_{s,r}
\end{equation}
where $\d$ is defined in \eqref{diffusivumsui}.
\end{lemma}

\begin{pf}
The lemma follows by
\begin{equation}\label{barbarossa}
\norma d X [ Y ] \norma_{s',r'} \leq 4^{n+2} \delta^{-1} \norma
X\norma_{s,r} \norma Y\norma_{s,r} \, ,
\end{equation}
 the analogous estimate for $ d Y [ X ] $ and \eqref{comm}.

We claim that, for each  $\varsigma\in\{+,-\}^n $,  the vector
field $X_\varsigma $ 
defined in \eqref{ficulle2} satisfies
\begin{equation}\label{federico}
\norma d X_\varsigma [ Y ] \norma_{s',r'} \leq
 2^{n+2}\delta^{-1}
\norma X_\varsigma\norma_{s,r} \norma Y\norma_{s,r}
\end{equation}
which 
implies  \eqref{barbarossa} because
\begin{eqnarray*}
\norma d X [ Y ] \norma_{s',r'} & \stackrel{\eqref{falco2}} \leq&
\sum_{\varsigma \in \{ +,-\}^n} \norma d X_\varsigma [ Y ]
\norma_{s',r'} \stackrel{\eqref{federico}}\leq \sum_{\varsigma \in
\{ +,-\}^n} 2^{n+2} \delta^{-1} \norma  X_\varsigma\norma_{s,r} \norma
Y\norma_{s,r}
\\
&\stackrel{\eqref{ProXs}}\leq& \sum_{\varsigma \in \{ +,-\}^n}
2^{n+2} \delta^{-1}
 \norma X\norma_{s,r} \norma Y\norma_{s,r}
\leq 4^{n+2} \delta^{-1} \norma X\norma_{s,r} \norma Y\norma_{s,r}\, .
\end{eqnarray*}
Let us prove \eqref{federico}.  First note that,
since $   \| X_\varsigma \|_{s,r}  \stackrel{ \eqref{ProXs}} \leq  \| X \|_{s,r} < + \infty $ and
$  \| Y \|_{s,r} < + \infty $ by assumption, Lemma \ref{XMvectA} implies that the vector fields
\be\label{nellosp}
X_\varsigma, M X_\varsigma , Y, M Y : D(s,r) \to E \, , \quad \forall  \varsigma \in  \{  +, - \}^n  \, ,
\ee
are analytic, as well as each component
$  X_\varsigma^{(i)}, M X_\varsigma^{(i)}, Y^{(i)}, M Y^{(i)} : D(s,r) \to \C  $, $ i \in {\cal J} $.

The key for  proving the lemma is 
the following chain of inequalities: 
\begin{eqnarray}
 dX_\varsigma [Y]^{(i)}  \prec M ( dX_\varsigma [Y])^{(i)}  \!\!  \!\!  \!\!  \!\!
&\stackrel{\eqref{maritozzo1}}=& \!\!    \!
M \Big( \sum_{j\in{\cal J}}
(\partial_j X_\varsigma^{(i)})  Y^{(j)} \Big) \nonumber \\
& \stackrel{{\rm Lemma} \ \ref{pollo}}\prec &  \!\!  \!
 \sum_{j\in{\cal J}}  M (\partial_j X_\varsigma^{(i)})  M Y^{(j)}
    \label{quela}\\
&\stackrel{\eqref{McompA}} =&  \!\!  \!
  \sum_{j\in{\cal J}}  q_\varsigma^{(j)}
  \partial_j \big( M X_\varsigma^{(i)}\big)  M Y^{(j)} \nonumber
 \stackrel{\eqref{maritozzo1}}  =  d\big( M X_\varsigma^{(i)}\big)  \big[\tilde Y_q] \label{tuttequante}
\end{eqnarray}
where
\be\label{Yq}
\tilde Y_q := (\tilde Y_q^{(j)} )_{j \in
{\cal J}} := (q_\varsigma^{(j)} M Y^{(j)} )_{j \in {\cal J}} \in E \, .
\ee
Actually, since $ | q_\varsigma^{(j)} | = 1 $ (see \eqref{q}),  
 then 
\be\label{ugua}
\|  \tilde Y_q (v) \|_E = \|  M Y (v) \|_E  \stackrel{\eqref{nellosp}} < + \infty \, , \quad \forall v \in D(s,r) \, .
\ee
In \eqref{quela} above we applied Lemma \ref{pollo} with
\be\label{newsett}
s\rightsquigarrow s'\,,\
r\rightsquigarrow r'\,,\
f^{(j)} \rightsquigarrow (\partial_j X_\varsigma^{(i)})   Y^{(j)} \, , \ g^{(j)}
\rightsquigarrow M (\partial_j X_\varsigma^{(i)})  M Y^{(j)}\, .
\ee
Let us verify that the hypotheses of Lemma \ref{pollo} hold:
\begin{enumerate}
\item $f^{(j)}\prec g^{(j)}$ follows by \eqref{sumproA} and since $ \| f^{(j)} \|_{s',r'}$, $ \| g^{(j)} \|_{s',r'} < + \infty $
because $\| X^{(i)}_{\varsigma} \|_{s,r} \leq \| X \|_{s,r} < + \infty $,
$ \| Y^{(j)} \|_{s,r} \leq \| Y \|_{s,r}  < + \infty $, and  Lemma \ref{baddi}.
\item
$ \|g^{(j)}\|_{s',r'}<\infty $ is proved above.
\item
We have $ \sum_{j\in{\cal J}} |g^{(j)}(v)|<\infty$, for all $  v\in D(s',r')$, because
\begin{eqnarray*}
\sum_{j\in{\cal J}}  |g^{(j)}(v)|  \!\!  \!\!  & \stackrel{\eqref{newsett}} = & \!\!  \!\!  \sum_{j\in{\cal J}}
| M (\partial_j X_\varsigma^{(i)}) (v) M Y^{(j)}(v) |
\!\!  \stackrel{\eqref{McompA}}  = \!\!
\sum_{j\in{\cal J}}  |q_\varsigma^{(j)}
  \partial_j \big( M X_\varsigma^{(i)}\big)(v)  M Y^{(j)} (v)| \\
& \stackrel{\eqref{q}} = & \!\! \!\!  \sum_{j\in{\cal J}}  |  \partial_j \big( M X_\varsigma^{(i)}\big)(v)  M Y^{(j)} (v) |
 \stackrel{\eqref{maritozzo}}  \leq  \| d M X_{\varsigma}^{(i)}(v) \|_{{\cal L}(E,\C)} \| MY (v)\|_E < + \infty
\end{eqnarray*}
by \eqref{nellosp}, \eqref{ugua}. Actually we  also proved that
$ g^{(j)} = q_\varsigma^{(j)}  \partial_j \big( M X_\varsigma^{(i)}\big)  M Y^{(j)} $.
\item The function
$$
g (v) := \sum_{j\in{\cal J}} g^{(j)} (v) =
 \sum_{j\in{\cal J}} q_\varsigma^{(j)}  \partial_j \big( M X_\varsigma^{(i)}\big)  M Y^{(j)}
 \stackrel{\eqref{maritozzo1}} =
d\big( M X_\varsigma^{(i)}\big)  \big[\tilde Y_q]
$$
since $ M X_\varsigma^{(i)}$ is differentiable (see \eqref{nellosp}) and $ \tilde Y_q \in E $ (see \eqref{ugua}).

Moreover the bound 
$ | g |_{s',r'} < \infty $ follows by
$$
| g |_{s',r'} =  | d\big( M X_\varsigma^{(i)}\big)  \big[\tilde Y_q] |_{s',r'} \leq
|d\big( M X_\varsigma\big)  \big[\tilde Y_q]|_{s',r'}
$$
and
\begin{eqnarray}
|d\big( M X_\varsigma\big)  \big[\tilde Y_q]|_{s',r'}  \!\!  \!\!  \!\! & \stackrel{\eqref{tiatira}} = & \!\!  \!\!
\sup_{v \in D(s',r')} \Big\|   d\big( M X_\varsigma \big) (v)
\big[\tilde Y_q(v)\big]\Big\|_{E,s',r'} \nonumber \\
& \leq &  \!\!  \!\!\sup_{v \in D(s',r')} \Big\|   d\big( M X_\varsigma \big) (v) \Big\|_{{\cal L}((E,s,r),(E,s',r'))}
\| \tilde Y_q(v) \|_{E,s,r} \nonumber \\
&\stackrel{\eqref{nitore}}\leq&
 \!\!  \!\! 4 \d^{-1} | M X_\varsigma |_{s,r}
\sup_{v\in D(s',r')} \| \tilde Y_q(v)\|_{E,s,r}
\nonumber \\
& \stackrel{\eqref{XMXA}, \eqref{ugua}} \leq &
 \!\!  \!\! 4\d^{-1} \|X_\varsigma\|_{s,r} \sup_{v \in
D(s',r')} \|  \big(M Y\big)(v)\|_{E,s,r} \nonumber
\\
&  \stackrel{\eqref{tiatira}} \leq &
  \!\!  \!\! 4 \d^{-1} \|X_\varsigma\|_{s,r} | MY|_{s,r}
  \stackrel{\eqref{equiv12A}}\leq
4 \d^{-1} \norma X_\varsigma\norma_{s,r}
 \norma  Y\norma_{s,r} \label{ellala} < + \infty
\end{eqnarray}
because 
 $ \norma Y \norma_{s,r} < + \infty $ and $ \norma X_\varsigma\norma_{s,r} \leq \norma X \norma_{s,r} < + \infty $
by assumption.
\end{enumerate}
Hence Lemma \ref{pollo} implies  
$$
dX_\varsigma^{(i)} [Y] \stackrel{\eqref{maritozzo1}}= \sum_j (\partial_j X_\varsigma^{(i)})   Y^{(j)}
 =: f \stackrel{{\rm Lemma} \ \ref{pollo}} \prec g := d\big( M X_\varsigma^{(i)}\big)  \big[\tilde Y_q] \, , \quad \forall i \in {\cal J} \, ,
$$
namely, by \eqref{precemaj} and Definition \ref{ordvetA},
\begin{equation}\label{adriano}
dX_\varsigma [Y] \prec M ( dX_\varsigma [Y]) \prec  d\big( M X_\varsigma\big)
\big[\tilde Y_q] \, .
\end{equation}
Hence \eqref{tuttequante} is  fully justified.  By \eqref{adriano} and \eqref{orderXYA} we get
\begin{eqnarray}\label{quasifina}
\norma  dX_\varsigma [Y]\norma_{s',r'} \, \leq \, \norma d\big( M
X_\varsigma\big)  \big[\tilde Y_q]\norma_{s',r'}
& \stackrel{\eqref{equiv12A}}{\leq} &
2^n  \Big|M\Big( d\big( M X_\varsigma\big)  \big[\tilde Y_q]\Big) \Big|_{s',r'} \nonumber  \\
& = & 2^n  \Big|  d\big( M X_\varsigma\big)  \big[\tilde Y_q] \Big|_{s',r'}
\end{eqnarray}
because $d\big( M X_\varsigma\big)  \big[\tilde Y_q]$
coincides with its majorant by \eqref{adriano}.
Finally  \eqref{federico} follows by \eqref{quasifina}, \eqref{ellala}.
\end{pf}

\subsection{Hamiltonian formalism}

Given a function 
$ H : D(s,r)\subset E \to \C $
we define the associated Hamiltonian vector field
\be\label{HAMVF}
X_H := ( \partial_y H, - \partial_x H, -\ii \partial_{\bar z} H,
\ii \partial_{z} H)
\ee
where the partial derivatives are defined as  in \eqref{ulmobis}.

For a subset of indices $ I \subset  \mathbb{I}  $, the bound \eqref{caligola} implies
\be\label{proiezaa}
\norma X_{\Pi_I H}\norma_{s,r}\leq \norma X_H
\norma_{s,r} \,.
\ee
The Poisson brackets are defined by
\begin{eqnarray}\label{Poissonbraket}
\{H,K\} & := & \{H,K\}^{x,y} +   \{H,K\}^{z, \bar z} \nonumber \\
& := & \Big( \partial_x H \cdot \partial_y K -  \partial_x K \cdot
\partial_y H \Big) + \ii \Big( \partial_{z} H \cdot \partial_{\bar z} K -
\partial_{\bar z} H \cdot \partial_{z} K \Big) \nonumber \\
&  = &
\partial_x H \cdot \partial_y K -  \partial_x K \cdot \partial_y H +
\ii \partial_{z^+} H \cdot \partial_{z^-} K -
 \ii  \partial_{z^-} H \cdot \partial_{z^+} K \nonumber \\
& = & \partial_x H \cdot \partial_y K -  \partial_x K \cdot
\partial_y H + \ii\sum_{\s=\pm,\, j\in\Z \setminus {\cal I}} \s \partial_{z_j^\s} H
\, \partial_{z_j^{-\s}}K \,
\end{eqnarray}
where ``$\,\cdot\,$'' denotes the standard pairing
$a \cdot b:=\sum_j a_j b_j$. We recall the Jacobi identity
\be\label{Jacobi}
\{\{K,G\}, H\} + \{\{G,H\}, K\} +\{\{H,K\}, G\} = 0 \, .
\ee
Along this paper we shall use the Lie algebra notations
\be\label{expLie}
{\rm ad}_F := \{ \ , F \} \, , \quad e^{{\rm ad}_F} := \sum_{k=0}^\infty \frac{{\rm ad}_F^k}{k!} \, .
\ee
Given a set of indices 
\begin{equation}\label{cC}
{\cal I}:=\{\pluto_1,\ldots, \pluto_n \} \subset \Z \, ,
\end{equation}
we define the {\it momentum}
$$
 {\cal M} :=  {\cal M_{\cal I}} := \sum_{ l =1}^n   {\mathtt j}_l \, y_l + \sum_{j \in \Zi } j z_j {\bar
z}_j = \sum_{l =1}^n  {\mathtt j}_l \, y_l + \sum_{j \in \Zi } j z_j^+ z_j^-  \, .
$$
We say that a function $ H $ satisfies momentum conservation if $ \{ H, {\cal M} \} = 0 $.

By \eqref{Poissonbraket},
any monomial $e^{\ii k\cdot x} y^i z^\alpha {\bar z}^{\b} $
is an eigenvector of the operator ${\rm ad}_{\mathcal M}$, namely
\be\label{Mexp}
\{ e^{\ii k\cdot x} y^i z^\alpha{\bar z}^\beta , {\cal M} \} =  \pi(k,\alpha,\beta)
e^{\ii k\cdot x} y^i z^\alpha{\bar z}^\beta 
\ee
where
\be\label{momento2}
\pi(k,\alpha,\beta) :=
\sum_{l =1}^n \pluto_l k_l + \sum_{j\in \Zi} j(\alpha_j-\beta_j)  \, .
\ee
We refer to $\pi(k,\alpha,\beta)$ as the \textsl{momentum of the monomial}
$e^{\ii k\cdot x} y^i z^\alpha {\bar z}^{\b} $.
A monomial
satisfies momentum conservation if and only if  $\pi(k,\a,\b)=0. $ Moreover,
a power series \eqref{formalpowerA} with $ \| f \|_{s,r} < + \infty $  satisfies momentum conservation if
and only if all its monomials have zero momentum.

\medskip

Let $ \mathcal O \subset \R^n $ be a subset of {\it parameters}, and
\begin{equation}\label{sardi}
f: D(s,r) \times \mathcal O\to \C  \qquad \mbox{with}\qquad
X_f: D(s,r)\times \mathcal O\to E\,.
\end{equation}
 For   $ \l > 0 $, we consider
 \begin{eqnarray}\label{filadelfia}
 | {X}_f|^{\l}_{s, r, {\cal O}} :=
 | {X}_f|^{\l}_{s, r} & := &
 \sup_{{\cal O}} | {X}_f |_{s,r} + \l | {X}_f |^{{\rm lip}}_{s,r} \\
& := &
 \sup_{\xi \in {\cal O}} | {X}_f (\xi) |_{s,r} + \l \sup_{\xi,\eta\in {\cal O},\ \xi\neq \eta}
\frac{| {X}_f(\xi) - {X}_f(\eta)|_{s,r}}{|\xi-\eta|} \, .
\nonumber
\end{eqnarray}
Note that $ | \cdot|^{\l}_{s, r}$ is only a semi-norm on spaces of functions $ f $
because the Hamiltonian vector field $ X_f = 0 $ when $ f $ is  constant.

\begin{definition}\label{Hregular}
A function $f$ as in \eqref{sardi} is called
\begin{itemize}
\item
{\bf regular}, if the sup-norm $|X_f|_{s,r,{\mathcal O}}
:= \sup_{\mathcal O}|X_f|_{s,r}<\infty,$ see \eqref{tiatira}.
\item
{\bf M-regular}, if the majorant norm  $\norma  X_f\norma_{s,r,{\mathcal O}}:=
\sup_{\mathcal O}\norma X_f\norma_{s,r}<\infty,$  see \eqref{normadueA}.
\item
{\bf $\l $-regular},  if the Lipschitz semi-norm $|X_f|^{\l}_{s,r,{\mathcal O}}<\infty $, see \eqref{filadelfia}.
\end{itemize}
We denote by  $ \Hsr $ the space of M-regular Hamiltonians
and by $ \Hsrn $ 
its subspace of functions satisfying momentum conservation.

When  $ {\cal I} = \emptyset $ (namely there are no $(x,y)$-variables)
we denote the space of M-regular
functions simply by  $ {\cal H}_r $, similarly  $\Hrn $, and
we drop $ s $ form the norms, i.e.
$|\cdot|_r,$ $\|\cdot\|_r,|\cdot|_{r,{\cal O}}$, etc.
 \end{definition}
Note that, by \eqref{XMXA} and \eqref{filadelfia}, we have
\begin{equation}\label{efeso}
{\rm M-regular}\quad \Longrightarrow\quad {\rm regular} \quad
\Longleftarrow\quad \l-{\rm regular} \, .
\end{equation}
If $ H $, $ F $ satisfy momentum conservation, the same holds for
$ \{  H, K\} $. Indeed
by the Jacobi identity \eqref{Jacobi},
\be\label{MH0}
\{ {\cal M},H\} = 0 \ \ {\rm and} \ \ \{{\cal M}, K \} = 0 \quad
\Longrightarrow \quad  \{{\cal M}, \{H,K\} \} = 0 \, .
\ee
For $ H, K\in {\cal H}_{s,r} $ 
we have
\be\label{cippalippaA}
X_{\{H,K\}} = d X_H [X_K] - d  X_K[X_H] = [X_H, X_K]
\ee
and the commutator Lemma \ref{settimiosevero} implies the fundamental lemma below.  
\begin{lemma}\label{cauchy}
Let $ H, K \in {\mathcal H}_{s,r} $. Then, for all $ r/2 \leq r' < r $,
$ s/2 \leq s' < s $ \be\label{commXHK} \norma X_{\{H,K\}}\norma_{s',r'}
=  \norma [ X_H, X_K] \norma_{s',r'} \leq 2^{2n+3} \delta^{-1}
\norma X_H \norma_{s,r} \norma X_K\norma_{s,r} \ee where $\d$ is
defined in \eqref{diffusivumsui}.
\end{lemma}

Unlike the sup-norm, the majorant norm of a function is
very sensitive to coordinate transformations. For our purposes, we only need
to consider close to identity canonical transformations that are generated by an $ M$-regular Hamiltonian flow.
We show below that the $ M $-regular functions are closed under this group
and we estimate the majorant norm of the transformed Hamiltonian vector field.

\begin{lemma} {\bf (Hamiltonian flow)} \label{Hamflow}
Let $ r/2 \leq r' < r $, $ s/2 \leq s' < s $, and
$ F \in \mathcal{H}_{s,r}$ with
\be\label{defeta}
\norma X_F \norma_{s,r} <\eta:=\d/ (2^{2n + 5 }e)
\ee
with $\d$ defined in \eqref{diffusivumsui}.
Then  the time $ 1$-hamiltonian flow
$$
 \Phi^1_F  : D(s', r') \to D(s,r)
 $$
is well defined, analytic, symplectic, and,
$ \forall H \in \mathcal{H}_{s,r} $, we have $ H\circ \Phi^1_F   \in
\mathcal{H}_{s',r'}  $ and
\begin{equation}\label{iluvatar}
\norma  X_{H\circ \Phi^1_F}\norma_{s', r'}\leq \frac{\norma
X_H\norma_{s,r} }{ 1 - \eta^{-1} \norma  X_F \norma_{s,r}} \,.
\end{equation}
Finally if $ F, H \in \Hsr^{\rm null} $ then
$H\circ \Phi^1_F\in\mathcal{H}_{s',r'}^{\rm null}.$
\end{lemma}
\begin{pf}
We estimate by Lie series the Hamiltonian vector field of
 \be\label{Htras}
 H' = H \circ \Phi^1_F = e^{\rm ad_F} H = \sum_{k=0}^\infty \frac { {\rm ad}_F^k H}{k!} =
\sum_{k=0}^\infty \frac {  H^{(k)}}{k!}\,,
\quad {\rm i.e.} \  \,
X_{H'} = \sum_{k=0}^\infty \frac { X_{H^{(k)}}}{k!}\,,
\ee
where
$ H^{(i)} := {\rm ad}_F^i (H)= {\rm ad}_F ( H^{(i-1)}) $,  $ H^{(0)}:=H $.

For each $ k \geq 0 $, divide the intervals $[s',s]$ and $ [ r' ,
r ] $ into $ k $ equal segments and set
$$
s_i := s - i \, \frac{s-s'}{k} \, ,\qquad r_i := r - i \,
\frac{r-r'}{k} \, ,  \qquad  i = 0,\ldots,k \, .
$$
By \eqref{commXHK} we have
\be\label{lap} \norma  X_{H^{(i)}}
\norma_{s_i,r_i} = \norma  [ X_F, X_{H^{(i-1)}}] \norma_{s_i,r_i}
\leq  2^{2n+3}\delta_i^{-1}\norma
X_{H^{(i-1)}}\norma_{s_{i-1},r_{i-1}}\norma
X_F\norma_{s_{i-1},r_{i-1}} \ee where
\begin{eqnarray}\label{PP}
\d_i &:=&\min\left\{ 1-\frac{s_i}{s_{i-1}} \,,\,
1-\frac{r_i}{r_{i-1}} \right\} \,
\geq
\frac{\d}{k}\,.
\end{eqnarray}
By \eqref{lap}-\eqref{PP} we deduce
$$
\norma X_{H^{(i)}} \norma_{s_i,r_i}   \leq
 2^{2n+3}k \d^{-1} \norma X_{H^{(i-1)}}\norma_{s_{i-1},r_{i-1}}\norma X_F\norma_{s_{i-1},r_{i-1}}
  \, ,  \quad  i =  1, \ldots, k \, .
$$
Iterating $k$-times, and using $ \norma X_F\norma_{s_{i-1},r_{i-1}} \leq 4 \norma X_F \norma_{s,r} $ (see \eqref{empire})
\be\label{XHk}
\norma X_{H^{(k)}}\norma_{s',r'} \leq (2^{2n+5}k\d^{-1} )^k \norma
X_H\norma_{s,r}\norma X_F\norma_{s,r}^k \, .
\ee
By \eqref{Htras}, using $ k^k\leq e^k k! $ and recalling the definition of $ \eta $ in \eqref{defeta},
we estimate
\begin{eqnarray*}
\norma X_{H'}\norma_{s', r'}
& \stackrel{\eqref{Htras}} \leq&
 \sum_{k=0}^\infty \frac { \norma
X_{H^{(k)}}\norma_{s',r'}}{k!} \stackrel{\eqref{XHk}} \leq \norma
X_H\norma_{s,r} \sum_{k=0}^\infty \frac { (2^{2n+5}k\d^{-1}\norma
X_F\norma_{s,r} )^k}{k!}
\\
&\leq & \norma X_H\norma_{s,r} \sum_{k=0}^\infty (\eta^{-1}\norma
X_F\norma_{s,r})^k
\stackrel{\eqref{defeta}} = \frac{\norma X_H\norma_{s,r} }{1 - \eta^{-1}\norma  X_F\norma_{s,r}}
\end{eqnarray*}
proving \eqref{iluvatar}.

Finally, if $ F $ and $ H $ satisfy momentum conservation then each
$ {\rm ad}_F^k H $, $ k \geq 1 $, satisfy momentum conservation.
For $ k = 1 $ it is proved in \eqref{MH0} and, for $ k > 1 $, it follows by
induction and the Jacobi identity \eqref{Jacobi}.
By \eqref{Htras} we conclude that
also $ H \circ \Phi^1_F $  satisfies momentum conservation.
\end{pf}

We conclude this section with two simple lemmata.

\begin{lemma}\label{Plowerb}
Let
$ P = \sum_{|k| \leq K,i,\a, \b} P_{k,i, \a, \b} e^{\ii k \cdot x} y^i z^{\a} \bar
z^{\b}  $ and $ |\D_{k, i, \a, \b}| \geq
\g \langle k \rangle^{-\tau} $, $ \forall |k| \leq K ,i, \a, \b $.
Then
$$
F := \sum_{|k| \leq K,i,\a, \b} \frac{P_{k,i, \a, \b}}{\D_{k,i,\a,
\b}}
 e^{\ii k \cdot x} y^i z^{\a} \bar z^{\b} \quad \
{\rm satisfies} \ \quad
\norma  X_F \norma_{s,r} \leq  \g^{-1} K^\t  \norma X_P \norma_{s,r} \, .
$$
\end{lemma}

\begin{pf}
By  Definition \ref{MNV} and $ |\D_{k, i, \a, \b}| \geq \g
K^{-\tau} $ for all $ | k | \leq K $.
\end{pf}

\begin{lemma}\label{diago}
Let $ P =  \sum_{j\in \Z\setminus \mathcal I} P_j z_j\bar z_j $ with $  \|X_{P}\|_{r}  < \infty $.
Then
$ |P_j|\leq  \|X_{P}\|_{r} $.
\end{lemma}
\begin{pf}
By \eqref{HAMVF} and Definition  \ref{MNV} we have
$$
 \|X_{P}\|^2_{r}= 2 \sup_{\|z \|_{a,p}<r}\sum_{h\in \Z\setminus \mathcal I}
 |P_h|^2 \frac{|z_h|^2}{r^2} e^{2a|h|}\langle h\rangle ^{2p}\geq |P_j|^2
$$
by evaluating at  $z^{(j)}_h:=   \delta_{jh} e^{-a|j|}\langle j\rangle ^{p}r/\sqrt 2$.
\end{pf}

\section{Quasi-T\"oplitz functions}\label{sec:3}\setcounter{equation}{0}

Let  $ N_0 \in \N $, $ \teta, \mu \in \R $ be parameters  such that
\begin{equation}\label{caracalla}
\cc < \teta, \mu <  \CC \, , \quad  \dueCC N_0^{L-1}+2 \kappa
N_0^{b-1} < \cc\,, \quad \kappa:=\max_{1\leq l \leq n} |\pluto_l |
\, ,
\end{equation}
(the $ \pluto_l $ are defined in \eqref{cC})
where
\begin{equation}\label{figaro}
   0< b< L<1\,.
\end{equation}
For $ N \geq N_0 $, we decompose
\be\label{splitt}
\ell^{a,p}_{\cal I} \times \ell^{a,p}_{\cal I}  = \ell_L^{a,p} \oplus \ell_R^{a,p}  \oplus \ell_H^{a,p}
\ee
where
$$
\ell_L^{a,p} :=
\ell_L^{a,p}(N) :=
\Big\{  w = ( z^+, z^-) \in \ell^{a,p}_{\cal I} \times  \ell^{a,p}_{\cal I} \, : \, z_j^\s = 0 \, ,  \ \s = \pm \, , \ \forall |j| \geq \CC N^L  \Big\}
$$
$$
\ell_R^{a,p}  :=
\ell_R^{a,p} (N) :=
\Big\{  w = ( z^+, z^-) \in \ell^{a,p}_{\cal I}  \times  \ell^{a,p}_{\cal I} \, : \, z_j^\s = 0 \, , \ \s = \pm \, , \  {\rm unless} \ \CC N^L < |j| <
\ccN \Big\}
$$
$$
\ell_H^{a,p} :=
\ell_H^{a,p}(N) :=
\Big\{  w = ( z^+, z^-) \in \ell^{a,p}_{\cal I} \times \ell^{a,p}_{\cal I}  \, : \, z_j^\s = 0 \, , \ \s = \pm \, , \ \forall |j| \leq
\ccN  \Big\} \, .
$$
Note that by \eqref{caracalla}-\eqref{figaro} the subspaces $ \ell_L^{a,p}  \cap \ell_H^{a,p} = 0 $ and $ \ell_R^{a,p} \neq 0 $.
Accordingly we decompose any
$$
w \in \ell^{a,p} \times \ell^{a,p}  \qquad {\rm as} \qquad  w = w_L + w_R + w_H
$$
and we call $ w_L  \in \ell_L^{a,p} $ the  ``low momentum variables"
and $ w_H  \in \ell_H^{a,p} $ the ``high momentum variables".

We split the Poisson brackets in \eqref{Poissonbraket} as
$$
\{\cdot,\cdot\} = \{ \cdot, \cdot \}^{x,y} +
 \{\cdot,\cdot\}^{L} +\{\cdot,\cdot\}^{R}  + \{\cdot,\cdot\}^{H}
$$
where
\be\label{altapoi}
\{ H,K\}^H:=
\ii\sum_{\s=\pm,\, |j|>cN} \s \partial_{z_j^\s} H
\, \partial_{z_j^{-\s}}K\, .
\ee
The other Poisson brackets  $ \{\cdot,\cdot\}^{L} $, $ \{\cdot,\cdot\}^{R} $ are defined  analogously
with respect to the splitting \eqref{splitt}.

\begin{lemma}\label{bimby}
Consider two monomials
$ {\mathfrak m} =c_{k,i,\a,\b} e^{\ii k\cdot x} y^i z^\alpha {\bar z}^{\b}$
and $ {\mathfrak m}' = c'_{k',i',\a',\b'} e^{\ii k'\cdot x} y^{i'} z^{\alpha'} {\bar z}^{\b'}$.
The momentum of
${\mathfrak m}  {\mathfrak m}' $, $ \{{\mathfrak m},{\mathfrak m}'\} $,
$ \{{\mathfrak m},{\mathfrak m}'\}^{x,y} $,  $ \{{\mathfrak m},{\mathfrak m}'\}^L $,
$ \{{\mathfrak m},{\mathfrak m}'\}^R $,  $ \{{\mathfrak m},{\mathfrak m}'\}^H $,
equals  the sum of the  momenta of each monomial $ {\mathfrak m} $, $ {\mathfrak m}' $.
\end{lemma}

\begin{pf} By \eqref{momento2}, \eqref{Poissonbraket}, and
$$
\pi(k+k',\a+\a',\b+\b')=\pi(k,\a,\b)+\pi(k',\a',\b') =\pi(k,\a-e_j,\b)+\pi(k',\a',\b'-e_j)\, ,
$$
for any $ j \in \Z $.
\end{pf}

We now define subspaces of ${\cal H}_{s,r}$
(recall Definition \ref{Hregular}).

\begin{definition}\label{LM} {\bf (Low-momentum)}
A monomial
$
e^{\ii k\cdot x} y^i z^\alpha {\bar z}^{\b}
$
is $(N,\mu)$-low momentum if
\begin{equation}\label{zerobis}
\sum_{j \in\Z\setminus{\cal I}} |j| (\alpha_j + {\b}_j ) < \mu N^L \, ,  \quad |k| < N^b \, .
\end{equation}
We denote by %
$$
{\cal L}_{s,r}(N,\mu) \subset {\cal H}_{s,r}
$$
 the subspace of functions
\begin{equation}\label{g}
g = \sum g_{k,i,\a,\b}e^{\ii k\cdot x} y^i z^\a {\bar z}^{\b} \in  {\cal H}_{s,r}
\end{equation}
whose monomials are $(N,\mu)$-low momentum.
The corresponding projection
\be\label{proiLNmu}
\Pi^L_{N,\mu} : {\cal H}_{s,r} \to {\cal L}_{s,r}(N,\mu)
\ee
is defined as $ \Pi^L_{N,\mu} := \Pi_I $ (see \eqref{toblerone})
 where $ I $ is the subset of $ \mathbb{I}$ (see \eqref{I})
satisfying \eqref{zerobis}.
Finally, given $h\in\Z$, we denote by 
$$
{\cal L}_{s,r}(N,\mu,h) \subset {\cal L}_{s,r}(N,\mu)
$$
the subspace  of functions whose monomials satisfy
\be\label{vincula}
\pi(k,\alpha,\beta) + h = 0 \,.
\ee
\end{definition}
By \eqref{zerobis}, \eqref{caracalla}-\eqref{figaro},
any function in ${\cal L}_{s,r}(N,\mu)$, $ \cc < \mu < \CC $,  only depends on $ x, y, w_L $
and therefore
\begin{equation}\label{scamuffo}
g,g'\in {\cal L}_{s,r}(N,\mu) \ \
\Longrightarrow\ \
g g' ,\  \{g,g'\}^{x,y}, \ \{g,g'\}^L \  \ {\rm do\ not\ depend \ on \ } w_H\,.
\end{equation}
Moreover, by \eqref{momento2}, \eqref{caracalla}, \eqref{zerobis},  if
\begin{equation}\label{narsete}
|h| \geq \mu N ^L+\kappa N^b \quad
\Longrightarrow \quad {\cal L}_{s,r}(N,\mu,h)=\emptyset \, .
\end{equation}

\begin{definition}\label{BL}{\bf ($(N,\teta, \mu)$-bilinear)}
We denote by
$$
{\cal B}_{s,r}( N, \teta, \mu ) \subset  \Hsrn
$$
the subspace of the $ (N,\teta,\mu) $-bilinear functions defined as
\be\label{cincinnato}
f : =\sum_{|m|,|n|> \teta N, \s,\s'= \pm }
f^{\s,\s'}_{m,n} (x,y,w_L) z_m^\s z_n^{\s'}
\quad {\rm with} \quad
f^{\s,\s'}_{m,n} \in {\cal L}_{s,r}( N, \mu, \s m +\s' n)
\ee
and we denote the projection
$$
\Pi_{N,\teta,\mu} : \Hsr \to {\cal B}_{s,r}( N, \teta, \mu )\, .
$$
Explicitely, for $ g \in \Hsr $ as in \eqref{g}, the coefficients in \eqref{cincinnato} of
$f:=\Pi_{N,\teta,\mu} g$ are
\begin{equation}\label{evaristo}
f^{\s,\s'}_{m,n}(x,y,w^L):=
\sum_{(k,i,\a,\b)\ {\rm s.t.}\ \eqref{zerobis} \ {\rm holds}
\atop
{\rm and}\ \pi(k,\a,\b)=-\s m-\s' n}
f^{\s,\s'}_{k,i, \a,\b,m,n} e^{\ii k\cdot x} y^i z^\a \bar z^\b
\end{equation}
where\begin{eqnarray}
&&
 f^{+,+}_{k,i, \a,\b,m,n}:=
(2-\d_{mn})^{-1} g_{k,i, \a +e_m+e_n,\b}\,,\qquad
 f^{+,-}_{k,i, \a,\b,m,n}:= g_{ k,i, \a +e_m,\b+e_n}\,,
\nonumber
\\
&&
 f^{-,-}_{k,i,\a,\b,m,n}:=
(2-\d_{mn})^{-1} g_{k, i,\a,\b +e_m+e_n}\,,\qquad
 f^{-,+}_{k,i,\a,\b,m,n}:= g_{ k, i,\a +e_n,\b+e_m}\,.\quad
\label{carriego}
\end{eqnarray}
\end{definition}

For parameters $ \cc < \teta < \teta' $, $ \CC > \mu > \mu' $, we have
$$
{\cal B}_{s,r}(N,\teta',\mu')\subset {\cal B}_{s,r}(N,\teta,\mu) \, .
$$

\begin{remark}\label{devemangiare}
The projection $\Pi_{N,\teta,\mu}$ can be written in the form
$\Pi_I$, see \eqref{toblerone}, for a suitable $I\subset \mathbb{I}.$
The representation in \eqref{cincinnato} is not unique.
It becomes unique if we impose the ``symmetric'' conditions
\begin{equation}\label{domiziano}
 f^{\s,\s'}_{m,n}= f^{\s',\s}_{n,m}\,.
\end{equation}
Note that the coefficients in \eqref{evaristo}-\eqref{carriego} satisfy \eqref{domiziano}.
\end{remark}

 \subsection{T\"oplitz functions}

Let $ N \geq N_0 $.

\begin{definition}\label{matteo_aa} {\bf (T\"oplitz)}
A function $ f \in {\cal B}_{s,r}( N, \teta, \mu )  $ is  $(N,\teta,\mu)$-{\em T\"oplitz} if
the coefficients in \eqref{cincinnato} have the form
\begin{equation}\label{marco}
f^{\s,\s'}_{m,n}=
f^{\s,\s'}\big(\mathtt s(m),\s m+\s' n\big)
\quad
{\rm for\ some\ \ }  f^{\s,\s'}(\varsigma,h)\in{\cal L}_{s,r}(N,\mu,h)\,,
\end{equation}
with $\mathtt s(m) := {\rm sign}(m)$, $\varsigma=+,-$  and  $h\in\Z$. We denote by
$$
{\cal T}_{s,r}  := \Ta_{s,r} ( N, \teta , \mu ) \subset {\cal B}_{s,r}( N, \teta, \mu )
$$
the space  of  the $ (N,\teta ,\mu)$-{\it T\"oplitz  functions}.
\end{definition}

For parameters
$ N' \geq N $, $ \teta' \geq \teta $, $ \mu' \leq \mu $, $ r' \leq r $,
$ s' \leq s $ we have
\begin{equation}\label{calippo}
\Ta_{s,r} ( N, \teta, \mu )\subseteq
\Ta_{s',r'} ( N', \teta', \mu' )\,.
\end{equation}

\begin{lemma}\label{poisb}
Consider   $ f, g \in {\mathcal T}_{s,r}(N, \teta, \mu ) $ and  $p \in {\cal L}_{s,r}(N,\mu_1,0) $ with
$\cc<\mu,\mu_1<\CC$. For all
$0<s'<s\,, \, 0<r'<r$
and  $ \teta'  \geq \teta , \mu' \leq \mu $ one has
\be\label{pro11}
\Pi_{N,\teta',\mu'} \{ f , p \}^L \, ,
\ \Pi_{N,\teta',\mu'} \{ f , p \}^{x,y}
  \in {\mathcal T}_{s',r'}(N, \teta', \mu' ) \, .
\ee
If moreover
\begin{equation}\label{ugo}
\mu N^L+\kappa N^b <( \teta' - \teta) N
\end{equation}
then
\be\label{pr2}
\Pi_{N,\teta',\mu'} \{f,g \}^H \in {\mathcal T}_{s',r'}(N, \teta', \mu' ) \, .
\ee
\end{lemma}

\begin{pf}
Write $ f \in {\cal T}_{s,r}(N, \teta, \mu) $ as in \eqref{cincinnato} where
$f_{m,n}^{\s,\s'}$ satisfy \eqref{marco} and
\eqref{domiziano},
namely
\be\label{ricordo}
f_{m,n}^{\s,\s'}=f_{n,m}^{\s',\s}=f^{\sigma,\sigma'}(\mathtt s(m),\sigma m+\sigma' n)
\in{\cal L}_{s,r}(N,\mu, \s m+\s' n) \, ,
\ee
similarly for $ g $.
\\[1mm]
{\sc Proof of \eqref{pro11}.}
Since the variables $ z_m^\s $, $ z_n^{\s'} $, $|m|,|n|> \teta N $, are high momentum,
$$
\{f_{m,n}^{\s,\s'} z_m^\s z_n^{\s'}\,, \, p  \}^L
=
\{f_{m,n}^{\s,\s'}\,, \, p  \}^L \, z_m^\s z_n^{\s'}
$$
and $\{f_{m,n}^{\s,\s'}\,, \ p  \}^L$ does not depend on
 $ w_H $ by  \eqref{scamuffo} (recall that $ f_{m,n}^{\s,\s'} $, $ p \in {\cal L}_{s,r}(N, \mu )$).
 The coefficient of
 $z_m^\sigma z_n^{\sigma'}$ in $\Pi_{N,\teta',\mu'}\{f,p\}^{L}$ is
$$
\Pi^L_{N,\mu'}\{f_{m,n}^{\s,\s'}\,, \ p  \}^L
\stackrel{\eqref{ricordo}} = \Pi^L_{N,\mu'}\{f^{\s,\s'}(\mathtt s(m),\s m +\s' n)\,, \, p  \}^L
\in
{\cal L}_{s',r'}(N,\mu',\s m + \s' n)
$$
using Lemma \ref{bimby} (recall that $p$ has zero momentum).
The proof that
$\Pi_{N,\teta',\mu'}\{f,p\}^{x,y}
 \in \mathcal{T}_{s',r'}(N, \teta', \mu')$
is analogous.
\\[1mm]
{\sc Proof of \eqref{pr2}.}
A direct computation, using \eqref{altapoi}, gives
$$
\{ f , g\}^H=\sum_{|m|,|n|> \teta N,\, \s,\s'=\pm} p_{m,n}^{\s,\s'}z_m^\s z_n^{\s'}
$$
with
\begin{equation}\label{valentiniano}
p_{m,n}^{\s,\s'}=
2\ii \sum_{|l|> \teta N\,,\ \s_1=\pm}
\s_1 \Big(
f_{m,l}^{\s,\s_1} g_{l,n}^{-\s_1,\s'} +
f_{n,l}^{\s',\s_1} g_{l,m}^{-\s_1,\s}
\Big)\,.
\end{equation}
By \eqref{scamuffo} the coefficient $p_{m,n}^{\s,\s'}$ does not depend on $w_H.$
Therefore
\begin{equation}\label{giustiniano}
\Pi_{N, \teta',\mu'} \{ f , g\}^H=
\sum_{|m|,|n|> \teta' N,\, \s,\s'=\pm} q_{m,n}^{\s,\s'}z_m^\s z_n^{\s'}
\quad \mbox{with} \quad
q_{m,n}^{\s,\s'} := \Pi^L_{N,\mu'} p_{m,n}^{\s,\s'}
\end{equation}
(recall \eqref{proiLNmu}). It results
$ q_{m,n}^{\s,\s'}\in {\cal L}_{s',r'} (N,\mu',\s m +\s' n) $
by \eqref{giustiniano}, \eqref{valentiniano},  and Lemma \ref{bimby}
since, i.e.,
$$
f_{m,l}^{\s,\s_1}\in {\cal L}_{s,r} (N,\mu,\s m +\s_1 l)\quad
\mbox{and} \quad
g_{l,n}^{-\s_1,\s'}\in {\cal L}_{s,r} (N,\mu,-\s_1 l +\s' n)\,.
$$
Hence the $ (N, \teta',\mu') $-bilinear
function $ \Pi_{N, \teta',\mu'} \{ f , g\}^H$ in \eqref{giustiniano}
is  written in the form \eqref{cincinnato}.
It remains to prove that 
it is $ (N,\teta',\mu') $-T\"oplitz, namely that for all  $ |m|, |n|> \teta' N $, $ \s,\s' = \pm $,
\begin{equation}\label{valente}
q^{\s,\s'}_{m,n} =  q^{\s,\s'}\big(\mathtt s(m),\s m+\s' n\big)
\quad
{\rm for\ some } \quad  q^{\s,\s'}(\varsigma,h)\in{\cal L}_{s,r}(N,\mu',h)\,.
\end{equation}
Let us consider  in \eqref{valentiniano}-\eqref{giustiniano} the term (with  $m,n,\s,\s',\s_1$ fixed)
\be\label{sommal}
\Pi^L_{N,\mu'}
\sum_{|l|> \teta N}
f_{m,l}^{\s,\s_1} g_{l,n}^{-\s_1,\s'}
\ee
(the other is analogous).
Since $f,g\in {\cal T}_{s,r} (N,\teta,\mu)$ we have 
\be\label{oggi}
f_{m,l}^{\s,\s_1}=f^{\s,\s_1}\big(\mathtt s(m),\s m+\s_1 l \big) \in {\cal L}_{s,r}(N, \mu, \s m  + \s_1 l)
\ee
\begin{equation}\label{lalla}
g_{l,n}^{-\s_1,\s'}
=g^{-\s_1,\s'}\big(\mathtt s(l),-\s_1 l+\s' n \big) \in {\cal L}_{s,r}(N, \mu, -\s_1 l + \s' n) \,.
\end{equation}
By \eqref{narsete}, \eqref{oggi}, \eqref{lalla},
if the coefficients  $ f_{m,l}^{\s,\s_1} $, $ g_{l,n}^{-\s_1,\s'} $ are not zero then
\begin{equation}\label{susina}
|\s m+\s_1 l |\,,\   |-\s_1 l+\s' n | < \mu N^L+\kappa N^b\, .
\end{equation}
By \eqref{susina}, \eqref{caracalla}, we get
$ cN > |\s m+\s_1 l |= $ $ |\s \s_1 \mathtt s(m)|m|+\mathtt s(l)|l| |$,
which implies, since $|m|> \teta' N > \ccN $ (see \eqref{giustiniano}),  that the sign
\begin{equation}\label{laura}
    \mathtt s(l)=-\s \s_1 \mathtt s(m)\,.
\end{equation}
Moreover 
$$
|l|\geq |m|-|\s m+\s_1 l | \stackrel{\eqref{susina}}> \teta' N -\mu
N^L-\kappa N^b \stackrel{\eqref{ugo}}> \teta N\, .
$$
This shows that the restriction $|l|> \teta N$ in the sum \eqref{sommal} is automatically met.
Then
\begin{eqnarray*}
\Pi^L_{N,\mu'}
\sum_{|l|> \teta N}
f_{m,l}^{\s,\s_1} g_{l,n}^{-\s_1,\s'}
&\stackrel{\eqref{lalla}}=&
\Pi^L_{N,\mu'}
\sum_{l\in\Z }
f^{\s,\s_1}\big(\mathtt s(m),\s m+\s_1 l \big)
g^{-\s_1,\s'}\big(\mathtt s(l),-\s_1 l+\s' n \big)\\
&=&
\Pi^L_{N,\mu'}
\sum_{j\in\Z}
f^{\s,\s_1}\big(\mathtt s(m), j \big)
g^{-\s_1,\s'}\big(\mathtt s(l),\s m+\s' n-j \big)
\\
&\stackrel{\eqref{laura}}=&
\Pi^L_{N,\mu'}
\sum_{j\in\Z}
f^{\s,\s_1}\big(\mathtt s(m), j \big)
g^{-\s_1,\s'}\big(-\s\s_1 \mathtt s(m),\s m+\s' n-j \big)
\end{eqnarray*}
depends only on $\mathtt s(m)$ and $\s m +\s' n$, i.e. \eqref{valente}.
 \end{pf}

\subsection{Quasi-T\"oplitz functions}

Given  $ f \in \Hsr $ and   $ {\ft} \in {\mathcal T}_{s,r}(N, \teta, \mu ) $  we set
\begin{equation}\label{limi_aa}
\fr := N (\Pi_{N, \teta,\mu} f - \ft)  \, .
\end{equation}
All the functions $ f \in {\cal H}_{s,r} $ below possibly depend on parameters $  \xi \in {\cal O } $, see \eqref{sardi}.
For simplicity we shall often omit this dependence and  denote $ \| \ \|_{s,r, {\cal O}} =  \| \ \|_{s,r}  $.

\begin{definition} \label{topbis_aa} {\bf (Quasi-T\"oplitz)}
A function $ f \in \Hsrn $  is called $ (N_0,  \teta, \mu) $-quasi-T\"oplitz  if
the quasi-T\"oplitz semi-norm
\begin{equation}\label{unobisbis_aa}
\| f   \|_{s,r}^T := \| f   \|_{s,r, N_0,  \teta, \mu}^T:= \sup_{N \geq N_0}
\Big[ \inf_{ {\ft} \in \mathcal T_{s,r}(N,  \teta,\mu)} \Big( \max \{\| {X}_f\|_{s,r},
\| {X}_{\ft}\|_{s,r},\| {X}_{\fr}\|_{s, r} \} \Big) \Big]
\end{equation}
is finite. We define
$$
\mathcal Q^T_{s,r} := {\mathcal Q}^T_{s,r} ( N_0,  \teta, \mu) :=
\Big\{  f\in {\cal H}_{s, r}^{\rm null}  \, : \,  \| f \|_{s,r,N_0,  \teta, \mu}^T < \infty \Big\} \, .
$$
\end{definition}

In other words, a function $ f $ is
$ (N_0,  \teta, \mu) $-quasi-T\"oplitz
with semi-norm $ \| f \|_{s,r}^T$ if,
for all $ N \geq N_0 $,  $ \forall \e > 0 $, there is $ \ft \in \Ta_{s,r} (N,  \teta, \mu) $ such that
\be\label{defto}
\Pi_{N, \teta,\mu} f = \ft + N^{-1} {\hat f}   \quad {\rm and } \quad
\| {X}_f \|_{s,r} \, , \ \| {X}_{\tilde f}\|_{s,r} \, , \ \|{X}_{\hat f}\|_{s,r} \leq \| f \|_{s,r}^T  + \e \, .
\ee
We call $ \ft \in \Ta_{s,r} (N,  \teta, \mu) $  a
``{\it T\"oplitz approximation}" of $ f $ and $ \hat f $ the ``{\it T\"oplitz-defect}".
Note that, by Definition  \ref{matteo_aa} and  \eqref{limi_aa}
$$
 \Pi_{N, \teta,\mu} \ft = \ft \, , \quad \Pi_{N, \teta,\mu} {\hat f} = {\hat f} \, .
$$
By the definition \eqref{unobisbis_aa} we get
\be\label{maggio}
\| {X}_f\|_{s,r} \leq \| f   \|_{s,r}^T
\ee
and we  complete \eqref{efeso} noting that
\begin{equation}\label{efesobis}
\mbox{{\rm quasi-T\"oplitz}}\quad \Longrightarrow\quad
{\rm M-regular}\quad \Longrightarrow\quad
{\rm regular} \quad \Longleftarrow\quad
 \l-{\rm regular} \, .
\end{equation}
Clearly, if $ f $ is $ (N_0,  \teta,\mu) $-T\"oplitz then $ f $ is $ (N_0,  \teta,\mu) $-quasi-T\"oplitz and
\be\label{clearly}
\| f \|^T_{s, r, N_0, \teta,\mu} = \| {X}_f \|_{s,r} \, .
\ee
Then we have the following inclusions
$$
{\mathcal T}_{s,r}
\subset
{{\mathcal Q}^T_{s,r} \, , \, {\mathcal B}_{s,r}}
\subset
\Hsrn
\subset
\Hsr\,.
$$
Note that neither ${\mathcal B}_{s,r}\subseteq{\mathcal Q}^T_{s,r}$
nor ${\mathcal B}_{s,r}\supseteq{\mathcal Q}^T_{s,r}$.
\begin{lemma}\label{primo}
For parameters $ N_1 \geq N_0 $, $ \mu_1 \leq \mu $,
$  \teta_1 \geq  \teta $, $ r_1 \leq r $, $ s_1 \leq s $, we have
$$
{\mathcal Q}^T_{s,r} ( N_0,  \teta, \mu) \subset  {\mathcal Q}^T_{s_1,r_1} ( N_1,  \teta_1, \mu_1)
$$
and
\be\label{inscatola}
\| f \|_{s_1,r_1,N_1,  \teta_1, \mu_1}^T \leq
\max \{ s / s_1 , ( r / r_1 )^2 \} \| f \|_{s,r,N_0,  \teta, \mu}^T \, .
\ee
\end{lemma}

\begin{pf}
By \eqref{defto}, for all $ N \geq N_1 \geq N_0 $ (since $  \teta_ 1 \geq  \teta $, $ \mu_1 \leq \mu  $)
$$
\Pi_{N,  \teta_1,\mu_1}  f =
\Pi_{N,  \teta_1,\mu_1}  \Pi_{N, \teta,\mu}  f = \Pi_{N,  \teta_1,\mu_1}   \ft + N^{-1} \Pi_{N, \teta_1,\mu_1}  {\hat f} \, .
$$
The function $ \Pi_{N, \teta_1,\mu_1}  \ft \in \Ta_{s_1,r_1} (N,  \teta_1, \mu_1 ) $ and 
$$
\| {X}_{ \Pi_{N,  \teta_1,\mu_1}   \tilde f} \|_{s_1,r_1}
\stackrel{\eqref{proiezaa}}  \leq  \| {X}_{\tilde f}\|_{s_1,r_1}
\stackrel{\eqref{defto}} \leq \| f \|_{s_1,r_1}^T + \e  \, ,
$$
$$
\| {X}_{ \Pi_{N,  \teta_1,\mu_1}  \hat f} \|_{s_1,r_1}
\stackrel{\eqref{proiezaa}}  \leq  \|{X}_{\hat f}\|_{s_1,r_1} \stackrel{\eqref{defto}}  \leq \| f \|_{s_1,r_1}^T + \e \, .
$$
Hence, $ \forall N \geq N_1 $,
$$
 \inf_{ {\ft} \in \mathcal T_{s_1,r_1}(N,  \teta_1,\mu_1)} \Big( \max \{\| {X}_f\|_{s_1,r_1},
 \| {X}_{\ft}\|_{s_1,r_1},\| {X}_{\fr}\|_{s_1, r_1} \} \Big) \leq \| f \|_{s_1,r_1}^T + \e \, ,
$$
applying \eqref{empire} we have \eqref{inscatola}, because $ \e > 0 $ is arbitrary.
\end{pf}

For $ f \in \Hsr $ we define its homogeneous component
of degree $  l \in \N $,
\be\label{phom}
f^{(l)} := \Pi^{(l)} f :=  \sum_{k\in\Z^n\,,\, 2|i|+|\a|+|\b|=l}
 f_{k,i,\a,\b} \,
e^{\ii k\cdot x} y^i z^\a \bar z^\b \, ,
\ee
and the projections 
\be\label{ProK}
f_K := \Pi_{|k| \leq K} f :=
\sum_{|k| \leq K,i,\a, \b} f_{k,i,\a,\b} \,
e^{\ii k \cdot x} y^i z^\a {\bar z}^\b \,,\ \quad
\Pi_{> K} f:=f-\Pi_{|k| \leq K} f\,.
\ee
We also set
\be\label{Pleq2}
f^{\leq 2}_K := \Pi_{|k| \leq K} f^{\leq 2} \, , \quad
f^{\leq 2} := f^{(0)} +  f^{(1)} + f^{(2)} \, .
\ee
The above projectors $  \Pi^{(l)}  $, $ \Pi_{|k| \leq K} $, $ \Pi_{> K} $ have the form
$\Pi_I$, see \eqref{toblerone}, for  suitable subsets $ I \subset \mathbb{I}.$

\begin{lemma} {\bf (Projections)}\label{proJ}
Let $ f \in {\cal Q}^T_{s,r} (N_0, \teta,\mu) $. Then, for all $ l \in \N $,  $ K \in \N $,
\be\label{fhT}
\| \Pi^{(l)} f \|^{T}_{s,r, N_0,  \teta, \mu} \leq \| f \|^{T}_{s,r, N_0, \teta, \mu}
\ee
\be\label{fhle2T}
\| f^{\leq 2} \|^{T}_{s,r, N_0, \teta, \mu} \, ,  \ \| f - f^{\leq 2}_K \|^{T}_{s,r, N_0, \teta, \mu} \leq \| f \|^{T}_{s,r, N_0, \teta, \mu}
\ee
\be\label{fKT}
\| \Pi_{|k| \leq K} f \|^T_{s,r, N_0, \teta, \mu} \leq \| f \|^T_{s,r, N_0, \teta, \mu}
\ee
\be\label{bische}
\| \Pi_{k=0}\Pi_{|\a|=|\b|=1}\Pi^{(2)} f \|^T_{r, N_0, \teta, \mu} \leq \| \Pi^{(2)} f \|^T_{s,r, N_0, \teta, \mu}
\ee
and, $ \forall\, 0 < s' < s $,
\begin{equation}\label{smoothT}
\| \Pi_{> K} f  \|_{s',r, N_0, \teta,\mu}^T \leq
e^{-K(s-s')}\frac{s}{s'} \|  f \|_{s,r, N_0, \teta,\mu}^T \, .
\end{equation}
\end{lemma}

\begin{pf}
We first note that
by \eqref{olinto} (recall also Remark \ref{devemangiare}) we have
\begin{equation}\label{massiminoiltrace}
\Pi^{(l)} \, \Pi_{N, \teta,\mu} g =  \Pi_{N, \teta,\mu} \, \Pi^{(l)} g \, ,
\quad\forall \, g\in \Hsr\, .
\end{equation}
 Then, applying  $ \Pi^{(l)} $ in \eqref{defto}, we deduce that,
$ \forall N \geq N_0 $,  $ \forall \e > 0 $,  there is $ \ft \in \Ta_{s,r} (N,  \teta, \mu) $ such that
\be\label{deftop}
\Pi^{(l)} \Pi_{N, \teta,\mu} f = \Pi_{N, \teta,\mu} \Pi^{(l)}  f =
\Pi^{(l)}  \ft + N^{-1} \Pi^{(l)}  {\hat f}
\ee
and, by \eqref{proiezaa}, \eqref{defto},
\be\label{boundpT}
\|   X_{ \Pi^{(l)}  f} \|_{s,r} \, ,  
\ \| X_{ \Pi^{(l)}  \tilde f} \|_{s,r} \, ,  
\ \| X_{ \Pi^{(l)}  \hat f} \|_{s,r} 
\leq \| f \|_{s,r}^T  + \e \, .
\ee
We claim that  $ \Pi^{(l)}{\tilde f} \in \Ta_{s,r} (N,  \teta, \mu) $, $ \forall l \geq 0 $. Hence
\eqref{deftop}-\eqref{boundpT} imply
$  \Pi^{(l)}  f \in {\cal Q}^T_{s,r} (N_0, \teta,\mu) $ and
$$
\|  \Pi^{(l)}  f \|_{s,r}^T  \leq \| f \|_{s,r}^T  + \e \, ,
$$
i.e. \eqref{fhT}. Let us prove our claim.
For $ l = 0, 1 $ the projection $ \Pi^{(l)} \tilde f = 0 $
because $ \tilde f \in \Ta_{s,r} (N,  \teta, \mu) $ is bilinear. 
For  $ l \geq 2 $, write $ \tilde f $  in the form \eqref{cincinnato} with coefficients
$\tilde f^{\s,\s'}_{m,n}$ satisfying \eqref{marco}.
Then also $ g :=  \Pi^{(l)} \tilde f  $ has the form
\eqref{cincinnato} with coefficients
$$
g^{\s,\s'}_{m,n}=\Pi^{(l-2)}  \tilde f^{\s,\s'}_{m,n}  
$$
which satisfy  \eqref{marco} noting that
$ \Pi^{(l)} {\cal L}_{s,r}(N,\mu,h)
\subset  {\cal L}_{s,r}(N,\mu,h) $.
Hence  $ g\in \Ta_{s,r} (N,  \teta, \mu) $, $ \forall l \geq 0 $, proving the claim.
The proof of \eqref{fhle2T}, \eqref{fKT}, \eqref{bische},  and \eqref{smoothT} are similar (use also \eqref{smoothl}).
\end{pf}

\begin{lemma}\label{briscola}
Assume that, $ \forall N \geq N_* $, we have the decomposition
\be\label{decomNN}
G = G'_N + G''_N \quad {\rm  with} \quad
\|G'_N\|_{s,r,N, \teta,\mu}^T \leq K_1 \, , \ \ N\|{X}_{\Pi_{N, \teta,\mu} G''_N}\|_{s,r}\leq K_2 \, .
\ee
Then $\|G\|_{s,r,N_*, \teta,\mu}^T\leq \max \{ \| {X}_G \|_{s,r} , K_1 + K_2\}.$
\end{lemma}

\begin{pf}
By assumption, $ \forall N \geq N_* $, we have $ \|G'_N\|_{s,r,N, \teta,\mu}^T \leq K_1$.
Then, $ \forall \e>0 $, there exist
$\tilde G'_N \in {\cal T}_{s,r}(N, \teta,\mu)$, $ \hat G'_N $, such that
\be\label{pocobo}
\Pi_{N, \teta,\mu} G'_N=\tilde G'_N+ N^{-1} \hat G'_N \quad {\rm and} \quad
\|{X}_{\tilde G'_N}\|_{s,r}, \|{X}_{\hat G'_N}\|_{s,r} \leq K_1+\e \, .
\ee
Therefore, $ \forall N \geq N_* $,
$$
\Pi_{N, \teta,\mu} G =\tilde G_N+ N^{-1} \hat G_N\,,\quad
\tilde G_N:=\tilde G_N'\,,\ \
\hat G_N:=\hat G_N'+ N\Pi_{N, \teta,\mu} G''_N
$$
where $ \tilde G_N \in  {\cal T}_{s,r}(N, \teta,\mu) $ and
\be\label{la1t}
\| {X}_{\tilde G_N}\|_{s,r}=\| {X}_{\tilde G'_N}\|_{s,r} \stackrel{\eqref{pocobo}} \leq K_1+\e,
\ee
\be\label{la2t}
\| {X}_{\hat G_N}\|_{s,r}
\leq \| {X}_{\hat G'_N}\|_{s,r}+N \| {X}_{\Pi_{N, \teta ,\mu} G''_N}\|_{s,r}
\stackrel{\eqref{pocobo}, \eqref{decomNN}} \leq K_1+ \e + K_2\,.
\ee
Then $ G \in {\cal Q}^T_{s,r, N_*,  \teta,\mu} $ and
\begin{eqnarray*}
\| G\|_{s,r,N_*, \teta,\mu}^T & \leq &
\sup_{N\geq N_*} \max \big\{
\| {X}_G\|_{s,r}, \| {X}_{\tilde G_N}\|_{s,r},
\| {X}_{\hat G_N}\|_{s,r}
\big\} \\
& \stackrel{\eqref{la1t}, \eqref{la2t}} \leq &
\max \{ \| {X}_G \|_{s,r}, K_1+ K_2+ \e \}\,.
\end{eqnarray*}
Since $\e>0$ is arbitrary the lemma follows.
\end{pf}

The Poisson bracket of two quasi-T\"oplitz functions is quasi-T\"oplitz.

\begin{proposition}\label{festa} {\bf (Poisson bracket)}
Assume that $ f^{(1)}, f^{(2)} \in {\mathcal Q}^T_{s,r}
(N_0,  \teta,\mu ) $ and $ N_1 \geq N_0 $,
$ \mu_1 \leq \mu $, $  \teta_1 \geq  \teta $, $ s/2 \leq s_1 < s $,
$ r/2 \leq r_1 < r $
satisfy
\be\label{Leviatan}
\kappa N_1^{b-L}  < \mu - \mu_1, \
\mu N_1^{L-1} + \kappa N_1^{b-1} <  \teta_1  -  \teta,
 \ 2N_1 e^{- N_1^b \frac{s-s_1}{2}} < 1, \
 b(s-s_1)N_1^b>2\,.
\ee
Then
$$
\{ f^{(1)}, f^{(2)} \}\in {\mathcal Q}^T_{s_1, r_1} (N_1,  \teta_1, \mu_1)
$$
and
\begin{equation}\label{poisbound2}
 \| \{ f^{(1)}, f^{(2)} \} \|^T_{ s_1, r_1, N_1,  \teta_1, \mu_1}
 \leq  C(n) \delta^{-1}  \| f^{(1)} \|^T_{s, r, N_0,  \teta, \mu}\| f^{(2)} \|^T_{s, r, N_0,  \teta, \mu}
\end{equation}
where $ C(n) \geq 1 $ and 
\be\label{deltaciccio}
\d := 
\min \Big\{ 1 - \frac{s_1}{s}, 1 - \frac{r_1}{r}  \Big\}\,.
\ee
\end{proposition}

The proof  is based on the following splitting Lemma for the Poisson brackets.

\begin{lemma}\label{nn} {\bf (Splitting lemma)}
Let $ f^{(1)}, f^{(2)} \in {\mathcal Q}_{s,r}^T (N_0,  \teta,\mu ) $ and \eqref{Leviatan} hold.
Then, for all $ N \geq N_1 $,
\begin{eqnarray} \label{splitting}
& &
\Pi_{N, \teta_1,\mu_1}\{ f^{(1)}, f^{(2)} \}  = \nonumber \\
&  &
\Pi_{N, \teta_1,\mu_1} \Big(  \Big\{ \Pi_{N, \teta,\mu} f^{(1)} ,\Pi_{N, \teta,\mu} f^{(2)} \Big\}^H
+  \Big\{ \Pi _{N, \teta,\mu} f^{(1)}, \Pi^L_{N,2\mu} f^{(2)} \Big\}^L   +
 \Big\{  \Pi^L_{N,2\mu} f^{(1)}, \Pi _{N, \teta,\mu}  f^{(2)} \Big\}^L  \nonumber \\
 & & \quad \qquad \quad \quad +  \,
\Big\{ \Pi _{N, \teta,\mu} f^{(1)}, \Pi^L_{N,\mu} f^{(2)} \Big\}^{x,y}   +
\Big\{  \Pi^L_{N,\mu} f^{(1)}, \Pi _{N, \teta,\mu}  f^{(2)} \Big\}^{x,y} \nonumber \\
& & \quad \qquad \quad \quad +  \,
\Big\{  \Pi_{|k| \geq N^b} f^{(1)},  f^{(2)} \Big\} +
\Big\{  \Pi_{|k| < N^b} f^{(1)},   \Pi_{|k| \geq N^b} f^{(2)} \Big\} \Big) \, .
\end{eqnarray}
\end{lemma}

\begin{pf}
We have
\begin{eqnarray}
\{ f^{(1)}, f^{(2)} \} & = & \{ \Pi_{|k| < N^b} f^{(1)},  \Pi_{|k| < N^b} f^{(2)} \}  \label{pr1}\\
& + &
\{  \Pi_{|k| \geq N^b} f^{(1)},  f^{(2)} \} + \{  \Pi_{|k| < N^b} f^{(1)}, \Pi_{|k| \geq N^b} f^{(2)} \} \nonumber \, .
\end{eqnarray}
The last two terms correspond to the last line in
  \eqref{splitting}. We now study the first term in the right hand side
  of \eqref{pr1}.
We replace each $f^{(i)}$, $ i = 1, 2 $, with  single monomials (with zero momentum)
and we analyze under which conditions 
the projection
$$
\Pi_{N, \teta_1,\mu_1} \Big\{ e^{\ii k^{(1)} \cdot x} y^{i^{(1)}} z^{\alpha^{(1)}}\bar z^{\beta^{(1)}},
e^{\ii k^{(2)} \cdot x} y^{i^{(2)}} z^{\alpha^{(2)}}\bar z^{\beta^{(2)}} \Big\} \, , \quad |k^{(1)}|, |k^{(2)}| < N^b \, ,
$$
is not zero.
By direct inspection, recalling the Definition \ref{BL} of $ \Pi_{N, \theta_1, \mu_1 } $ and the expression \eqref{Poissonbraket} of the Poisson brackets
$ \{ \, , \, \} =  \{ \, , \, \}^{x,y} +  \{ \, , \, \}^{z, \bar z} $,
one of the following situations (apart from a trivial permutation of the indexes $1,2$) must hold:
 \begin{enumerate}
 \item
 one has $ z^{\alpha^{(1)}}\bar z^{\beta^{(1)}}= z^{\tilde\alpha^{(1)}}\bar z^{\tilde\beta^{(1)}} z_m^{\sigma}z_j^{\sigma_1}$
and $ z^{\alpha^{(2)}}\bar z^{\beta^{(2)}}= z^{\tilde\alpha^{(2)}}\bar z^{\tilde\beta^{(2)}} z_n^{\sigma'} z_j^{-\sigma_1}$ where
 $ |m|,|n| \geq  \teta_1 N $,  $ \s, \s_1, \s'  = \pm $,
  and   $z^{\tilde\alpha^{(1)}}\bar z^{\tilde\beta^{(1)}} z^{\tilde\alpha^{(2)}}\bar z^{\tilde\beta^{(2)}}$
  is of $(N,\mu_1)$-low momentum. We consider the Poisson bracket  $ \{ \, ,  \}^{z, \bar z} $
  (in the variables $  (z_j^+,z_j^-) $) of the monomials.

  \item one has $ z^{\alpha^{(1)}}\bar z^{\beta^{(1)}}= z^{\tilde\alpha^{(1)}}\bar z^{\tilde\beta^{(1)}} z_m^{\sigma}
  z_n^{\sigma'} z_j^{\sigma_1}$ and
  $ z^{\alpha^{(2)}}\bar z^{\beta^{(2)}}= z^{\tilde\alpha^{(2)}}\bar z^{\tilde\beta^{(2)}} z_j^{-\sigma_1}$ where
  $ |m|, |n| \geq  \teta_1 N $ and
  $ z^{\tilde\alpha^{(1)}}\bar z^{\tilde\beta^{(1)}} z^{\tilde\alpha^{(2)}}\bar z^{\tilde\beta^{(2)}}$ is of $(N,\mu_1)$--low momentum.
  We  consider the Poisson bracket  $\{ \, ,  \}^{z, \bar z} $. 
 \item
 one has $ z^{\alpha^{(1)}}\bar z^{\beta^{(1)}}= z^{\tilde\alpha^{(1)}}\bar z^{\tilde\beta^{(1)}} z_m^{\sigma} z_n^{\sigma'}$
 and  $ z^{\alpha^{(2)}}\bar z^{\beta^{(2)}}= z^{\tilde\alpha^{(2)}}\bar z^{\tilde\beta^{(2)}} $, where
 $ |m|,|n| \geq  \teta_1 N $
 and   $ z^{\tilde\alpha^{(1)}}\bar z^{\tilde\beta^{(1)} } z^{\tilde\alpha^{(2)}}\bar z^{\tilde\beta^{(2)}}$
 is of $(N,\mu_1)$-low momentum.
  We consider the Poisson bracket  $\{ \, , \,  \}^{x, y} $, i.e. in the variables $ ( x, y ) $.
 \end{enumerate}
Note that when we consider the $\{ \, , \,  \}^{x, y} $ Poisson bracket,  the  case
$$
z^{\alpha^{(1)}}\bar z^{\beta^{(1)}}= z^{\tilde\alpha^{(1)}}\bar z^{\tilde\beta^{(1)}}  z_m^{\sigma} \quad
{\rm and} \quad z^{\alpha^{(2)}}\bar z^{\beta^{(2)}}=
  z^{\tilde\alpha^{(2)}}\bar z^{\tilde\beta^{(2)}} z_n^{\sigma'}\, , \  \  |m|,|n| \geq  \teta_1 N \, ,
$$
 and   $ z^{\tilde\alpha^{(1)}}\bar z^{\tilde\beta^{(1)} }
  z^{\tilde\alpha^{(2)}}\bar z^{\tilde\beta^{(2)}}$
 is of $(N,\mu_1)$-low momentum, does not appear.
Indeed, the momentum conservation
$-\s m=\pi(\tilde \a^{(1)},\tilde \b^{(1)},k^{(1)})$,  \eqref{momento2} and $ |k^{(1)} | < N^b $, give
$$
\teta_1 N< |m|\leq
\sum_{l\in \Z \setminus {\cal I}} | l | (|\tilde \a^{(1)}_l|+|\tilde \b^{(1)}_l|)
+\kappa N^b
\leq \mu_1 N^L+\kappa N^b\,,
$$
which contradicts \eqref{caracalla}.
\\[1mm]
{\sc Case 1.} The momentum conservation of each monomial gives
\be\label{conm}
\sigma_1j = -\sigma m-\pi(\tilde\alpha^{(1)},\tilde\beta^{(1)}, k^{(1)}) =
\sigma' n + \pi(\tilde\alpha^{(2)},\tilde\beta^{(2)}, k^{(2)}) \, .
\ee
Since  $z^{\tilde\alpha^{(1)}}\bar z^{\tilde\beta^{(1)}} z^{\tilde\alpha^{(2)}}\bar z^{\tilde\beta^{(2)}}$
  is of $(N,\mu_1)$-low momentum (Definition \ref{LM}),
$$
\sum_{l\in \Z \setminus {\cal I}  }|l|(\tilde\alpha^{(1)}_l+\tilde\beta^{(1)}_l+\tilde\alpha^{(2)}_l+
\tilde\beta^{(2)}_l)\leq \mu_1 N^L \ \
\Longrightarrow  \ \ \sum_{l\in \Z \setminus {\cal I}  }|l|(\tilde\alpha^{(i)}_l+\tilde\beta^{(i)}_l)\leq \mu_1 N^L \, , \ i = 1,2 \, ,
$$
which implies, by \eqref{conm}, \eqref{momento2}, $ |k^{(1)} | < N^b $,
$ |j| \geq  \teta_1 N - \mu_1 N^L - \kappa N^b  >   \teta N $ by  \eqref{Leviatan}.
Hence $|m|,|n|,|j| >  \teta N $.
Then $ e^{\ii k^{(h)} \cdot x } y^{i^{(h)}} z^{\alpha^{(h)}}
\bar z^{\beta^{(h)}} $, $ h = 1,2 $,
are  $(N, \teta,\mu)$-bilinear. Moreover the $(z_j, {\bar z}_j ) $ are high momentum variables, namely
$ \{  \, , \, \}^{z, \bar z} = \{ \, , \, \}^H $, see \eqref{altapoi}.
As $m,n$ run over all  $\Z \setminus {\cal I} $ with $ |m|,|n| \geq  \teta_1 N$,
 we obtain the first term in formula \eqref{splitting}.
\\[1mm]
{\sc Case 2.}  The momentum conservation of the second monomial reads
\be\label{como2}
  - \sigma_1 j = - \pi(\tilde\alpha^{(2)},\tilde\beta^{(2)}, k^{(2)}) \, .
\ee
  Then, using also  \eqref{momento2}, $ |k^{(2)} | < N^b $, that
  $ z^{\tilde\alpha^{(1)}} \bar z^{\tilde\beta^{(1)}} z^{\tilde\alpha^{(2)}}\bar z^{\tilde\beta^{(2)}}$
  is of $ ( N , \mu_1 ) $-low momentum,
  $$
  |j|+ \sum_{l \in \Z \setminus {\cal I} } |l|(\tilde\alpha^{(1)}_l+\tilde\beta^{(1)}_l) 
  \stackrel{\eqref{como2}} =
  |\pi(\tilde\alpha^{(2)},\tilde\beta^{(2)}, k^{(2)})|
  + \sum_{l\in \Z  \setminus {\cal I} } |l|(\tilde\alpha^{(1)}_l+\tilde\beta^{(1)}_l) \leq
$$
$$
\sum_{ l\in \Z \setminus {\cal I} } |l|(\tilde\alpha^{(1)}_l+\tilde\beta^{(1)}_l+\tilde\alpha^{(2)}_l+\tilde\beta^{(2)}_l )
 + \kappa N^b  \leq  \mu_1 N^L  + \kappa N^b \stackrel{\eqref{Leviatan}} < \mu N^L \, .
$$
Then  $ z^{\tilde\alpha^{(1)}}\bar z^{\tilde\beta^{(1)}} z_j^{\s_1} $
  is of $(N,\mu_1)$-low momentum and the first monomial
$$
e^{\ii k^{(1)} \cdot x } y^{i^{(1)}} z^{\alpha^{(1)}}\bar z^{\beta^{(1)}} =  e^{\ii k^{(1)} \cdot x } y^{i^{(1)}}
z^{\tilde\alpha^{(1)}}\bar z^{\tilde\beta^{(1)}} z_j^{\s_1} z_m^\s z_n^{\s'}
$$
is $(N, \teta,\mu)$-bilinear ($ \mu_1 \leq \mu $). The second monomial
$$
e^{\ii k^{(2)} \cdot x } y^{i^{(2)}} z^{\alpha^{(2)}}\bar z^{\beta^{(2)}}
= e^{\ii k^{(2)} \cdot x } y^{i^{(2)}}  z^{\tilde\alpha^{(2)}}\bar z^{\tilde\beta^{(2)}} z_j^{-\s_1}
$$
is $ ( N , 2 \mu ) $-low-momentum because,  arguing as above, 
\begin{eqnarray*}
  |j|+ \sum_{l } |l|(\tilde\alpha^{(2)}_l+\tilde\beta^{(2)}_l)
& \stackrel{ \eqref{como2}}   = &
 |\pi(\tilde\alpha^{(2)},\tilde\beta^{(2)}, k^{(2)})| + \sum_{l} |l|(\tilde\alpha^{(2)}_l+\tilde\beta^{(2)}_l) \\
& \leq & 2 \mu_1 N^L + \kappa N^b \stackrel{\eqref{Leviatan}} < 2 \mu N^L \, .
\end{eqnarray*}
The $ (z_j, {\bar z}_j ) $ are  low momentum variables,  
namely $  \{  \, , \, \}^{z, \bar z} = \{ \, , \, \}^L  $, and
we obtain the second and third  contribution in formula \eqref{splitting}.
\\[1mm]
{\sc Case 3.} We have, for $ i = 1, 2 $, that
$$
\sum_{l} |l| ( {\tilde \a}_l^{(i)} + {\tilde \b}_l^{(i)} ) \leq
\sum_{l} |l| ( {\tilde \a}_l^{(1)} + {\tilde \b}_l^{(1)} + {\tilde \a}_l^{(2)} + {\tilde \b}_l^{(2)} ) \leq \mu_1 N^L \leq \mu N^L \ .
$$
Then $ e^{\ii k^{(1)} \cdot x } y^{i^{(1)}} z^{\alpha^{(1)}}\bar z^{\beta^{(1)}} $
is $(N, \teta,\mu)$-bilinear and $ e^{\ii k^{(2)} \cdot x } y^{i^{(2)}} z^{\alpha^{(2)}}\bar z^{\beta^{(2)}} $
is $ ( N , \mu ) $-low-momentum.
We obtain the fourth and fifth
 contribution in formula \eqref{splitting}.
\end{pf}

\medskip

\noindent
{\sc Proof of Proposition \ref{festa}.}
Since $ f^{(i)} \in {\mathcal Q}^T_{s,r} (N_0,  \teta, \mu) $, $  i = 1, 2 $,
for all $ N \geq N_1 \geq N_0  $ there exist  $ {\tilde f}^{(i)} \in \Ta_{s,r} (N,  \teta, \mu)$ and  $ {\hat f}^{(i)}$
such that (see \eqref{defto})
\be\label{f1f2pi}
\Pi_{N,  \teta, \mu} f^{(i)} = {\tilde f}^{(i)} + N^{-1} {\hat f}^{(i)} \, , \quad i = 1, 2 \, ,
\ee
and
\be\label{b12}
\| {X}_{f^{(i)}} \|_{s,r}, \
\| {X}_{{\tilde f}^{(i)}} \|_{s,r}, \  \| {X}_{{\hat f}^{(i)}} \|_{s,r} \leq 2 \| f^{(i)} \|^T_{s,r} \, .
\ee
In order to show that
$ \{f^{(1)},f^{(2)}\} \in {\mathcal Q}^T_{s_1,r_1} (N_1,  \teta_1, \mu_1 ) $ and prove
\eqref{poisbound2} we have to provide a decomposition
$$
\Pi_{N,  \teta_1,\mu_1} \{ f^{(1)}, f^{(2)} \} =  \ft^{(1,2)}+ N^{-1} \fr^{(1,2)} \, , \quad  \forall N \geq N_1 \, ,
$$
so that $ {\ft}^{(1,2)}\in \mathcal{T}_{s_1, r_1} (N,  \teta_1, \mu_1 ) $ and
\begin{equation}\label{stizza}
\| {X}_{   \{ f^{(1)}, f^{(2)} \}   }   \|_{s_1, r_1}, \
\| {X}_{\ft^{(1,2)}}\|_{s_1, r_1}, \| {X}_{\fr^{(1,2)}}\|_{s_1, r_1} <
C (n)\delta^{-1} \| f^{(1)} \|_{s,r}^T \|  f^{(2)} \|_{s,r}^T
\end{equation}
(for brevity  we omit the indices $ N_1, \theta_1, \mu_1, N_0, \theta, \mu $).
By \eqref{commXHK} we have  ($ \d $ is defined in \eqref{deltaciccio})
$$
\|  {X}_{ \{ f^{(1)}, f^{(2)} \} } \|_{s_1,r_1}
\leq
 2^{2n+3} \d^{-1} \|  {X}_{f^{(1)}} \|_{s,r} \|  {X}_{f^{(2)}} \|_{s,r} \, .
$$
Considering \eqref{f1f2pi} and \eqref{splitting}, we define
 the candidate T\"oplitz approximation
\begin{eqnarray} \label{split}
\ft^{(1,2)} & := &
\Pi_{N, \teta_1,\mu_1} \Big(  \Big\{ {\tilde f}^{(1)}, {\tilde  f}^{(2)} \Big\}^H
+  \Big\{ {\tilde f}^{(1)}, \Pi^L_{N,2\mu} f^{(2)} \Big\}^L   +
 \Big\{  \Pi^L_{N,2\mu} f^{(1)},  {\tilde  f}^{(2)} \Big\}^L  \nonumber \\
 & & \quad  \quad \quad +  \,
\Big\{ {\tilde f}^{(1)}, \Pi^L_{N,\mu} f^{(2)} \Big\}^{x,y}   +
\Big\{  \Pi^L_{N,\mu} f^{(1)},  {\tilde f}^{(2)} \Big\}^{x,y}
\end{eqnarray}
and T\"oplitz-defect
\begin{eqnarray}
\fr^{(1,2)}
&:=&
N \Big(\Pi_{N, \teta_1,\mu_1}  \{f^{(1)},f^{(2)} \} - \ft^{(1,2)} \Big)
\, .
\end{eqnarray}
Lemma \ref{poisb} and \eqref{Leviatan} imply that $ \ft^{(1,2)}
\in \mathcal{T}_{s_1,r_1} (N,  \teta_1, \mu_1) $.
The estimate \eqref{stizza} for $\ft^{(1,2)}$
follows by \eqref{split}, \eqref{commXHK},
\eqref{proiezaa},
\eqref{b12}.
Next
\begin{eqnarray*}
\fr^{(1,2)} & = &
\Pi_{N, \teta_1,\mu_1} \Big(  \Big\{ {\tilde f}^{(1)} , {\hat f}^{(2)} \Big\}^H +
 \Big\{ {\hat f}^{(1)} , {\tilde f}^{(2)} \Big\}^H +  N^{-1} \Big\{ {\hat f}^{(1)} , {\hat f}^{(2)} \Big\}^H   \\
 & + &  \Big\{  {\hat f}^{(1)}, \Pi^L_{N,2\mu} f^{(2)} \Big\}^L   +
 \Big\{  \Pi^L_{N,2\mu} f^{(1)},  {\hat f}^{(2)} \Big\}^L  \nonumber \\
 &  + &
 \Big\{ {\hat f}^{(1)}, \Pi^L_{N,\mu} f^{(2)} \Big\}^{x,y}   +
\Big\{  \Pi^L_{N,\mu} f^{(1)}, {\hat f}^{(2)} \Big\}^{x,y} \\
& + &
N \Big\{  \Pi_{|k| \geq N^b} f^{(1)},  f^{(2)} \Big\} +
N \Big\{  \Pi_{|k| < N^b} f^{(1)},   \Pi_{|k| \geq N^b} f^{(2)} \Big\} \Big)
\end{eqnarray*}
and the bound  \eqref{stizza}  follows again by
\eqref{commXHK},  \eqref{proiezaa},
 \eqref{b12}, \eqref{smoothl}, \eqref{Leviatan}.
Let consider only the term 
$N \Big\{  \Pi_{|k| \geq N^b} f^{(1)},  f^{(2)} \Big\}=:g $,
the last one being analogous.
We first use Lemma \ref{cauchy} with
$r'\rightsquigarrow r_1$, $r \rightsquigarrow r,$
$s'\rightsquigarrow s_1$ and $s\rightsquigarrow s_1+\s/2$,
where $\s:=s-s_1.$
Since
$ \Big( 1 - \frac{s_1 }{s_1 + \s / 2 } \Big)^{-1}
\leq 2 \Big( 1 - \frac{s_1}{s } \Big)^{-1} \leq 2 \d^{-1}
$
with the $ \d $ in \eqref{deltaciccio}, by \eqref{commXHK}
we get
\begin{eqnarray*}
\| {X}_g\|_{s_1,r_1} & \leq &  C(n)\d^{-1} N
\|  {X}_{\Pi_{|k| \geq N^b} f^{(1)}}\|_{s_1+\s/2,r}
\|  {X}_{f^{(2)}}\|_{s,r} \\
&  \stackrel{\eqref{smoothl}} \leq &
 C(n)\d^{-1} N\frac{s}{s_1} e^{-N^b (s -s_1)/2} \| {X}_{f^{(1)}}\|_{s,r}
\|  {X}_{f^{(2)}}\|_{s,r} \\
&  \stackrel{\eqref{Leviatan}} \leq &
C(n)\d^{-1}  \| {X}_{f^{(1)}}\|_{s,r} \|  {X}_{f^{(2)}}\|_{s,r} \,,
\end{eqnarray*}
for every $ N \geq N_1 $.
The proof of Proposition \ref{festa} is complete.
\rule{2mm}{2mm}

\medskip

The quasi-T\"oplitz character of a function is preserved under the flow
generated by a quasi-T\"oplitz Hamiltonian.

\begin{proposition}\label{main} {\bf (Lie transform)}
Let $ f, g \in {\mathcal Q}^T_{s,r} (N_0,  \teta,\mu) $ and let
$ s/2 \leq s' < s $, $ r/2\leq r' < r $.
 There is $  c(n) >  0  $ 
  such that, if
\be\label{piccof}
 \| f \|_{s,r, N_0, \teta, \mu}^T  \leq c(n) \, \d  \, ,
\ee
with $ \d $ defined in  \eqref{diffusivumsui},
then the hamiltonian flow of $ f $ at time $ t = 1 $,
$ \Phi_f^1 : D(s',r') \to D(s,r) $
is well defined, analytic and symplectic, and,
for  
\begin{equation}\label{stoppa}
    N_0' \geq \max \{N_0, \bar N \} \, , \quad
    \bar N := \exp \Big(
\max\Big\{
\frac{2}{b},\frac{1}{L-b},\frac{1}{1-L},8
\Big\} \Big)\, ,
\end{equation}
$($recall \eqref{figaro}$)$,
 $ \mu' < \mu $,
$  \teta' >  \teta $,  satisfying
\begin{equation}\label{Giobbe}
    \kappa  (N_0')^{b-L} \ln N_0'  \leq \mu - \mu' \,, \
(\CC  + \kappa) (N_0')^{L-1} \ln N_0' \leq   \teta'  -  \teta \, ,
\
 2(N_0')^{-b}\ln^2 N_0' \leq b(s-s')\,,
\end{equation}
we have
$ e^{ {\rm ad}_f } g \in {\mathcal Q}^T_{ s' , r' } ( N_0',  \teta', \mu' ) $
and
\be\label{gPhif}
\| e^{{\rm ad}_f} g \|^T_{s', r', N_0',  \teta', \mu' } \leq
2 \| g \|^T_{s,r, N_0,  \teta, \mu} 
\, .
\ee
Moreover, for  $ h = 0,1, 2 $, and  coefficients $ 0 \leq b_j\leq 1/j!$, $j\in\mathbb{N},$
\be\label{gPhif12}
\Big\| \sum_{j \geq h}  b_j \,{\rm ad}_f^j (g)  \Big\|^T_{s', r',  N_0', \teta', \mu'} \leq
2 (C\d^{-1} \| f \|_{s,r,N_0, \teta, \mu}^T )^{h}
\| g \|^T_{s,r, N_0,  \teta, \mu} \, . 
\ee
\end{proposition}
Note that \eqref{gPhif} is \eqref{gPhif12} with $ h = 0 $, $ b_j := 1/j!$

\smallskip
\begin{pf}
Let us prove \eqref{gPhif12}.   We define 
$$
G^{(0)} := g \, ,  \quad  G^{(j)} := {\rm ad}_f^j (g) := {\rm ad}_f ( G^{(j-1)}) =
\{f , G^{(j-1)}\} \, ,  \  j \geq 1 \, ,
$$
and we split, for $ h = 0, 1, 2 $, 
\be\label{expone}
G^{\geq h} :=  \sum_{j \geq h} b_j G^{(j)} = \sum_{j=h}^{J-1}  b_j G^{(j)} +
\sum_{j \geq J} b_j G^{(j)} =: G_{< J}^{\geq h} + G_{\geq J} \, .
\ee
As in \eqref{XHk} we deduce
\be\label{Giterat}
\| {X}_{G^{(j)}} \|_{s',r'} \leq ( C(n)  j \d^{-1} )^j \| {X}_f \|_{s,r}^j  \| {X}_g \|_{s,r} \, , \quad \forall  j   \geq 0  \, ,
\ee
where $ \delta $ is defined in \eqref{diffusivumsui}.
Let
\be\label{defet}
\eta := C(n)  e \d^{-1} \| {X}_f \|_{s,r} < 1 / (2e)
\ee
(namely take $c(n) $ small in \eqref{piccof}).
By \ref{Giterat}, using $j^j b_j\leq j^j  / j! < e^j $, we get
\be\label{code}
\| {X}_{G_{\geq J} } \|_{s', r'} \leq
\sum_{j \geq J} b_j ( C(n)  j  \d^{-1} \| {X}_f \|_{s,r}  )^j  \| {X}_g \|_{s,r}
\leq 2 \eta^J \| {X}_g \|_{s,r} \, .
\ee
In particular, for $ J = h=0,1,2 $, we get
\begin{equation}\label{rubamazzo}
\|{X}_{G^{\geq h}}\|_{s',r'} \leq 2 \eta^h\| {X}_g \|_{s,r}\,.
\end{equation}

For any $N\geq N_0'$
 we choose
\be\label{choiceJN}
J :=  J(N) :=  \ln N \, ,
\ee
and  we  set
$$
G'_N:=G_{<J}^{\geq h}\,,\quad G''_N :=G_{\geq J} \, , \quad G^{\geq h} = G'_N + G''_N \,.
$$
Then \eqref{gPhif12} follows by Lemma \ref{briscola}
(with $N_*\rightsquigarrow N_0', s \rightsquigarrow s', r \rightsquigarrow r', $
$ \teta \rightsquigarrow  \teta', \mu \rightsquigarrow \mu'$) and  \eqref{rubamazzo},
 once we show that
\begin{equation}\label{scopone}
\|G'_N\|_{s',r',N, \teta',\mu'}^T\leq \frac32 \eta^h \| g \|^T_{s,r} \, ,
\quad
N\|{X}_{G''_N}\|_{s',r'}\leq \frac12 \eta^h \| g \|^T_{s,r}
\end{equation}
with $h=0,1,2$
(for simplicity $\| g \|^T_{s,r}:=\| g \|^T_{s,r,N_0, \teta,\mu}$).

For all $ N \geq N_0'\geq e^8 $ (recall \eqref{stoppa}),
\begin{eqnarray}\label{GNhat}
N \| {X}_{G_{\geq J} } \|_{s',r'} & \stackrel{\eqref{code}} \leq &
N 2 \eta^J \| {X}_g \|_{s,r}  \leq  \eta^h ( N 2 \eta^{J-h}) \| g \|_{s,r}^T
\nonumber \\
&   \stackrel{\eqref{defet}} \leq &
\eta^h 2^{-J + h+1} e^h N  e^{-J} \| g \|_{s,r}^{T}  \leq  \frac{\eta^h}{2} \| g \|^T_{s,r} \, ,
\end{eqnarray}
proving the second inequality in \eqref{scopone}.
Let us prove the first  inequality in \eqref{scopone}.
\\[1mm]
{\sc Claim}: {\it $ \forall j = 1, \ldots, J - 1 $, we have
$G^{(j)}   \in  {\mathcal Q}_{s',r'}^T (N,  \teta', \mu' ) $ and
\be\label{adfjg}
\| G^{(j)} \|^T_{ r',  s', N,  \teta'  , \m' }  \leq \| g \|^T_{s,r} ( C'  j \d^{-1} \| f \|^T_{s,r} )^j
\ee
$($for simplicity $ \| f \|^T_{s,r} := \| f \|^T_{s,r, N_0,  \teta, \mu}) $.} This claim implies (using $j^j b_j< e^j$)
\begin{eqnarray*}
\Big\| \sum_{j=h}^{J-1} b_j \,G^{(j)}\Big\|^T_{s', r',  N, \teta',\mu'}
& \stackrel{\eqref{adfjg}} \leq &
\sum_{j=h}^{J-1}   b_j \| g \|^T_{s,r}  (C'  j \d^{-1} \| f \|^T_{s,r} )^j  \\
& \stackrel{\eqref{defet}}  \leq &  \| g \|^T_{s,r} \sum_{j=h}^{+\infty} \eta^j
\leq \frac32 \eta^h \|g\|^T_{s,r}
\end{eqnarray*}
 for $ c $ small enough in \eqref{piccof}.
This  proves the first  inequality in \eqref{scopone}.

\smallskip

Let us prove the claim. Fix $ 0 \leq j \leq J -1 $. We define, $ \forall i = 0, \ldots , j $,
\be\label{defsmt}
\mu_i := \mu - i \, \frac{\mu - \mu'}{j}, \  \teta_i :=  \teta + i \frac{ \teta' -  \teta}{j}, \
r_i := r -  i \, \frac{r - r'}{j}, \  s_i := s -  i \, \frac{s - s'}{j} \, ,
\ee
and we  prove inductively that, for all $ i = 0, \ldots , j $,
\be\label{agi}
\| {\rm ad}_f^i (g) \|^T_{s_i, r_i, N,  \teta_i, \mu_i }  \leq ( C'  j \d^{-1} \|f \|^T_{s,r} )^i
\| g \|^T_{s,r} \, ,
\ee
which, for $ i = j $, gives \eqref{adfjg}.
For $ i = 0 $, formula
\eqref{agi}  follows because $ g \in {\mathcal Q}^T_{s,r} (N_0,  \teta, \mu) $ and Lemma \ref{primo}.

Now assume that \eqref{agi} holds for $i$ and prove it for $i+1$.
We want to apply Proposition \ref{festa} to the functions $f$ and
${\rm ad}_f^i (g)$ with $N_1 \rightsquigarrow N$, $s \rightsquigarrow
s_i, s_1 \rightsquigarrow s_{i+1},$ $ \teta \rightsquigarrow  \teta_i,
 \teta_1 \rightsquigarrow  \teta_{i+1},$ etc. We have to verify conditions
\eqref{Leviatan} that reads
\begin{eqnarray}
&&\kappa N^{b-L}  < \mu_{i} - \mu_{i+1} \,, \quad \mu_{i}
N^{L-1} + \kappa N^{b-1} <   \teta_{i+1}  -  \teta_i \, ,  \label{eccolequi1}
\\
&& \ 2N e^{- N^b \frac{s_i - s_{i+1}}{2}} < 1 \,, \quad b(s_i -
s_{i+1}) N^b >2\,. \label{eccolequi2}
\end{eqnarray}
Since, by \eqref{defsmt},
$$
\mu_i - \mu_{i+1} = \frac{\mu - \mu'}{j} \, , \quad
\teta_{i+1} - \teta_{i} = \frac{\teta - \teta'}{j} \, , \quad s_i - s_{i+1} = \frac{s - s'}{j}
$$
and  $ j < J = \ln N$ (see \eqref{choiceJN}), $ 0 < b < L < 1 $ (recall \eqref{figaro}),
$ \mu'\leq \mu \leq \CC $,  the above conditions \eqref{eccolequi1}-\eqref{eccolequi2} are implied by
\begin{eqnarray}
&&\kappa  N^{b-L} \ln N  < \mu - \mu' \,, \quad
(\CC  + \kappa) N^{L-1} \ln N <   \teta'  -  \teta \, ,
\nonumber\\
&&
\ 2N e^{- N^b (s-s')/2\ln N} < 1\,,
\quad
b(s-s') N^b > 2\ln N\, . \label{last2}
\end{eqnarray}
The last two conditions \eqref{last2} are implied by
$ b(s-s') N^b > 2\ln^2 N $ and  since $N\geq e^{1/1-b}$ (recall \eqref{stoppa}).
Recollecting we have to verify
\begin{equation}\label{uomonero}
\kappa  N^{b-L} \ln N  \leq \mu - \mu' \,, \ \
(\CC  + \kappa) N^{L-1} \ln N \leq   \teta'  -  \teta \, ,
\ \  2N^{-b}\ln^2 N \leq b(s-s') \,.
\end{equation}
Since the function $ N \mapsto N^{-\g}\ln N$ is decreasing
for $ N \geq e^{1/\g} $,
we have that \eqref{uomonero} follows by \eqref{stoppa}-\eqref{Giobbe}.
Therefore Proposition \ref{festa} implies that
$ {\rm ad}_f^{i+1} (g) \in  {\mathcal Q}^T_{s_{i+1},r_{i+1}} (N,  \teta_{i+1}, \mu_{i+1} ) $
and, by
\eqref{poisbound2}, \eqref{inscatola},  we get
\begin{equation}\label{polpettone}
\| {\rm ad}_f^{i+1} (g) \|^T_{s_{i+1}, r_{i+1}, N,  \teta_{i+1}, \mu_{i+1} }
\leq
C'  \d_i^{-1} \| f \|^T_{s, r} \| {\rm ad}_f^{i} (g)  \|^T_{s_i, r_i, N,  \teta_i, \mu_i}
\end{equation}
where
\be\label{deltaii}
\d_i := \min\left\{ 1 - \frac{s_{i+1}}{s_i}\,,\,
1-\frac{r_{i+1}}{r_i}   \right\}
 \geq  \frac{\d}{j}
\ee
and $ \d $ is defined in \eqref{diffusivumsui}.
Then 
\begin{eqnarray*}
\| {\rm ad}_f^{i+1} (g) \|^T_{s_{i+1}, r_{i+1}, N,  \teta_{i+1}, \mu_{i+1} }
& \stackrel{\eqref{polpettone}, \eqref{deltaii}} \leq &
C'  j \d^{-1}    \| f \|^T_{s, r, N_0,  \teta, \mu}  \| {\rm ad}_f^{i} (g)  \|^T_{s_i, r_i, N,  \teta_i, \mu_i} \\
& \stackrel{\eqref{agi}} \leq &( C'  j \d^{-1} \|f \|^T_{s,r} )^{i+1}  \| g \|^T_{s,r}
\end{eqnarray*}
proving \eqref{agi}
 by  induction.
\end{pf}

\section{An abstract KAM theorem}\label{sec:4}\setcounter{equation}{0}

We consider a family of integrable Hamiltonians
\be\label{NormalN}
{\cal N} := {\cal N}(x,y,z,\bar z; \xi) :=  e(\xi) +  \o (\xi) \cdot y + \Om (\xi) \cdot z {\bar z}
\ee
defined on
$ \mathbb{T}^n_s  \times \mathbb{C}^n \times \ell^{a,p}_{\cal I} \times \ell^{a,p}_{\cal I}  $,
where $ {\cal I} $ is defined in \eqref{cC},
 the tangential  frequencies $ \om := (\om_1, \ldots, \om_n) $ 
and the normal frequencies
$ \Om := (\Om_j)_{j \in \Z \setminus {\cal I}} $ 
depend on $ n $-parameters
$$
\xi \in {\cal O} \subset \mathbb{R}^n \, . 
$$
For each $ \xi $ there is an invariant $ n$-torus
$ {\cal T}_0 = \mathbb{T}^n \times \{0 \} \times \{0 \} \times \{0 \} $
with frequency $ \om ( \xi ) $. In its normal space, the origin $ (z,\bar z ) = 0  $ is an elliptic fixed point with
proper frequencies $ \Om (\xi ) $.
The aim is to prove the persistence of a large portion of this family  of linearly stable
tori under small analytic perturbations $ H = {\cal N} + P $.

\smallskip

\begin{description}
\item $ {\bf (A1)} $ {\sc Parameter dependence}.
The map $ \om : {\cal O} \to {\mathbb R}^n $,
$ \xi \mapsto \omega ( \xi )$,  is Lipschitz continuous.
\end{description}
With in mind the application to NLW we assume
\begin{description}
\item $ {\bf (A2)} $
{\sc Frequency asymptotics}. We have
\begin{equation}\label{carbonara}
\Om_j (\xi) =   \sqrt{j^2+ \mm}  + a (\xi) \in {\mathbb R}\, , \ \  j \in \Z \setminus {\cal I} \, ,
\end{equation}
for some Lipschitz continuous functions $ a (\xi) \in \R $.
 \end{description}

By (A$1$) and (A$2$),  the Lipschitz semi-norms 
of the frequency maps satisfy, for some $1\leq M_0<\infty $,
\begin{equation}\label{M}
|\o|^{\rm lip} \, , \ | \O |^{\rm lip}_{\infty} \leq M_0
\end{equation}

 where the Lipschitz semi-norm is
\be\label{lipfre}
| \O |^{\rm lip}_{\infty} :=
| \O |^{\rm lip}_{\infty,{\cal O}}:=
\sup_{\xi,\eta\in {\cal O}, \xi\neq\eta} \frac{ |\Om(\xi)-
\Om(\eta)|_{\infty}}{| \xi - \eta |}
\ee
and $ | z |_\infty := \sup_{j \in \Z \setminus {\cal I}} | z_j | $. Note that
by the Kirszbraun theorem (see e.g. \cite{KP}) applied  componentwise
we can extend $\o,\O$ on the whole
 $ \R^n$ with the same  bound \eqref{M}.

\begin{description}
\item $ {\bf (A3)} $  {\sc Regularity}.
The perturbation $ P : D(s,r) \times {\cal O} \to \C $
is  $ \l $-regular (see Definition \ref{Hregular}).
\end{description}
In order to obtain the asymptotic expansion  \eqref{freco}  for the perturbed frequencies we also assume
\begin{description}
\item $ {\bf (A4)} $  {\sc Quasi-T\"oplitz}. The perturbation $ P $ (preserves momentum and)
is quasi-T\"oplitz (see Definition \ref{topbis_aa}).
\end{description}
Thanks to the conservation of momentum we restrict to the set of indices
\begin{eqnarray}\label{masotti}
{\bf I } :=  \Big\{ && (k,l) \in \Z^n \times \Z^\infty, (k,l) \neq (0,0) \, , |l| \leq 2,  \ {\rm where} \\
& & {\rm or} \  l  = 0\, , \ k \cdot \pluto = 0 \, ,  \nonumber \\
& & {\rm or} \ l  = \s e_m \, , m \in \Z \setminus {\cal I} \, , \ k \cdot \pluto + \s m = 0 \, ,  \nonumber  \\
& &   {\rm or} \ l  = \s e_m + \s'  e_n \, , m,n \in \Z \setminus {\cal I}  \, ,   \ k \cdot \pluto + \s m + \s' n = 0 \nonumber
\Big\} \, .
\end{eqnarray}
\noindent
Let
\be\label{P0Pbar}
P = P_{00}(x) + \bar P (x,y, z, \bar z ) \quad {\rm where} \quad
\bar P (x, 0, 0, 0)  = 0 \, .
\ee
\begin{theorem} {\bf (KAM theorem)} \label{thm:IBKAM}
Suppose that $ H = {\cal N} + P  $ satisfies $ {\rm (A1)} $-${\rm (A4)}$ with
$ s,r >  0 $, $ \cc < \teta,\mu < \CC $, $ N > 0 $.
Let $ \g>0 $ be a small  parameter and  set
\be\label{KAMconditionT}
\e := \max \Big\{ \gamma^{- 2/3} | {X}_{P_{00}}|_{s,r}^\l \, ,\ \gamma^{- 2/3}
\| X_{P_{00}} \|_{s, r} \, ,     \
 \g^{-1} | {X}_{\bar P} |_{s,r}^\l\,,\
 \g^{-1} \| \bar P \|_{s, r, N,\teta, \mu}^T
\Big\}\,,\quad   \l := \g / M_0 \,.
\ee
If $\e$ is small enough, then  there exist:
\\[1mm] $ \bullet $ {\bf (Frequencies)}
Lipschitz functions $ \o^\infty:  \R^n \to {\mathbb R}^n $,
$ \O^\infty :\R^n \to  \ell_\infty $, $ a^\infty_{\pm}:
\R^n\to {\mathbb R} $, such that
\begin{equation}\label{muflone}
| \o^\infty - \o |+\l | \o^\infty - \o |^{\rm lip} \, ,
\ \  | \O^\infty -  \O  |_{\infty}+
\l| \O^\infty -  \O  |_{\infty}^{\rm lip}\
 \leq \ C \g \e \, , \quad | a^\infty_\pm | \leq C \g \e \, ,
\end{equation}
\be\label{freco}
\sup_{\xi \in \R^n }|\O^\infty_j(\xi) -  \O_j(\xi) - a^\infty_{\mathtt s(j)}(\xi) |
\leq   \g^{2/3} \e \, \frac{C}{|j|} \, , \quad \forall |j| \geq C_\star \g^{-1/3} \, .
\ee
\\[1mm] $ \bullet $  {\bf (KAM normal form)}
A Lipschitz family of analytic symplectic maps
\be\label{eq:Psi}
\Phi : D(s/4, r/4)\times  {\cal O}_\infty \ni(x_\infty,y_\infty,w_\infty;\xi)  \mapsto (x,y,w)\in D(s, r)
\ee
close to the identity where
\begin{eqnarray}\label{Cantorinf}
{\cal O}_\infty & := & \Big\{ \xi \in {\cal O} \ \, : \,
\  |\o^\infty(\xi ) \cdot k +
\O^\infty (\xi) \cdot l | \geq \frac{2 \g}{1+|k|^\tau}\, , \, \forall \, (k,l) \in {\bf I}
\ {\rm   \ defined \ in \ } \eqref{masotti} \, ,
\nonumber  \\
& & \ |\o^\infty(\xi)\cdot k+p|\geq \frac{2\g^{2/3}}{1+|k|^\t}, \, \forall  k\in \Z^n,\, p\in\Z\,, \, (k,p) \neq (0,0) \, , \, \t > 1 / b  \
{\rm see} \, \eqref{figaro},
 \nonumber
\\
& & \ |\o(\xi)\cdot k|
\geq \frac{2\g^{2/3}}{1+|k|^n}\,, \ \forall\,
0<|k|< \g^{-1/(7n)}\
\Big\}
\end{eqnarray}
such that, $\forall \xi\in \mathcal O_\infty$:
\be\label{Hnew}
H^\infty (\cdot ; \xi ) :=  H \circ \Phi (\cdot ; \xi) =
\o^\infty(\xi) \cdot y_\infty + \O^\infty(\xi) \cdot z_\infty \bar z_\infty + P^\infty
\quad {has} \quad P^\infty_{\leq 2} = 0 \, .
\ee
Then, $ \forall \xi \in {\cal O}_\infty $, the map $  x_\infty \mapsto  \Phi (x_\infty, 0,0;\xi)$
is a real analytic embedding of an elliptic, $ n$-dimensional torus with frequency $ \o^\infty (\xi ) $
for the system with Hamiltonian $ H $.
\end{theorem}

The main novelty of Theorem \ref{thm:IBKAM}
is the  asymptotic decay  \eqref{freco} of the perturbed frequencies.
 In order to prove  \eqref{freco} we use the
 quasi-T\"oplitz property (A4) of the perturbation.
 The reason for introducing in  \eqref{KAMconditionT} conditions
for both the Lipschitz-sup  and the T\"oplitz-norms 
is the following.  For the measure estimates, we need the usual Lipschitz dependence
of the perturbed frequencies
with respect to the parameters, see \eqref{muflone}. This is derived as in \cite{Po2} and \cite{BB10}.
On the other hand,
we do not need (in section \ref{sec:meas}) a Lipschitz estimate on $ a^\infty_\pm $ (that, in any case, could  be obtained).
For this reason,
we do not introduce the Lipschitz dependence in the T\"oplitz norm.

In the next Theorem \ref{thm:measure} we verify the second order Melnikov non-resonance conditions thanks to
\begin{enumerate}
\item
the  asymptotic decay  \eqref{freco} of the perturbed frequencies,
\item
the restriction to indices $ (k,l) \in   {\bf I}  $ in \eqref{Cantorinf} 
which follows by momentum conservation, see (A4).
\end{enumerate}

As in \cite{BB10}, the Cantor set of ``good" parameters
 $ {\cal O}_\infty $ in (\ref{Cantorinf}),
is expressed in terms of the final frequencies $ \om^\infty (\xi) $, $ \O^\infty (\xi) $
(and of the initial tangential frequencies $ \om (\xi) $)
and not inductively as, for example,  in \cite{Po2}.
 This simplifies the measure estimates.

\begin{theorem}
{\bf (Measure estimate)} \label{thm:measure}
Let $ {\cal O} := [\rho/2, \rho]^n $, $ \rho >  0 $. Suppose
\be\label{omegaxi}
\om (\xi) = \bar \om + A \xi \, , \ \bar \om \in \R^n \, , \ A \in {\rm Mat} (n\times n) \ \, ,
\quad \Om_j (\xi ) = \sqrt{j^2 + \mm} + \vec a \cdot \xi \, , \ a \in \R^n
\ee
and assume
the non-degeneracy condition:
\begin{equation}\label{pota}
 A  \ {\rm  \ invertible}  \quad \ {\rm and } \quad   \
2 (A^{-1})^T \vec  a \notin \Z^n \setminus \{ 0 \} \, .
\end{equation}
Then,  the Cantor like set  $ {\cal O}_\infty $
defined in
\eqref{Cantorinf}, with exponent
\be\label{antiochia}
\t > \max\{2 n + 1, 1/ b \}
\ee
($ b $ is fixed in \eqref{figaro}),
satisfies
\be\label{consolatrixafflictorum}
| {\cal O} \setminus {\cal O}_\infty | \leq C(\t) \rho^{n-1} \g^{2/3} \, .
\ee
\end{theorem}

\noindent
Theorem \ref{thm:measure} is proved in section \ref{sec:meas}.
The asymptotic estimate \eqref{freco} is used for
the key inclusion \eqref{claimrq1}.

\section{Proof of the KAM Theorem \ref{thm:IBKAM}}\setcounter{equation}{0} \label{sec:5}

In the following by $ a \lessdot b $
we mean that there exists  $ c  >  0 $ depending
only on $n,\mm, \kappa$ such that
$a\leq  c b$.

\subsection{First step}

We perform a preliminary change of variables to improve the
smallness conditions.
For all $\xi$ in
\begin{equation}\label{Ozero}
\mathcal{O}_*:= \Big\{\ \xi\in\mathcal O \ :\
|\o(\xi)\cdot k|
\geq \frac{\g^{2/3}}{1+|k|^n}\,, \ \forall\,
0<|k|< \g^{-1/(7n)}\
\Big\}
\end{equation}
we consider the solution
\be\label{def:F0}
F_{00}(x) := \sum_{0<|k|< \g^{-1/(7n)}}
\frac{P_{00,k}}{\ii \o(\xi) \cdot k} e^{\ii k \cdot x}
\ee
of the homological equation
\be\label{F0P0}
- {\rm ad}_{\mathcal N}  F_{00} + \Pi_{|k|<\g^{-1/(7n)}} P_{00}(x) = \langle P_{00}\rangle \, .
\ee
Here $P_{00}$ is defined in \eqref{P0Pbar} and $\langle \cdot \rangle$ denotes the mean value on the angles.
Note that for any  function $ F_{00} (x) $ we have
$ \| F_{00} \|^T_{s,r} = \| X_{F_{00}} \|_{s,r} $, see Definition \ref{topbis_aa}.
We want to apply Proposition \ref{main}  with
$ s, r, s', r'  \rightsquigarrow 3s/4, 3r/4, s/2, r/2 $. The condition \eqref{piccof} is verified
because 
$$
\| F_{00} \|_{3s/4,r}^T = \| X_{F_{00}} \|_{3s/4,r} \stackrel{\eqref{def:F0}, \eqref{Ozero}, \eqref{caligola}}
\leq C(n, s) \g^{-2/3} \| X_{P_{00}} \|_{s,r} \stackrel{\eqref{KAMconditionT} }  \leq C(n,s) \e
$$
and $ \e$ is sufficiently small. Hence  the time--one flow
\begin{equation}\label{Phi00}
\Phi_{00} :=  \Phi^1_{F_{00}} : D(s_0,r_0)\times \mathcal{O}_*\to D(s,r)
    \quad {\rm with}\quad s_0 := s/2 \, , \  r_0 := r/2 \, ,
\end{equation}
is well defined, analytic, symplectic. Let
 $ \mu_0 <  \mu $, $ \teta_0 >  \teta $, $ N_0 >  N $ large enough,
 so that \eqref{Giobbe} is satisfied with $ s, r, N_0,  \teta, \mu, \rightsquigarrow  s, r, N,  \teta ,
 \mu $ and $ s', r' , N_0',  \teta' , \mu' \rightsquigarrow s_0, r_0, N_0, \teta_0, \mu_0 $. Note that here $N_0$ is independent of $\g$.
Hence \eqref{gPhif} implies
\be\label{zione}
\| e^{{\rm ad}_{F_{00}} }  \bar P \|_{s_0,r_0, N_0, \teta_0, \mu_0}^T \leq 2  \| \bar P \|_{s, r, N, \teta,\mu}^{T}  \, .
\ee
Noting that $ e^{{\rm ad}_{F_{00}}} P_{00} = P_{00} $ and
$ e^{{\rm ad}_{F_{00}}} {\cal N} =  {\cal N} + {\rm ad}_{F_{00}}  {\cal N} $
the new Hamiltonian  is
\begin{eqnarray}
H^0  :=  e^{{\rm ad}_{F_{00}}} H & = &  e^{{\rm ad}_{F_{00}} } {\cal N} +   e^{{\rm ad}_{F_{00}} }   P_{00} +  e^{{\rm ad}_{F_{00}} }  \bar P
 =    {\cal N} + {\rm ad}_{F_{00}}  {\cal N}  +  P_{00} +  e^{{\rm ad}_{F_{00}} }  \bar P  \nonumber \\
& \stackrel{\eqref{F0P0}}  =& \big( \langle P_{00} \rangle +  {\cal N}\big) +
\Big(\Pi_{|k|\geq \g^{-1/(7n)}}P_{00}+ e^{{\rm ad}_{F_{00}} }  \bar P \Big)=: {\cal N}_0 + P_0 \, . \label{accazero}
\end{eqnarray}
By \eqref{smoothl} (and since $ P_{00} (x)$ depends only on $ x $)
\be\label{restialti}
\Big\| \Pi_{|k|\geq \g^{-1/(7n)}}P_{00} \Big\|_{3s/4,r}^T\leq 4 e^{-s\g^{-1/(7n)}/4} \|X_{P_{00}}\|_{s,r}\stackrel{\eqref{KAMconditionT}} \leq 4 \g^{2/3}e^{-s\g^{-1/(7n)}/4} \e \leq \g \e\,,
\ee
for $\g$ small.
By \eqref{restialti}, \eqref{zione} and \eqref{KAMconditionT} we get
\be\label{piccolezza0}
\| P_0 \|_{s_0, r_0, N_0, \teta_0,\mu_0}^{T} <  3 \g \e \, .
\ee
In the same way, since
$
|X_{F_{00}}|^\lambda_{3s/4,r}\leq C(n, s) \g^{-2/3} | X_{P_{00}} |_{s,r}^\lambda
$, we also obtain the Lipschitz estimate
\be\label{piccolezza0lip}
| X_{P_0} |_{s_0, r_0}^\l <  3 \g \e \, .
\ee

\subsection{KAM step}

We now consider the generic KAM step for an Hamiltonian
\be\label{def:R}
H =  {\cal N} +  P =  {\cal N} +  P^{\leq 2}_K + ( P - P^{\leq 2}_K  ) 
\ee
where  $ P^{\leq 2}_K $ are defined as in \eqref{Pleq2}.

\subsubsection{Homological equation}

\begin{lemma}\label{paperino}
Assume that
\be\label{Omasin}
 | \O_j -\sqrt{j^2+\mm}- a_{\mathtt s(j)}|\leq \frac{\g }{|j|} \, ,
\quad \forall\, |j| \geq j_* \, ,
\ee
for some $a_+,a_-\in\R.$ Let
$ \D_{k,m,n} := \o\cdot k + \O_m -  \O_n $, $ \tilde\D_{k,m,n} := \o\cdot k + |m| -  |n| $.

If  $ |m|,|n| \geq \max\{ j_*,\sqrt \mm \}$ and  $ \mathtt s (m) = \mathtt s(n) $,
then \be\label{laurina} |\D_{k,m,n} - \tilde\D_{k,m,n} |  \leq
\frac{\mm}{2} \frac{|m  - n|}{|n| |m|} +\g  \Big(
\frac{1}{|m|}+\frac{1}{|n|} \Big)+ \frac{\mm^2}{2}
\left(\frac{1}{|m|^3}+\frac{1}{|n|^3} \right)\,. \ee
\end{lemma}

\begin{pf}
For $0\leq x\leq 1$ we have $|\sqrt{1+x}-1-x/2|\leq x^2/2 $.
Setting $x:=\mm/n^2 $ (which is $ \leq 1$) and using
\eqref{Omasin}, we get
$$
\left| \O_n - |n| - \frac{\mm}{2|n|} - a_{\mathtt s(n)} \right| \leq
\frac{\g }{|n|} + \frac{\mm^2}{2|n|^3} \,.
$$
An analogous estimates holds for $\O_m $. Since $|\D_{k,m,n} - \tilde
\D_{k,m,n} |=|\O_m-|m|-\O_n + |n||$ the estimate \eqref{laurina} follows
noting that $a_{\mathtt s(m)}=a_{\mathtt s(n)}$.
\end{pf}

For a monomial $ \mathfrak m_{k,i,\a,\b} := e^{\ii k \cdot x} y^i z^\a {\bar z}^\b $
we set
\be\label{operatorQ}
[  \mathfrak m_{k,i,\a,\b} ] :=
\begin{cases}
 \mathfrak m_{k,i,\a,\b}
\qquad {\rm if} \ \  \ \, k = 0 \, , \ \a = \b  \cr
 0  \quad \qquad  \qquad \quad \ {\rm otherwise.}
\end{cases}
\ee
The following key proposition proves that the solution of the homological equation with a   quasi-T\"oplitz
datum is quasi-T\"oplitz.

\begin{proposition}\label{FP}
{\bf (Homological equation)} Let $ K \in \N $.
For all $ \xi \in {\cal O} $ such that 
\begin{equation}\label{foederisarca}
| \o(\xi) \cdot k +  \O(\xi) \cdot l | \geq \frac{\g}{ \langle k \rangle^\tau}\,,
\ \ \forall (k,l) \in {\bf I} \ (see \  \eqref{masotti}), \,  |k| \leq K \, ,
\end{equation}
then $ \forall P^{(h)}_K \in {\cal H}_{s,r}^{\rm null} $,
$ h = 0,1,2 $ $($see \eqref{phom}, \eqref{ProK}$)$,
the homological equations
\be\label{homoh}
- {\rm ad}_{\cal N} F^{(h)}_K  +
P^{(h)}_K =  [P^{(h)}_K] \,,\qquad h=0,1,2\,,
\ee
have a unique solution of the same form 
$  F^{(h)}_K \in {\cal H}_{s,r}^{\rm null}$ with $ [F^{(h)}_K ] =
0 $ and
\be\label{la1} \| {X}_{F^{(h)}_K} \|_{s,r} <  \g^{-1} K^{\t}  \|
{X}_{P^{(h)}_K}\|_{s,r}\,, \qquad | {X}_{F^{(h)}_K} |_{s,r}^\l \lessdot  \g^{-1} K^{\t+1}  |
{X}_{P^{(h)}_K}|_{s,r}^\l
\ee
where $ 2\g\l^{-1}\geq |\o|^{\rm lip} $, $ |\O|_\infty^{\rm lip} $.
In particular $ F_K^{\leq 2} := F_K^{(0)}
+ F_K^{(1)} + F_K^{(2)} $ solves
\be\label{homo}
- {\rm ad}_{\cal N} F_K^{\leq 2}  + P^{\leq 2}_K =  [P^{\leq 2}_K]\,.
\ee
Assume now that $ P^{(h)}_K  \in {\mathcal Q}^T_{s,r} (N_0,  \teta, \mu ) $ and
$ \O(\xi) $ satisfies \eqref{Omasin} for all $|j|\geq \teta N_0^*$ where
\be\label{N0large}
N_0^*:= \max\Big\{N_0 \, , \ \hat{c}\g^{-1/3} K^{\tau +1} \Big\}
\ee
for a constant $\hat{c} := \hat{c} (\mm , \kappa ) \geq 1 $. 
Then,
$ \forall \xi\in{\cal O}$ such that
\be\label{s+s-}
|\o (\xi)
\cdot k + p | \geq \frac{\g^{2/3}}{ \langle k \rangle^{\tau}} \,
, \ \ \forall  |k| \leq K, \ p \in \Z \, ,
\ee
we have $ F^{(h)}_K \in {\mathcal Q}^T_{s,r} (N_0^* , \teta, \mu) $,  $   h = 0, 1, 2 $,
and \be\label{Fijstimeh} \| F^{(h)}_K \|^T_{s,r, N_0^*,  \teta,
\mu} \leq 4\hat{c}
 \g^{-1} K^{2 \t}
\| P^{(h)}_K \|_{s,r, N_0,  \teta, \mu}^T \, . \ee
\end{proposition}

\begin{pf}
The solution of the homological equation \eqref{homoh} is
$$
F^{(h)}_K := - \ii \!\!\!\!\!\!\!\!\!\!\!\!\!\! \sum_{|k| \leq K, (k,i,\a,\b) \neq (0,i,\a,\a) \\
\atop 2 i + |\a|+ |\b| = h} \frac{P_{k,i,\a,\b}}{\Delta_{k,i,\a,\b} } e^{\ii k \cdot x} y^i z^\a \bar z^{\b} \, , \
\quad \Delta_{k,i,\a,\b}:=  \o(\xi) \cdot k +  \O(\xi) \cdot (\a-\b) \, .
$$
The divisors  $\Delta_{k,i,\a,\b}\neq 0, $ $ \forall (k,i,\a,\b)\neq (0,i,\a,\a)$, because 
$(k,i,\a,\b)\neq (0,i,\a,\a)$ is equivalent to $(k,\a-\b)\in {\bf I}$, and   the bounds \eqref{foederisarca} hold.
Then the first estimates in \eqref{la1}  follows by Lemma \ref{Plowerb}. The Lipsichtz estimate in \eqref{la1}
is standard, see e.g. Lemma 1 (and the first comment after the statement)
of \cite{Po2}. We just note that the Melnikov condition used in \cite{Po2} follows by \eqref{foederisarca} and momentum
consevation, e.g.
$$
| \om \cdot k + \Om_m - \Om_n | \stackrel{\eqref{foederisarca}}
\geq \frac{\g}{ \langle k \rangle^\t} \stackrel{\eqref{zmc}} = \frac{\g |m-n|}{| \pluto \cdot k | \langle k \rangle^\t}
\geq \g \frac{|m-n|}{\kappa \langle k \rangle^{\t+1}} \, .
$$
For the T\"oplitz estimate notice that the  cases $h=0,1$ are trivial since $\Pi_{N, \teta, \mu} F^{\leq 1}_K = 0 $.
When $ h = 2 $ 
 we first consider the subtlest case when $ P^{(2)}_K $ contains only the monomials with $ i = 0 $, $|\a | = | \b | = 1 $ (see \eqref{phom}),
namely \be\label{sump1}{\mathcal P} :=  P^{(2)}_K=  \sum_{|k| \leq K, m,n
\in \Z \setminus {\cal I}} P_{k, m, n} e^{\ii k \cdot x} z_m {\bar
z}_n \, , \ee
and, because of the conservation of momentum,
 the indices $ k, m , n $ in \eqref{sump1} 
 are restricted to
\be\label{zmc}
 \pluto \cdot k + m - n = 0 \, .
 \ee
The unique solution
$ F^{(2)}_K $ of  \eqref{homoh} with $ [F^{(2)}_K ] = 0 $ 
is
\be\label{pionab} {\mathcal F}:=
F^{(2)}_K :=  - \ii \!\!  \!\! \!\!  \!\! \!  \sum_{|k| \leq K, (k, m,n) \neq (0,m,m)}
 \!\! \frac{P_{k, m, n}   }{  \Delta_{k,m,n}  } e^{\ii k \cdot x} z_m {\bar z}_n \, ,
\ \ \Delta_{k,m,n}  := \o(\xi) \cdot k +  \O_m (\xi) - \O_n (\xi) \,
\ee
Note that by \eqref{foederisarca} and
\eqref{zmc} we have $ \Delta_{k,m,n}  \neq 0 $ if and only if  $ (k,m,n) \neq (0,m,m) $.

 \smallskip
Let us prove \eqref{Fijstimeh}.
For all $ N \geq N_0^* $
\be\label{PIFF}
\Pi_{N, \teta, \mu} \mathcal F =  - \ii \sum_{|k| \leq K, |m|,|n| >  \teta N}
\frac{P_{k, m, n}}{  \Delta_{k,m,n} } e^{\ii k \cdot x} z_m {\bar z}_n ,
\ee
and note that $e^{\ii k \cdot x}$ is $(N,\mu)$-low
momentum since  $ |k| \leq K < (N_0^*)^b \leq N^b $ by \eqref{N0large} and $ \t > 1 / b $.
By assumption $ \mathcal P \in {\cal Q}^T_{s,r, N_0,  \teta, \mu} $ and so,
recalling formula \eqref{deftop}, we may  write,
$ \forall N \geq  N_0^* \geq N_0 $,
\be\label{Pappro} \Pi_{N, \teta,\mu} \mathcal P =
\tilde {\mathcal P} + N^{-1} \hat {\mathcal P} \quad {\rm with} \quad \tilde {\mathcal P} := \sum_{|k| \leq K,
|m|,|n| >  \teta N} {\tilde P}_{k, m - n}  e^{\ii k \cdot x} z_m
{\bar z}_n \in {\cal T}_{s,r} (N,  \teta,\mu) \ee and
\be\label{Pgra} \| {X}_{\mathcal P} \|_{s,r}, \| {X}_{\tilde {\mathcal P}} \|_{s,r}, \|
{X}_{\hat {\mathcal P}} \|_{s,r} \leq 2 \| {\mathcal P} \|^T_{s,r} \, . \ee We now prove
that
\be\label{Ftil}
\tilde {\mathcal F} := \sum_{|k| \leq K, |m|,|n| >
\teta N} \frac{{\tilde P}_{k, m - n}}{ \tilde \D_{k, m,n}}  e^{\ii
k \cdot x} z_m {\bar z}_n \, , \quad {\tilde \Delta}_{k,m,n}  := \o (\xi)  \cdot k +  |m| - |n| \, ,
\ee
is a T\"oplitz approximation of $ \mathcal F $. Since
$ |m|, |n | >  \teta N \geq $ $  \teta N_0^* > $
$ \ccN_0^* \stackrel{\eqref{N0large}} > \kappa \, K \geq | \pluto \cdot  k | $ by  \eqref{caracalla},
we deduce by  \eqref{zmc} that $ m ,n $ have the same sign. Then
$$
\tilde \D_{k, m,n} =  \o (\xi)  \cdot k +  |m| - |n| =  \om (\xi)  \cdot
k + \mathtt s(m)(m - n) \, , \quad \mathtt s(m) := {\rm sign} (m) \, ,
$$
and $ \tilde {\mathcal F} $  in \eqref{Ftil} is  $ (N,  \teta,
\mu)$-T\"oplitz (see \eqref{marco}). Moreover, since $ |m| - |n |
\in \Z $, by \eqref{s+s-},  we get
\be\label{masciarelli}
|\tilde \D_{k, m,n}| \geq
\g^{2/3} \langle k \rangle^{- \tau} \, ,
\quad \forall  |k| \leq K, \ m, n,
\ee
and Lemma \ref{Plowerb} and \eqref{Ftil}  imply
\be\label{la2}
\| {X}_{\tilde{\mathcal F} } \|_{s,r} \leq  \g^{-2/3} K^{\t}  \| {X}_{\tilde {\mathcal P} } \|_{s,r} \, .
\ee
The T\"oplitz defect is
\begin{eqnarray}
N ^{-1} \hat {\mathcal F}  & := & \Pi_{N,\teta, \mu} {\mathcal F}  - \tilde {\mathcal F}  \label{N-1F}  \\
& \stackrel{\eqref{PIFF}, \eqref{Ftil}} = & \sum_{|k| \leq K,
|m|,|n| >  \teta N} \Big( \frac{P_{k, m, n}}{\D_{k, m,n}} -
\frac{{\tilde P}_{k, m-n}}{\tilde \D_{k, m,n}} \Big)
 e^{\ii k \cdot x} z_m {\bar z}_n  \nonumber  \\
 & = &  \sum_{|k| \leq K, |m|, |n| >  \teta N}
\Big[\Big( \frac{P_{k, m, n}}{\D_{k, m,n}} - \frac{P_{k,
m,n}}{\tilde \D_{k, m,n}}  \Big) +  \Big( \frac{P_{k, m, n} -
{\tilde P}_{k, m-n}}{\tilde \D_{k, m,n}}  \Big) \Big]
 e^{\ii k \cdot x} z_m {\bar z}_n   \nonumber \\
 & \stackrel{\eqref{Pappro}}
 = &  \sum_{|k| \leq K, |m|, |n| >  \teta N}
\Big[P_{k, m, n} \Big( \frac{ \tilde \D_{k, m,n} - \D_{k, m,n} }{
\D_{k, m,n} \tilde \D_{k, m,n} } \Big) +  N^{-1}  \frac{{\hat
P}_{k, m, n}}{ \tilde \D_{k, m,n} }  \Big]
 e^{\ii k \cdot x} z_m {\bar z}_n \, . \nonumber
\end{eqnarray}
By
\eqref{laurina}, $ |m|, |n| \geq  \teta N \geq \ccN $, and $|m-n|\leq \kappa K$ (see  \eqref{zmc}) we get,
taking $\hat{c} $ large enough,
\begin{eqnarray}
| \tilde \D_{k, m,n} - \D_{k, m,n} |
& \leq &
 \frac{\mm \kappa K}{2 \ccN^2}
 +\frac{2 \g }{\ccN} +\frac{\mm^2 }{\ccN^3}
\leq
 \frac{\hat{c}}{4 N}\left(\frac{K}{N} +\g\right)
\stackrel{\eqref{N0large}}\leq \min \left\{
\frac{\hat{c}\g^{1/3}}{2N}\,,\ \frac{\g^{2/3}}{2K^{\t}} \right\}
  \, .
  \label{DDkmn}
\end{eqnarray}
Hence
\be\label{daba} |  \D_{k, m,n}  |    \geq | \tilde \D_{k, m,n} | - | \tilde \D_{k, m,n} - \D_{k, m,n} |
\stackrel{\eqref{masciarelli}, \eqref{DDkmn}} \geq
\frac{\g^{2/3}}{\langle k \rangle^{\t}} -  \frac{\g^{2/3}}{2
K^{\t}} \geq \frac{\g^{2/3}}{2 \langle k \rangle^{\t}} \, .
\ee
Therefore \eqref{DDkmn}, \eqref{masciarelli}, \eqref{daba} imply
$$
\frac{ |\tilde \D_{k, m,n} - \D_{k, m,n}|}{| \D_{k, m,n}| |\tilde \D_{k, m,n}|}
 \leq \frac{\hat{c}\g^{1/3}}{2N}
  \frac{2 \langle k \rangle^{\t}}{\g^{2/3}}    \frac{ \langle k \rangle^{\t}}{\g^{2/3}}
  \leq   \frac{\hat{c}}{N \g} K^{2\t}
  $$
and \eqref{N-1F}, \eqref{masciarelli}, and 
Lemma \ref{Plowerb},  imply
\be\label{la3} \| {X}_{\hat {\mathcal F} }
\|_{s,r} \leq
 \hat{c} \g^{-1} K^{2\t} \| {X}_{\mathcal P}  \|_{s,r}
   +  \g^{-2/3} K^{\t}  \| {X}_{\hat {\mathcal P} } \|_{s,r}
\stackrel{\eqref{Pgra}} \leq 4 \hat{c} \g^{-1} K^{2 \t}  \| {\mathcal P}
\|_{s,r}^T  \, .
\ee
In conclusion \eqref{la1}, \eqref{la2}, \eqref{la3}  prove \eqref{Fijstimeh}  for $\mathcal F$.

Let us briefly discuss the case  when $ h = 2 $ 
 and  $ P^{(2)}_K $ contains only the monomials with $ i = 0 $, $|\a | =2$, $|\b|=0$ or viceversa (see \eqref{phom}).
Denoting
\be\label{sump1bis}{\mathcal P} :=  P^{(2)}_K :=   \sum_{|k| \leq K, m,n
\in \Z \setminus {\cal I}} P_{k, m, n} e^{\ii k \cdot x} z_m { z}_n \, ,
\ee
we have
$$
\Pi_{N, \teta, \mu} \mathcal F =  - \ii \sum_{|k| \leq K, |m|,|n| >  \teta N}
\frac{P_{k, m, n}}{ \o\cdot k+\O_m+\O_n}e^{\ii k \cdot x} z_m { z}_n
$$
where $ |\o\cdot k + \O_m + \O_n | > (|m|+|n|) / 2  > \teta N / 2 $ since $|m|,|n| >  \teta N$ and $|k| \leq K<N^b$.
In this case we may take as T\"oplitz approximation  $ \tilde {\mathcal F} = 0 $.
\end{pf}

\subsubsection{The new Hamiltonian $ H^+ $}

Let $ F =F^{\leq 2}_K$ be the solution  of the homological equation \eqref{homo}.
If, for $ s/2 \leq s_+ < s $, $ r/2\leq r_+ < r $,  the condition
\be\label{piccoF}
\| F \|_{s,r, N_0^*, \teta, \mu}^T
\leq c(n) \, \d_+ \, , \quad
\d_+ := \min \Big\{  1 - \frac{s_+}{s}, 1 - \frac{r_+}{r}  \Big\}
\ee
holds (see \eqref{piccof}), then Proposition \ref{main}
(with $s'\rightsquigarrow s_+, r'\rightsquigarrow r_+,N_0\rightsquigarrow N_0^* $ defined in \eqref{N0large})
implies that the Hamiltonian flow
$ \Phi^1_F : D(s_+,r_+) \to D(s,r) $
is well defined, analytic and symplectic.
We transform the Hamiltonian $ H $ in \eqref{def:R}, obtaining
\begin{eqnarray*}
H^+ := e^{{\rm ad}_F} H & \stackrel{\eqref{expLie}} =
& H + {\rm ad}_F (H) + \sum_{j \geq 2} \frac{1}{j!}{\rm ad}_F^j (H)  \\
& \stackrel{\eqref{def:R}} = & {\cal N} + P^{\leq 2}_K + ( P - P^{\leq 2}_K ) + 
{\rm ad}_F   {\cal N} + {\rm ad}_{F} P +  \sum_{j \geq 2} \frac{1}{j!}{\rm ad}_F^j (H)   \\
& \stackrel{\eqref{homo}} = &
{\cal N} + [P^{\leq 2}_K] + 
   P - P^{\leq 2}_K +
{\rm ad}_{F} P +
\sum_{j \geq 2} \frac{1}{j!}{\rm ad}_F^j (H) :=
 {\cal N}^+ + P^+
\end{eqnarray*}
with new normal form
$$
{\cal N}^+  := {\cal N} + \hat {\cal N} \, , \quad
 \hat {\cal N} := [P^{\leq 2}_K] = \hat e + \hat \o \cdot y + \hat \O z \cdot \bar z
$$
\begin{equation}\label{triciclo}
\hat \o_i:= \partial_{y_i |\, y=0, z = \bar z = 0 } \langle P \rangle \, ,
\ i=1,\dots n\,,\ \   \hat \O := ({\hat \Omega}_j )_{j \in \Z \setminus {\cal I}} \, , \
{\hat \Omega}_j := [P]_j :=
\partial^2_{z_j \bar z_j\,|\, y=0, z = \bar z = 0 } \langle P \rangle
\end{equation}
(the $ \langle \ \rangle $ denotes the average with respect to the angles $ x $) and new perturbation
\be\label{newNP}
P^+ :=  P - P^{\leq 2}_K
+ {\rm ad}_{F} P^{\leq 2} + {\rm ad}_{F} P^{\geq 3}  + \sum_{j \geq 2} \frac{1}{j!}{\rm ad}_F^j (H)
\ee
having decomposed $ P = P^{\leq 2} + P^{\geq 3} $ with $  P^{\geq 3} := \sum_{h \geq 3} P^{(h)} $, see \eqref{phom}.

\subsubsection{The new normal form $ {\cal N}^+ $}\label{sec:normalf}

\begin{lemma}\label{frodo}
Let $ P \in {\mathcal Q}^T_{s,r}(N_0, \teta,\mu)$  with  $ \cc <  \teta , \mu < \CC $, $ N_0 \geq 9 $.
Then
\be\label{supomega}
|\hat \o|, | \hat \O|_\infty \leq 2 \| P^{ (2)} \|_{s, r,N_0, \teta,\mu}^T
\ee
and there exist $ \hat a_\pm \in\mathbb{R} $
satisfying
$$
|\hat a_\pm |  \leq 2 \| P^{ (2)} \|_{s, r, N_0, \teta, \mu }^T
$$
such that
\begin{equation}\label{baggins}
|\hat\O_j-\hat a_{\mathtt s(j)}|\leq \frac{40}{|j|} \| P^{ (2)} \|_{s, r,N_0,\teta,\mu}^T  \, , \quad
\forall\, |j| \geq \CC(N_0+1) \, .  
\end{equation}
Moreover
$ |\hat \o|^{\rm lip} $,  $ |\hat \O|^{\rm lip}_\infty \lessdot |X_{P^{(2)}}|^{\rm lip}_{s,r} $.
\end{lemma}
\noindent
Lemma \ref{frodo} is based on the following elementary Lemma, whose proof is postponed.

\begin{lemma}\label{weleda}
Suppose that,
$ \forall N \geq N_0 \geq 9 $,  $ j \geq  \teta N $,
\begin{equation}\label{purelan}
\O_j = a_N + b_{N,j} N^{-1}  \ \ {\rm with }\ \
a_N, b_{N,j} \in\mathbb{R}\,,\ \ |a_{N}|\leq c_1\,, \ |b_{N,j}|\leq c_1\,,
\end{equation}
for some  $ c_1 > 0 $ (independent of $ j $).  Then there exists $ a \in \mathbb{R} $,
satisfying $ |a| \leq c_1 $, such that
\be\label{corbezzolo}
| \O_j - a | \leq \frac{20 c_1}{|j|} \, , \quad
\forall\, |j| \geq \CC(N_0+1) \, .
\ee
\end{lemma}

\noindent
\textsc{proof of Lemma \ref{frodo}.}
The estimate on $\hat\o$ is trivial.  Regarding $ \hat\O $ we set (recall \eqref{phom}, \eqref{bische})
$$
P^{(2)}_0:=\Pi_{k=0} \Pi_{ |\a|=|\b|=1}\Pi^{(2)} P = \sum_j [P]_j z_j \bar z_j
$$
since, by the momentum conservation \eqref{momento2}, all the monomials in $P^{(2)}_0$  have
$ \a = \b = e_j $. 
Note that $[P]_j$ is defined in \eqref{triciclo}. By  Lemma \ref{diago}
\be\label{piro}
| [P]_j | \leq
 \|X_{P^{(2)}_0}\|_{r} \stackrel{\eqref{unobisbis_aa}} \leq  \| P^{(2)}_0 \|_{r}^T
  {\stackrel{\eqref{bische}}\leq}  \| P^{(2)} \|_{s, r}^T \, . 
\ee
We now prove \eqref{baggins} for  $ j > 0 $ (the case $ j < 0 $ is similar).
Since $ P^{(2)}_0 \in  \mathcal{Q}^T_{r} (N, \teta,\mu) $,
for all $  N \geq N_0 $,  we may write
$\Pi_{N, \teta,\m} P^{(2)}_0 = {\tilde P}^{(2)}_{0,N} + N^{-1} {\hat P}^{(2)}_{0,N} $ with
$$
{\tilde P}^{(2)}_{0,N}:= \sum_{j>\teta N} \tilde P_j z_j\bar z_j \in \mathcal{T}_{r} (N, \teta,\mu)
\,, \quad  {\hat P}^{(2)}_{0,N}:= \sum_{j>\teta N} \hat P_j z_j\bar z_j
$$
and
\be\label{bou3}
\| {X}_{P^{(2)}_0}\|_{r}\, ,
\| {X}_{{\tilde P}^{(2)}_{0,N} }\|_{r}\, , \, \| {X}_{{\hat P}^{(2)}_{0,N} }\|_{r} \leq 2  \| P^{(2)}_0 \|_{r}^T\leq
2 \| P^{(2)} \|_{s, r}^T \, .
\ee
For $|j|>\teta N$, since all the quadratic forms in \eqref{bou3} are diagonal, we have
$$
\hat \Omega_j = [P]_j= \tilde P_j +N^{-1} \hat P_j:=
a_{N,+} + N^{-1} b_{N,j}
$$
where
$ a_{N,+} := {\tilde P} _j $ is   independent  of $ j > 0 $ because
$ {\tilde P}^{(2)}_{0,N} \in \mathcal{T}_{r}(N, \teta,\mu) $ (see \eqref{marco}). Applying Lemma
\ref{diago} to ${\tilde P}^{(2)}_{0,N}$ and ${\hat P}^{(2)}_{0,N}$, we obtain
$$
|a_{N,+}|
\leq  \| {X}_{{\tilde P}^{(2)}_{0,N}}\|_{s, r} \stackrel{\eqref{bou3}}   \leq  2  \| P^{(2)} \|_{s, r}^T \, , \quad
|b_{N,j}| = | \hat P_j  |  \leq  \| {X}_{{\hat P}^{(2)}_{0,N}}\|_{r} \stackrel{\eqref{bou3}} \leq 2 \| P^{(2)} \|_{r}^T \, .
$$
Hence the assumptions of
 Lemma \ref{weleda} are satisfied with $ c_1 = 2 \|P^{(2)} \|^T_{s, r} $ and  \eqref{baggins}  follows.

The final Lipschitz estimate 
is standard, see e.g. \cite{BB10}, \cite{Po2}.
 \rule{2mm}{2mm}

\medskip

\begin{pfn}\textsc{of Lemma \ref{weleda}.}
For all $ N_1 > N \geq N_0 $, $ j \geq  \teta N_1 $  we get,  by \eqref{purelan},
\begin{equation}\label{an1n2}
|a_{N} - a_{N_1}| = |  b_{N_1,j} N_1^{-1} - b_{N,j} N^{-1}  |
\leq 2 c_1 N^{-1}\,.
\end{equation}
Therefore $ a_N $ is a Cauchy sequence. Let $ a := \lim_{N\to + \infty} a_N $ be its limit.
Since $|a_N| \leq c_1 $ we have $ |a| \leq c_1 $.
Moreover, letting $ N_1 \to + \infty $ in  \eqref{an1n2}, we derive
$ | a - a_N| \leq 2 c_1 N^{-1} $, $ \forall  N \geq N_0 $,
and, using also  \eqref{purelan},
\be\label{disuf}
| \O_j - a | \leq |\O_j - a_N| + |a_N -a | \leq 3 c_1 N^{-1}, \quad
\forall\, N \geq N_0 \, , \ j \geq \CC N \, . 
\ee
For all $ j \geq  \CC (N_0+1) $ let $ N :=[  j /\CC ]$ (where $[\cdot]$ denotes the integer part).
Since $ N \geq N_0 $, $j \geq \CC N $,
$$
|\O_j-a| \stackrel{\eqref{disuf}}  \leq \frac{3 c_1}{[j/\CC]} \leq \frac{3 c_1}{(j/\CC)-1} \leq
\frac{\treCC c_1 }{j}
\Big( 1 + \frac{1}{N_0} \Big) \leq \frac{20 c_1 }{j}
$$
for all $  j \geq  \CC (N_0 +1) $.
\end{pfn}

\subsubsection{The new perturbation $ P^+ $} \label{sec:P+}

We introduce, for  $ h = 0,1,2 $,
\begin{eqnarray}\label{xhx}
&&\varepsilon^{(h)} :=
\g^{-1}\max\Big\{ \| P^{(h)} \|_{s,r, N_0, \teta,\mu}^T \,,\
| {X}_{P^{(h)}} |_{s,r}^\l
\Big\}\, ,
\qquad
\bar\varepsilon := \sum_{h=0}^2 \varepsilon^{(h)} \, , \quad
\\
&&\Theta :=
\g^{-1}\max\Big\{ \| P \|_{s,r, N_0, \teta,\mu}^T \,,\
| {X}_{P}|_{s,r}^\l
\Big\}\, ,
\nonumber
\end{eqnarray}
($\l$ defined in \eqref{KAMconditionT}) and the corresponding quantities 
for $ P^+ $ with indices $ r_+ $, $ s_+, N_0^+,\teta_+,\mu_+ $.
The $ P^{(h)}$ denote the homogeneous components of $ P $
of degree $ h $ (see \eqref{phom}). 

\begin{proposition} {\bf (KAM step)} \label{KAMstep}
Suppose  $ (s,r,N_0, \teta,\mu) $,
$ (s_+,r_+,N_0^+,\teta_+,\mu_+)$
satisfy $s/2\leq s_+<s$, $r/2\leq  r_+  < r $,
\be\label{lap1}
N_0^+>\max\{N_0^*,\bar N    \}  \
\mbox{(recall\  \eqref{N0large},  \eqref{stoppa})}
 \, ,  \quad 2(N_0^+)^{-b}\ln^2 N_0^+ \leq b(s-s_+) \, ,
\ee
\begin{equation}\label{Giobbe+}
 \kappa  (N_0^+)^{b-L} \ln N_0^+  \leq \mu - \mu_+, \quad
 (\CC  + \kappa) (N_0^+)^{L-1} \ln N_0^+ \leq   \teta_+  -  \teta \, .
 \end{equation}

Assume that
\be\label{piccoloKAM}
\bar\e K^{ \bar \tau} \d_+^{-1} \leq c \ {\rm small \ enough}  \, , \quad \Theta \leq 1 \, ,
\ee
where ${\bar\t} := 2 \t + n + 1$ and $ \d_+ $ is defined in \eqref{piccoF}.
Suppose also that \eqref{Omasin} holds  for  $|j|\geq \teta N_0^*$.

Then,  for all $ \xi \in {\cal O} $ satisfying \eqref{foederisarca},\eqref{s+s-},
denoting by $ F := F^{\leq 2}_K $ the solution of the homological equation \eqref{homo},
the Hamiltonian flow $ \Phi^1_F : D(s_+,r_+) \to D(s,r) $,
and the transformed Hamiltonian
$$
H^+ := H \circ \Phi^1_F = e^{{\rm ad}_F} H = {\cal N}_+ + P_+
$$
satisfies
\begin{eqnarray}
\varepsilon^{(0)}_+  & \lessdot &
\d_+^{- 2} K^{2 {\bar\t}} \bar\varepsilon^2 +   \varepsilon^{(0)}  e^{- (s - s_+) K} \nonumber \\
\varepsilon^{(1)}_+  & \lessdot &
 \d_+^{- 2 } K^{2 {\bar\t}}  \big(  \varepsilon^{(0)} +   \bar\varepsilon^2 \big)  +   \varepsilon^{(1)}  e^{- (s - s_+) K}
 \nonumber \\
 \varepsilon^{(2)}_+  & \lessdot &   \d_+^{- 2 } K^{2 {\bar\t}}
  \big( \varepsilon^{(0)}+ \varepsilon^{(1)} +  \bar\varepsilon^2   \big) +
   \varepsilon^{(2)}  e^{- (s - s_+) K }
 \label{kamst}
 \end{eqnarray}
 \be\label{ultimax+}
\Theta_+ \, \leq  \Theta ( 1 + C \d_+^{- 2 } K^{2 {\bar\t}} \bar\e )  \, .
\ee
\end{proposition}

We focus on the quasi-T\"oplitz estimates,
the Lipschitz ones follow formally in the same way.
The proof splits in several lemmas
where we analyze each term of $P^+ $ in \eqref{newNP}.
We note first that
\be\label{Pe2}
\| P^{\leq 2}_K \|_{s,r, N_0, \teta, \mu}^T \stackrel{\eqref{fKT}}
\leq \| P^{\leq 2} \|_{s,r, N_0, \teta, \mu}^T \stackrel{\eqref{Pleq2}, \eqref{xhx}} \leq \g \bar\varepsilon \, .
\ee
Moreover, the solution $ F = F^{(0)} +  F^{(1)} +  F^{(2)} $
of the homological equation \eqref{homo} (for brevity $ F^{(h)} \equiv F^{(h)}_K  $ and $ F \equiv F^{\leq 2}_K $)
satisfies,
by \eqref{Fijstimeh} (with $ N_0^*  $ defined in \eqref{N0large}),  \eqref{fKT},  \eqref{xhx},  
\be\label{PF}
\| F^{(h)} \|_{s,r,N_0^*, \teta,\mu}^T \lessdot K^{{\bar\t}} \varepsilon^{(h)} \, , \ h = 0,1,2, \
\quad \| F \|_{s,r, N_0^*, \teta,\mu}^T \lessdot K^{{\bar\t}} \bar\varepsilon \, .
\ee
Hence \eqref{piccoloKAM} and \eqref{PF} imply condition \eqref{piccoF} and therefore
$ \Phi^1_F : D(s_+, r_+ ) \to D(s,r) $ is well defined. We now estimate the terms of the new
perturbation $ P^+ $ in  \eqref{newNP}.

\begin{lemma}\label{2pez}
$$
\Big\|  {\rm ad}_F ( P^{\leq 2} ) \Big\|_{s_+, r_+, N_0^+, \teta_+,\mu_+}^T  +
\Big\| \sum_{j \geq 2} \frac{1}{j!}{\rm ad}_F^j (H) \Big\|_{s_+, r_+, N_0^+, \teta_+,\mu_+}^T \lessdot
\d_+^{-2} \gamma  K^{2 {\bar\t}} \bar\varepsilon^2 \, .
$$
\end{lemma}

\begin{pf}
We have
\begin{eqnarray*}
\sum_{j \geq 2} \frac{1}{j!}{\rm ad}_F^j (H) & = & \sum_{j \geq 2} \frac{1}{j!}{\rm ad}_F^j ({\cal N} + P) =
\sum_{j \geq 2} \frac{1}{j!}{\rm ad}_F^{j-1} ( {\rm ad}_F {\cal N}) +  \sum_{j \geq 2} \frac{1}{j!}{\rm ad}_F^j (P)  \\
& \stackrel{\eqref{homo} } = &
\sum_{j \geq 2} \frac{1}{j!}{\rm ad}_F^{j-1} ( [P^{\leq 2}_K] - P^{\leq 2}_K)
+  \sum_{j \geq 2} \frac{1}{j!}{\rm ad}_F^j (P) \, .
\end{eqnarray*}
By
\eqref{lap1}, \eqref{Giobbe+} and \eqref{piccoF}
we can apply Proposition \ref{main}
with
$ N_0,N_0',s',r',\teta',\mu',\d \rightsquigarrow $ $ N_0^*,N_0^+,s_+, $ $ r_+,\teta_+,\mu_+,\d_+ $.
We get (recall $ N_0^* \geq N_0 $)
\begin{eqnarray}\label{uno}
\Big\| \sum_{j \geq 2} \frac{1}{j!}{\rm ad}_F^j (P) \Big\|_{s_+,r_+, N_0^+, \teta_+, \mu_+}^T
& \stackrel{\eqref{gPhif12}, \eqref{inscatola}}\lessdot &
 \Big( \d_+^{-1} \| F \|_{s,r, N_0^*,  \teta, \mu}^T \Big)^2 \| P \|_{s,r, N_0, \teta,\mu}^T \nonumber \\
& \stackrel{\eqref{PF}, \eqref{xhx} } \lessdot  &  \d_+^{-2} K^{2 {\bar\t}} \bar\varepsilon^2 \g  \, \Theta
\end{eqnarray}
and, similarly,
\begin{eqnarray}
\Big\| \sum_{j \geq 2} \frac{1}{j!}{\rm ad}_F^{j-1} (P^{\leq 2}_K) \Big\|_{s_+,r_+, N_0^+, \teta_+, \mu_+}^T
& = &
\Big\| \sum_{j \geq 1} \frac{1}{(j+1)!}{\rm ad}_F^{j} (P^{\leq 2}_K) \Big\|_{s_+,r_+, N_0^+, \teta_+, \mu_+}^T \nonumber \\
& \stackrel{\eqref{gPhif12}} \lessdot &   \d_+^{-1} \| F \|_{s,r, N_0^*,  \teta, \mu}^T  \| P^{\leq 2}_K \|_{s,r, N_0,  \teta, \mu}^T  \nonumber \\
& \stackrel{\eqref{PF}, \eqref{Pe2}}  \lessdot &
   \d_+^{-1}  K^{\bar\t} \g \bar\varepsilon^2  \, .  \label{due}
\end{eqnarray}
Finally, by Proposition \ref{festa}, applied with
\be\label{pampersprogressi}
N_0,N_1,s_1,r_1,\teta_1,\mu_1,\d
\rightsquigarrow
N_0^*,N_0^+,s_+,r_+,\teta_+,\mu_+,\d_+\,,
\ee
we get
\begin{eqnarray}
\Big\|  {\rm ad}_F ( P^{\leq 2} ) \Big\|_{s_+, r_+, N_0^+, \teta_+, \mu_+}^T \nonumber
& \stackrel{\eqref{poisbound2}} \lessdot &
\d_+^{-1} \| F \|_{s,r, N_0^*,  \teta, \mu}^T \| P^{\leq 2} \|_{s,r, N_0,  \teta, \mu}^T \\
& \stackrel{\eqref{PF}, \eqref{Pe2}}\lessdot &  \d_+^{-1} K^{\bar\t} \gamma \, \bar\varepsilon^2  \, . \label{tre}
\end{eqnarray}

The bounds \eqref{uno}, \eqref{due},  \eqref{tre},  and $ \Theta \leq 1 $ (see \eqref{piccoloKAM}),
prove the lemma.
\end{pf}

\begin{lemma}
\eqref{ultimax+} holds.
\end{lemma}

\begin{pf}
By Proposition \ref{festa} (applied with \eqref{pampersprogressi}) we have
\begin{eqnarray}
\Big\|  {\rm ad}_F ( P^{\geq 3} ) \Big\|_{s_+, r_+, N_0^+, \teta_+, \mu_+}^T \nonumber
& \lessdot &
 \d_+^{-1} \| F \|_{s,r, N_0^*,  \teta, \mu}^T \| P^{\geq 3} \|_{s,r, N_0,  \teta, \mu}^T \\
& \stackrel{  \eqref{PF}, \eqref{fhle2T}, \eqref{xhx}}\lessdot &
 \d_+^{-1} K^{\bar\t} \gamma \,  \bar\varepsilon \, \Theta \, , \label{tre0}
\end{eqnarray}
and \eqref{ultimax+} follows by  \eqref{newNP},  \eqref{fhle2T}, \eqref{inscatola},
 \eqref{xhx}
\eqref{tre0},
Lemma \ref{2pez} and $ \bar\varepsilon \leq  3 \Theta $ (which follows by \eqref{xhx} and \eqref{fhT}).
\end{pf}

We now consider $  P^{(h)}_+ $, $ h = 0,1, 2 $.
The term $ {\rm ad}_F P^{\geq 3}  $ in \eqref{newNP}
does not contribute to  $ P^{(0)}_+ $.
On the contrary, its contribution to $ P^{(1)}_+ $ is
\be\label{da1}
\{  F^{(0)} , P^{(3)} \}
\ee
and to $ P^{(2)}_+ $ is
\be\label{da2}
\{  F^{(1)} , P^{(3)} \} +  \{  F^{(0)} , P^{(4)} \} \, .
\ee
\begin{lemma}\label{contrib12}
$ \|  \{  F^{(0)} , P^{(3)} \}\|_{s_+,r_+,N_0^+,\teta_+,\mu_+}^T \lessdot \d_+^{-1} \g K^{\bar\t} \varepsilon^{(0)} \Theta $
and
$$ \Big\|  \{  F^{(1)} , P^{(3)} \} +  \{  F^{(0)} , P^{(4)} \}
\Big\|_{s_+,r_+,N_0^+,\teta_+,\mu_+}^T \lessdot \d_+^{-1}  K^{\bar\t}
\g ( \varepsilon^{(0)} + \varepsilon^{(1)} ) \Theta \, .
$$
\end{lemma}

\begin{pf}
By \eqref{poisbound2} (applied with \eqref{pampersprogressi}),
\eqref{PF}, \eqref{xhx} and \eqref{fhT}.
\end{pf}

The contribution of $ P - P^{\leq 2}_K  $ in  \eqref{newNP} to  $ P^{(h)}_+ $, $ h = 0,1,2 $,
is $ P^{(h)}_{> K} $.

\begin{lemma}\label{lem:ultra}
$ \| P^{(h)}_{> K} \|_{s_+, r_+, N_0^+, \teta_+, \mu_+}^T \leq 2 e^{- K (s - s_+)} \g \varepsilon^{(h)} $
\end{lemma}

\begin{pf}
By  \eqref{smoothT} and \eqref{xhx}.
\end{pf}

\textsc{Proof of Proposition \ref{KAMstep} completed.}
Finally,  \eqref{kamst} follows by
\eqref{newNP}, Lemmata \ref{2pez}, \ref{contrib12} (and \eqref{da1}-\eqref{da2}), Lemma \ref{lem:ultra} and
$ \Theta \leq 1 $.

\subsection{KAM iteration}

\begin{lemma}\label{lem:quadr}
Suppose that $ \varepsilon^{(0)}_{i}, \varepsilon^{(1)}_{i}, \varepsilon^{(2)}_{i}\geq 0,$
$ i = 0, \ldots , \nu $,
satisfy
\begin{eqnarray}
\varepsilon^{(0)}_{i+1}  & \leq &
C_* \,  {\mathtt K}^{i}  \, \bar\varepsilon_i^2 +  C_* \varepsilon_i^{(0)}  \, e^{- K_* 2^{i}}  \label{trio} \\
\varepsilon^{(1)}_{i +1}  & \leq &  C_* \, {\mathtt K}^{i} \big( \,  \varepsilon_i^{(0)} +   \bar\varepsilon_i^2 \big)
+   C_*  \varepsilon_i^{(1)}  \, e^{- K_* 2^{i}}
 \nonumber \\
 \varepsilon^{(2)}_{i+1}  & \leq & C_* \, {\mathtt K}^{i}   \big(  \varepsilon_i^{(0)}+ \varepsilon_i^{(1)}  +  \bar\varepsilon_i^2   \big) +
  \,  C_*  \varepsilon_i^{(2)} \, e^{- K_* 2^{i}} \, , \quad  i = 0, \ldots , \nu -1 \, ,  \nonumber
  \end{eqnarray}
    where $ \bar\varepsilon_i := \varepsilon^{(0)}_{i} + \varepsilon^{(1)}_{i} + \varepsilon^{(2)}_{i} $, for some ${\mathtt K}, C_* , K_* > 1 $.
 Then there exist $ \bar \e_\star<1, C_\star  >  0 $, $ \chi \in (1,2) $,
 depending only on  $ {\mathtt K}, C_*, K_* $
 (and not on $\nu$ and satisfying $1\leq  C_\star e^{- K_*}  $),  such that, if
\be\label{statem}
\bar\varepsilon_0 \leq \bar\e_\star
\ \quad \Longrightarrow \quad \bar\varepsilon_i
\leq
C_\star \,
\bar\e_0 \, e^{ -  K_* \chi^i } \, ,
\ \ \forall i = 0, \ldots, \nu \, .
\ee
\end{lemma}

\begin{pf}
Iterating three times \eqref{trio}
we get
\begin{eqnarray}
\bar\e_{j+3}
&\leq&
c_1 C_*^{c_1} {\mathtt K}^{c_1 j}
\Big(
\e_{j+2}^{(0)} + \e_{j+2}^{(1)}
+\bar\e_{j+2}^2+ \bar\e_{j+2}
e^{- K_* 2^{j+2}}
\Big)
\nonumber
\\
&\leq&
c_2 C_*^{c_2} {\mathtt K}^{c_2 j}
\Big(
\e_{j+1}^{(0)}
+\bar\e_{j+1}^2+\bar\e_{j+1}^4+ \bar\e_{j+1}
e^{- K_* 2^{j+1}}
\Big)
\nonumber
\\
&\leq&
c_3 C_*^{c_3} {\mathtt K}^{c_3 j}
\Big(
\bar\e_{j}^2+\bar\e_{j}^8 + \bar\e_{j}
e^{- K_* 2^{j}}
\Big)\,,  \qquad\qquad\qquad
\forall\, 0\leq j\leq \nu-3\,,
\label{elea}
\end{eqnarray}
for suitable constants
$1<c_1<c_2<c_3.$

We first claim that \eqref{statem} holds
with $\chi:=6/5$ for all $i=3 j\leq \nu.$
Setting $a_j:=\bar \e_{3j}$, we prove  that
there exist $C_\star$ large and $\bar\e_\star$
small (as in the statement) such that if $a_0\leq \bar\e_\star$
then
$$
{\bf (S)}_j \qquad
a_j  \leq c_4^{j+1} a_0  e^{ - K_* \tilde\chi^{3j} } \, , \quad  \forall\,
 0\leq j\leq \nu/3
$$
for a suitable  $c_4=c_4({\mathtt K}, C_*, K_* )\geq 1$
large enough
and $\tilde \chi <2^{1/3}$, e.g. $ \tilde \chi := 5/4 $.
We proceed by induction.
The statement $ {\bf (S)}_0 $
is trivial.
Now suppose  $ {\bf (S)}_{j} $ holds true.
Note that $a_j\leq 1$
taking
$ \bar \e_\star \leq \min_{j\geq 0}   e^{K_* \tilde\chi^{3j} } / c_4^{j+1}  $.
 Then $ {\bf (S)}_{j+1} $ follows by
\begin{eqnarray*}
a_{j+1} & = & \bar\e_{3j+3}
\stackrel{\eqref{elea}}\leq
c_3 C_*^{c_3} {\mathtt K}^{3c_3 j}
\Big(
a_j^2+a_j^8 + a_j
e^{- K_* 2^{3j}}
\Big) \stackrel{a_j\leq 1}\leq
2 c_3 C_*^{c_3} {\mathtt K}^{3c_3 j}
\Big(
a_j^2+ a_j
e^{- K_* 2^{3j}}
\Big)
\\
&\stackrel{({\bf S})_j}\leq&
2 c_3 C_*^{c_3} {\mathtt K}^{3c_3 j}
\Big(
(c_4^{j+1} a_0
 e^{ - K_* \tilde\chi^{3j} })^2+
 (c_4^{j+1} a_0
 e^{ - K_* \tilde\chi^{3j} })
e^{- K_* 2^{3j}}   \Big) \leq c_4^{j+2} a_0  e^{ - K_* \tilde\chi^{3j+3} }
\end{eqnarray*}
since
$
4c_3 C_*^{c_3} {\mathtt K}^{3c_3 j}
 (c_4^{j+1} a_0
 e^{ - K_* \tilde\chi^{3j} })
e^{- K_* 2^{3j}}
 \leq
 c_4^{j+2} a_0
 e^{ - K_* \tilde\chi^{3j+3} }
$
taking $c_4$ large enough (use $\tilde \chi<2$)
and
$$
4 c_3 C_*^{c_3} {\mathtt K}^{3c_3 j}
(c_4^{j+1} a_0
 e^{ - K_* \tilde\chi^{3j} })^2
 \leq
 c_4^{j+2} a_0
 e^{ - K_* \tilde\chi^{3j+3} }
$$
taking $a_0\leq \bar\e_\star $ small enough. We have proved inductively $({\bf S})_j$.
Then \eqref{statem} for $i=3j$
follows since $ 6/ 5 =: \chi < \tilde \chi := 5 / 4 $
and taking $C_\star$ large enough.
The cases $i=3j+1$ and $i=3j+2$
follow analogously
noting that $\bar\e_1,\bar\e_2$
can be made small
by \eqref{trio}
taking $\bar\e_\star$ small.
\end{pf}

For $ \nu \in \mathbb{N} $,  we define
\begin{eqnarray}
&\bullet\ &  s_{\nu+1} := s_\nu - s_0 2^{-\nu-2} \searrow  \frac{s_0}{2},
 \quad \ \,
 r_{\nu+1} := r_\nu -
r_0 2^{-\nu-2} \searrow  \frac{r_0}{2} ,  \quad
D_\nu := D(s_\nu,r_\nu) \,,
\nonumber
\\
&\bullet\ &
K_\nu := K_0 4^\nu\,,\quad N_\nu := N_0  2^{\nu \rho}  \ {\rm with} \
 N_0 :=  \hat{c}\g^{-1/3}K_0^{\t +1}\, ,\quad
 \rho := \max \Big\{ 2(\t+1),
 \frac{1}{L-b}, \frac{1}{1-L} \Big\}
 \,,
\nonumber
\\
&\bullet\ &
\mu_{\nu+1} := \mu_\nu - \mu_0 2^{-\nu-2}
\searrow  \frac{\mu_0}{2} , \quad
 \teta_{\nu+1} := \teta_\nu +  \teta_0 2^{-\nu-2}
\nearrow  3 \frac{\teta_0}{2}
\,.
\label{merdina}
\end{eqnarray}

We consider  $ H^0 = {\cal N}_0 + P_0 : D_0 \times  {\cal O}_* \to \C $ 
with $ {\cal N}_0 := e_0 + \om^{(0)} (\xi ) \cdot y + \Omega^{(0)} (\xi ) \cdot z \bar z $.
We  suppose that $\o^{(0)}$ and $\O^{(0)}$ are defined
 on the \textsl{whole}
 $  \R^n$ (using in case the  Kirszbraun extension theorem),
 that $\Omega^{(0)}$ satisfies \eqref{carbonara} and
 $ | \o^{(0)} |^{\rm lip} $,
$ | \Omega^{(0)} |_\infty^{\rm lip}  \leq  M_0 $
on $\R^n$.
Let $ {\cal O}_0 \subseteq \{  \xi \in {\cal O}_*
 \ : \ B_{\g/M_0}(\xi) \subset {\cal O}_* \} $ where $ {\cal O}_* $ is defined in \eqref{Ozero} and
$ B_r (\xi) $ denotes the open ball in $ \R^n $ of center $ \xi $ and radius $ r > 0 $.

\begin{lemma}{\bf (Iterative lemma)}\label{lem:iter}
Let $ H^0 $ be as above and let
 $ \bar\varepsilon_0 $, $ \Theta_0 $
be defined as in \eqref{xhx} for $ P_0 $.
Then there are  $ K_0 > 0  $ large enough, $ \epsilon_0 > 0 $ small enough, such that, if
\be\label{tolkien}
 \bar\e_0, \Theta_0 \leq \epsilon_0 \, ,
\ee
  then
\\[1mm]
$ {\bf (S1)_{\nu}} $  $ \forall 0 \leq i \leq \nu $, there exist
$ \om^{(i)} $, $ \Omega^{(i)} $, $ a^{(i)}_\pm $ defined for all $ \xi \in \R^n $, satisfying
\begin{eqnarray}\label{muflonej}
&&| \o^{(i)} - \o^{(0)} |+\l | \o^{(i)} - \o^{(0)} |^{\rm lip},
 |  \Omega^{(i)} -  \O^{(0)}  |_{\infty}+
\l |  \Omega^{(i)} -  \O^{(0)}  |_{\infty}^{\rm lip} \leq C (1 - 2^{-i}) \g {\bar \e}_0
 \\
 &&| a^{(i)}_\pm  | \leq C (1 - 2^{-i}) \g {\bar \e}_0 \, , \ | \o^{(i)} |^{\rm lip} \,, \ \  | \Omega^{(i)} |_\infty^{\rm lip}  \leq  (2 - 2^{-i}) M_0 \, .
 \label{muflonejlip}
\end{eqnarray}
There exists $ H^i := {\cal N}_i + P_i  \, : \, D_i \times  {\cal O}_i \to \C $
with  $ {\cal N}_i := e_i + \om^{(i)} (\xi ) \cdot y +
 \Omega^{(i)} (\xi ) \cdot z \bar z $
 in normal form, where, for $ i > 0 $,
\begin{eqnarray}\label{fa}
 {\cal O}_{i} & := &
\Big\{\xi \in  {\cal O}_{i-1} \,  : \,
|\o^{(i-1)} (\xi)\cdot k + \O^{(i-1)}(\xi)\cdot l |\geq
(1-2^{-i})\frac{2 \g }{1+|k|^\t}, \forall (k,l) \in {\bf I} \, ,
 \ |k|\leq K_{i-1} \, , \nonumber  \\
& & \ |\o^{(i-1)} (\xi)\cdot k + p  |\geq
(1-2^{-i}) \frac{2 \g^{2/3}}{1+|k|^\t} \, ,  \forall (k,p) \neq (0,0) \, ,  \ |k|\leq K_{i-1} , \ p \in \Z  \Big\} \, .
\end{eqnarray}
Moreover, $\forall\, 1\leq i\leq\nu$, $H^i=H^{i-1}\circ \Phi^i$ where
$\Phi^i:D_i\times{\cal O}_i \to D_{i-1}$
is a (Lipschitz) family (in $\xi \in{\cal O}_i$)
of close-to-the-identity analytic symplectic maps.
Setting, for $ h = 0,1,2 $,
\begin{eqnarray}\label{rosamystica}
&&
\e_i^{(h)} :=
\g^{-1}\max\Big\{ \| P_i^{(h)} \|_{s_i, r_i,  N_i,\teta_i,\mu_i}^T \,,\
| {X}_{P_i^{(h)}} |_{s_i,r_i}^{\l}
\Big\}\,,\qquad
\bar\varepsilon_i := \sum_{h = 0}^2 \varepsilon_i^{(h)} \, ,
\\
\nonumber
&&\Theta_i :=
\g^{-1}\max\Big\{ \| P_i \|_{s_i, r_i,  N_i,\teta_i,\mu_i}^T \,,\
| {X}_{P_i} |_{s_i,r_i}^{\l}
\Big\}\, ,
\end{eqnarray}
$\forall\, 1\leq i\leq \nu$ and $\forall\, \xi \in \R^n$
\begin{eqnarray}\label{frecoj}
&&| {\o}^{(i)}(\xi)
- {\o}^{(i-1)}(\xi)|\,, \
|  {\Om}^{(i)}(\xi)
- {\Om}^{(i-1)}(\xi) |_\infty \, ,
\ |a_\pm^{(i)}(\xi)
- a_\pm^{(i-1)}(\xi)|
 \leq 2 \g \bar\e_{i-1}  \, ,
 \nonumber
 \\
&&
| {\Om}_j^{(i)}(\xi) -  a^{(i)}_{\mathtt s(j)}(\xi)
- {\Om}_j^{(i-1)}(\xi) +  a^{(i-1)}_{\mathtt s(j)}(\xi)|
  \leq 40 \g \frac{\bar\e_{i-1} }{|j|} \, , \quad \forall |j| \geq \CC (N_{i-1}+1)  \, .
\end{eqnarray}
$ {\bf (S2)_{\nu}} $  $ \forall 0 \leq i \leq \nu -1 $,
the $ \varepsilon_i^{(0)},   \varepsilon_i^{(1)},  \varepsilon_i^{(2)}  $ satisfy  (\ref{trio}) with
$ \mathtt K = 4^{2\bar \t +1} $, $ \bar \t := 2 \t + n + 1 $, $ C_*= 4 K_0^{2 {\bar\t}}$, $  K_*=  s_0 K_0 / 4 $.
\\[1mm]
$ {\bf (S3)_{\nu}} $
$ \forall 0 \leq i \leq \nu  $, we have $ \bar\varepsilon_i
\leq  C_\star \bar\e_0  e^{- K_* \chi^{i}} $ and
$ \Theta_i \leq 2  \Theta_0 $ (recall that $ C_\star e^{- K_*} \geq 1 $, see Lemma \ref{lem:quadr}).
\end{lemma}

\begin{pf}
The statement $ {\bf (S1)}_0  $ follows by  the hypotheses setting $ a^{(0)}_\pm (\xi) := 0 $, $ \forall \xi \in \R^n $.
$ {\bf (S2)}_0  $ is empty. $ {\bf (S3)}_0  $ is trivial. We then proceed by induction.
\\[1mm]
$ {\bf (S1)}_{\nu+1} $.
We denote $ \hat\o^{(\nu)} := \nabla_y\langle P_\nu (\xi)\rangle\vert_{y=0,z=\bar z=0} $ and
$  \hat\Om^{(\nu)}_j(\xi) := \partial^2_{z_j \bar z_j\,|\, y=0, z = \bar z = 0 } \langle P_\nu (\xi) \rangle $,
see \eqref{triciclo}, for all $ \xi \in {\cal O}_\nu $ if $ \nu \geq 1$ and $  \xi \in {\cal O}_*  $ (see \eqref{Ozero}) if $ \nu = 0 $.
By Lemma \ref{frodo} and \eqref{rosamystica}
there exist constants $ \hat a^{(\nu)}_\pm(\xi) \in \R $ such that
\be\label{alfonsina}
|  \hat\om^{(\nu)}(\xi) |\,, |  \hat\Om^{(\nu)}(\xi) |_\infty \, , \, |\hat a^{(\nu)}_\pm(\xi)|
 \leq 2 \g \bar\e_\nu  \, ,  \ \
 | {\hat \Om}_j^{(\nu)}(\xi) - \hat a^{(\nu)}_{\mathtt s(j)}(\xi) |
  \leq 40 \g \frac{\bar\e_\nu }{|j|} \, , \   \forall |j| \geq \CC (N_{\nu}+1)  \, ,
\ee
uniformly in $ \xi \in{\cal O}_\nu $ (resp. $ {\cal O}_* $ if $ \nu = 0 $), and
\be\label{Caprese}
| \hat\om^{(\nu)} |^{\rm lip} \, , \ | \hat\Om^{(\nu)} |_\infty^{\rm lip} \leq C {\bar \e}_\nu \, .
\ee
Let
\begin{equation}\label{salamina}
\eta_0 := \l = \g/M_0 \, , \quad
\eta_\nu:=\g / (2^{\nu+3}M_0 K_{\nu-1}^{\t+1})\,, \ \nu \geq 1 \, .
\end{equation}
We claim that, for $ \nu \geq 1 $, the $ \eta_\nu$-neighborhood  of ${\cal O}_{\nu+1}$
\begin{equation}\label{salamona}
  \tilde{\cal O}_{\nu+1}:=
   \bigcup_{\xi\in{\cal O}_{\nu+1}} \Big\{ \tilde\xi\in \R^n \ : \
  \tilde \xi=\xi+ \hat \xi\,,
   \  |\hat \xi|< \eta_\nu \Big\}
    \ \subseteq  \ {\cal O}_\nu \, .
\end{equation}
Note that  the definitions of $ {\cal O}_0 $, $ {\cal O}_1 $ in \eqref{fa}, and
 \eqref{salamina} imply $ \tilde{\cal O}_{1} \subset {\cal O}_* $.
Recalling \eqref{fa}, we have to prove that for $ \nu \geq 1 $,
for every $ \tilde \xi=\xi+ \hat \xi $, $ \xi \in{\cal O}_{\nu+1}$, $ |\hat \xi|\leq \eta_\nu $,  we have
\begin{equation}\label{termopili}
|\o^{(\nu-1)} (\tilde\xi)\cdot k + \O^{(\nu-1)}(\tilde\xi)\cdot l |\geq
(1-2^{-\nu})\frac{2 \g }{1+|k|^\t}, \quad  \forall (k,l) \in {\bf I} \, , \ |k|\leq K_{\nu-1} \, ,
\end{equation}
and the analogous estimate for
$|\o^{(\nu-1)} (\tilde\xi)\cdot k + p  |$.
By  the expression \eqref{defnu+1} (at the previous step) for $  \o^{(\nu)} $,  $ \O^{(\nu)} $, and since $ \chi_{\nu-1} \in [0,1] $,
we get
\begin{eqnarray*}
&& |\o^{(\nu-1)} (\tilde\xi)\cdot k + \O^{(\nu-1)}(\tilde\xi)\cdot l |
\geq
|\o^{(\nu)} (\tilde\xi)\cdot k + \O^{(\nu)}(\tilde\xi)\cdot l |
-|\chi_{\nu-1}(\tilde\xi)|
|\hat\o^{(\nu-1)} (\tilde\xi)\cdot k +
 \hat\O^{(\nu-1)}(\tilde\xi)\cdot l |
 \\
 &&\stackrel{\eqref{alfonsina}}\geq |\o^{(\nu)} (\xi)\cdot k + \O^{(\nu)}(\xi)\cdot l |
- \Big|(\o^{(\nu)} (\tilde\xi)-\o^{(\nu)} (\xi))\cdot k
+ (\O^{(\nu)} (\tilde\xi)-\O^{(\nu)} (\xi))\cdot l \Big|
-2\g \bar\e_{\nu-1} (K_{\nu-1}+2)
\\
&&
\stackrel{\xi\in{\cal O}_{\nu+1}, \eqref{fa}, (S1)_\nu}{\geq}
(1-2^{-\nu-1})\frac{2\g}{1+|k|^\t}
-(K_{\nu-1}+2)2 M_0 \eta_\nu - 2\g \bar\e_{\nu-1} (K_{\nu-1}+2)
\\
&&
\stackrel{\eqref{salamina}, (S3)_\nu}{\geq}
(1-2^{-\nu})\frac{2\g}{1+|k|^\t}
\end{eqnarray*}
taking $\epsilon_0$ small enough, and \eqref{termopili} follows.
The  estimate for $|\o^{(\nu-1)} (\tilde\xi)\cdot k + p  |$
follows similarly.

We define a smooth cut-off function $ \chi_\nu : \R^n \to [0,1] $ which takes value $ 1 $ on
 $ {\cal O}_{\nu+1}  $  and  value $ 0 $ outside $ {\tilde {\cal O}}_{\nu+1}  $.
 Thanks to \eqref{salamona} and recalling \eqref{salamina} we can construct $ \chi_\nu $, $ \nu \geq  0 $,
in such a way that
\be\label{estchinu}
 | \chi_\nu |^{\rm lip} \lessdot \g^{-1}M_0 2^\nu K_{\nu-1}^{\t + 1}  
\ee
where  $ K_{-1} := 1 $. We extend $ \hat\om^{(\nu)}, \hat\Om^{(\nu)},
\hat a^{(\nu)}_\pm$ to zero outside ${\cal O}_\nu$ for $ \nu \geq 1 $ and, for $ \nu = 0 $
outside ${\cal O}_\star $.  Then we define on the {\sl whole} $\R^n$
\be\label{defnu+1}
 \o^{(\nu+1)} := \o^{(\nu)} +  \chi_\nu  {\hat \o}^{(\nu)}
 \,,\quad  \Om^{(\nu+1)} := \Om^{(\nu)} +  \chi_\nu  {\hat \O}^{(\nu)}   \,,
 \quad a^{(\nu+1)}_{\pm} := a^{(\nu)}_\pm +  \chi_\nu {\hat a^{(\nu)}}_\pm \, .
 \ee
By $ \eqref{estchinu}, \eqref{Caprese},  \eqref{alfonsina}$, we get
$$
| \o^{(\nu+1)} -\o^{(\nu)}|^{\rm lip}  \leq
 | \chi_\nu |^{\rm lip} | \hat \o^{(\nu)} | +
| \chi_\nu || \hat \o^{(\nu)} |^{\rm lip} \leq
  CK_{\nu-1}^{\t+1} M_0 {\bar \e}_\nu  + C {\bar \e}_\nu \leq   2^{-\nu-1} M_0
$$
by $({\bf S3})_\nu$ and $\bar \e_0$ small enough.
Similarly for $| \O^{(\nu+1)} - \O^{(\nu)}|^{\rm lip}_\infty.$
Recalling also \eqref{alfonsina}, we get
\eqref{muflonej} and \eqref{muflonejlip}
with $ i = \nu + 1 $.
Moreover
 \eqref{alfonsina}-\eqref{defnu+1}  imply \eqref{frecoj} for $ i = \nu + 1 $
  and $\forall |j|>  \CC( N_\nu+1)$.

We wish to apply the KAM step Proposition \ref{KAMstep}
with $\mathcal N= \mathcal N_\nu \,, P=P_\nu, N_0= N_\nu,\theta=\theta_\nu\dots$ and
$N^+_0= N_{\nu+1}, \theta_+=\theta_{\nu+1},\dots $   Our definitions in \eqref{merdina} (and $ \t > 1/b $) imply that
the conditions\footnote{For example the first inequality in \eqref{lap1}
reads $ N_{\nu+1} \geq \max \{ N_\nu, \hat{c} \g^{-1/3} K_\nu^{\t+1}, \bar N \} $.}
\eqref{lap1}-\eqref{Giobbe+} are satisfied, for all
$ \nu \in \N $, taking $ K_0  $ large enough.
Moreover, since
\be\label{deltanu}
\d^+=
\d_{\nu+1} :=
\min \Big\{  1 - \frac{s_{\nu+1}}{ s_{\nu}} ,
 1 - \frac{r_{\nu+1}}{r_{\nu}} \Big\}
 \quad
{\rm so \ that } \quad 2^{-\nu-2} \leq \d_{\nu+1}  \leq 2^{-\nu-1} \, ,
\ee
and $ (S3)_\nu $  the condition \eqref{piccoloKAM}
 is satisfied, for $ \bar\e_0 \leq \epsilon_0 $ small enough,  $ \forall \nu \in \N $.
By \eqref{frecoj}, the condition \eqref{Omasin}
holds for $ |j | \geq \teta_\nu N_\nu $,
and  \eqref{foederisarca} and \eqref{s+s-} hold  for all $\xi\in {\cal O}_{\nu+1}$ (it is the definition of
${\cal O}_{\nu+1} $, see \eqref{fa}). Hence  Proposition \ref{KAMstep} applies. For all $ \xi \in {\cal O}_{\nu +1} $
the Hamiltonian flow
$ \Phi^{\nu+1} := \Phi^1_{F_\nu} : D_{\nu+1} \times  {\cal O}_{\nu+ 1} \to D_\nu $ and we define
$$
H^{\nu+1} := H^\nu \circ  \Phi^{\nu+1}  = e^{\rm ad_{F_\nu}} H^\nu
= {\cal N}_{\nu + 1}  + P_{\nu+1}  : D_{\nu+1} \times  {\cal O}_{\nu+ 1} \to \C \, .
$$
$ {\bf (S2)}_{\nu+1} $ follows by \eqref{kamst} and \eqref{merdina}.
\\[1mm]
$ {\bf (S3)}_{\nu+1} $. By $ (S2)_\nu $ we can apply Lemma
\ref{lem:quadr} and
 \eqref{statem} implies
$ \bar\varepsilon_{\nu+1} \leq C_\star \bar\e_0 e^{- K_* \chi^{\nu+1}} $.
Moreover,  for $ \epsilon_0 $ small enough,
$
\Theta_{\nu+1}   \stackrel{\eqref{ultimax+}} \leq
\Theta_0  \Pi_{i=0}^{\nu}  \Big( 1 + C \d_{i+1}^{-2} K_i^{2 \bar \tau}
\bar\varepsilon_i \Big)
 \stackrel{\eqref{deltanu}, (S3)_\nu}
\leq 2 \Theta_0 \, .
$
\end{pf}

\noindent
{\bf Proof of the KAM Theorem \ref{thm:IBKAM} completed.}
We apply the iterative Lemma \ref{lem:iter} to the Hamiltonian
$H^0$ in \eqref{accazero} where $\o^{(0)}=\o$ and $\O^{(0)}=\O$ are  defined in \eqref{NormalN}.
We choose
\be\label{defO0}
\mathcal O_0:= \Big\{\xi\in \mathcal O\,:\quad |\o(\xi)\cdot k|\geq \frac{2\g^{2/3}}{1+|k|^n}\,, \ \forall\,
0<|k|< \g^{-1/(7n)} \Big\}
\ee
so that $ {\cal O}_0 \subseteq \{  \xi \in {\cal O}_*
 \ : \ B_{\g/M_0}(\xi) \subset {\cal O}_* \} $, see \eqref{Ozero} and \eqref{M}.
The smallness assumption \eqref{tolkien}
holds by  \eqref{piccolezza0}-\eqref{piccolezza0lip} (use also Lemma \ref{proJ}) and $ \e $ small enough. 
Then the iterative Lemma \ref{lem:iter} applies.
Let us define
$$
\o^\infty :=\lim_{\nu \to \infty} \o^{(\nu)}
\, , \quad \O^\infty := \lim_{\nu \to \infty} \O^{(\nu)} \, , \quad a^\infty_\pm :=\lim_{\nu \to \infty} a^{(\nu)}_\pm \, .
$$
It could happen that $ {\cal O}_{\nu_0} = \emptyset $ for some $ \nu_0 $.
In such a case $ {\cal O}_\infty = \emptyset $ and the  iterative process stops after finitely many steps.
However, we can always set
$ \omega^{(\nu)}  :=  \omega^{(\nu_0)} $, $  \Omega^{(\nu)}  := \Omega^{(\nu_0)} $,
$ a^{(\nu)}_\pm  := a^{(\nu_0)}_\pm $,
$ \forall \nu \geq \nu_0 $, and $ \o^\infty $, $ \O^\infty $, $ a^\infty_\pm $ are always well defined.

The bounds \eqref{muflone} follow by \eqref{muflonej} (with a different constant $C$).
We now prove \eqref{freco}. We consider the case $ j > 0 $.
For all $ \forall \nu \geq 0 $, $ j \geq  6 (N_\nu + 1)  $, we have (recall that $ a^{(0)}_+ = 0 $)
\begin{eqnarray*}
|\Om_j^\infty - \Om_j^{(0)} - a^\infty_+ |
&\leq&
\sum_{0\leq i\leq\nu}
| \Om_j^{(i+1)}  -  a^{(i+1)}_+ - \Om_j^{(i)}  + a^{(i)}_+  |
+
\sum_{i > \nu}
| \Om_j^{(i+1)}   - \Om_j^{(i)}|  + |  a^{(i+1)}_+ - a^{(i)}_+  |
\\
  & \stackrel{\eqref{frecoj}} \leq &
 40 \g  \sum_{0\leq i \leq \nu} \frac{{\bar \e}_i}{j}
 +4 \g \sum_{i > \nu}  \bar\e_i \stackrel{ ({\bf S3})_\nu} \lessdot
 \frac{\bar\e_0 \g}{j}  + \g \sum_{i > \nu}  \bar\e_i \, .
\end{eqnarray*}
Therefore, $ \forall \nu \geq 0$,
$ 6 (N_\nu +1) \leq j < 6 ( N_{\nu+1} + 1)  $,
$$
|\Om_j^\infty - \Om_j^{(0)} - a^\infty_+ |
\lessdot
\frac{\bar\e_0 \g }{j} + \g  \frac{ N_{\nu+1}}{j}
\sum_{i > \nu}  \bar\e_i \stackrel{\eqref{merdina}} \lessdot
\frac{ \bar\e_0 \g}{j} +  \frac{\g }{j}  \g^{-1/3} K_0^{\t+1} 2^{\rho (\nu +1)}
\sum_{i > \nu}  \bar\e_i
$$
and  \eqref{freco} follows by $(S3)_\nu $.

The symplectic transformation $\Phi$ in \eqref{eq:Psi} is defined by
$$
\Phi:=\lim_{\nu\to\infty}\Phi_{00}\circ
\Phi_0\circ \Phi_1 \circ \cdots \circ \Phi^\nu
$$
with $\Phi_{00}$ defined in \eqref{Phi00}.
We now verify that $\Phi$ is defined on
$\mathcal O_\infty $, see \eqref{Cantorinf}.
\begin{lemma}\label{L:new}
$ \mathcal O_\infty \subset \cap_i \mathcal O_i $ (defined in \eqref{fa}).
\end{lemma}
\begin{pf}
We have $ \mathcal O_\infty \subseteq \mathcal O_0  $ by  \eqref{Cantorinf} and \eqref{defO0}.
For $ i \geq 1$, if $ \xi \in   \mathcal O_\infty $ then, for all $ |k | \leq K_i $, 
$|l| \leq  2 $,
\begin{eqnarray*}
&&| \o^{(i)} (\xi)\cdot k + \O^{(i)}(\xi)\cdot l |
\\
&&  \geq  |\o^\infty (\xi)
\cdot k +  \O^{\infty}(\xi)\cdot l |-
|k|\sum_{n\geq i}| \o^{(n+1)}(\xi)-\o^{(n)} (\xi)|
-  2 \sum_{n\geq i}|\O^{(n+1)}(\xi)-\O^{(n)}(\xi)|_\infty
\\
&&\stackrel{\eqref{Cantorinf}, \eqref{frecoj}}{\geq}
\frac{2\g}{1+|k|^{\t}}- K_i 2\g \sum_{n \geq i }\bar \e_n
-4\g\sum_{ n\geq i}\bar \e_n
\geq (1-2^{-i})\frac{2\g}{1+|k|^{\t}}
\end{eqnarray*}
by the
definition of $ K_i $ in \eqref{merdina}, $(S3)_\nu $ and $ \e $ small enough.
The other estimate is analogous.
\end{pf}

Finally $P^\infty_{\leq 2}=0$ (see \eqref{Hnew})
follows by $\bar \e_i\to 0$ as $i\to\infty.$ This concludes the proof of Theorem  \ref{thm:IBKAM}.

\section{Measure estimates: proof of Theorem \ref{thm:measure}}\setcounter{equation}{0} \label{sec:meas}

We have to estimate the measure of
\be\label{santadorotea}
{\cal O} \setminus {\cal O}_\infty = \bigcup_{(k,l) \in \L_0 \cup \L_1 \cup \L_2^+\cup \L_2^-}
{\cal R}_{kl}(\g) \bigcup_{(k,p) \in \Z^{n+1} \setminus \{0\}}  \tilde {\cal R}_{kp}(\g^{2/3})
\bigcup (\mathcal O \setminus \mathcal O_0)
\ee
where
\be\label{stimaRkl}
{\cal R}_{kl}(\g) :=
{\cal R}_{kl}^\t (\g) :=
 \Big\{  \xi \in {\cal O} \, : \,
|\om^\infty (\xi) \cdot k + \Om^\infty (\xi )\cdot l| < \frac{2\g}{1+ |k|^{\tau}} \Big\}
\ee
\be\label{defRkl}
 \tilde {\cal R}_{kp}(\g^{2/3}) := \tilde {\cal R}_{kp}^\t(\g^{2/3}) :=
 \Big\{  \xi \in {\cal O} \, : \,  |\om^\infty (\xi) \cdot k + p| < \frac{2\g^{2/3}}{1+ |k|^{\tau}} \Big\}
\ee
and
\be\label{Lh}
\quad \L_h := \Big\{ (k, l) \in {\bf I} \ ({\rm see} \, \, \eqref{masotti})  
\, , \, |l| = h  \Big\}\, , \  \  h = 0, 1, 2 \, ,
\qquad
\L_2=\L_2^+\cup\L_2^-\,,
\ee
$$
\L_2^+ := \Big\{  (k, l) \in \Lambda_2  \, , \, l = \pm (e_i + e_j) \Big\}\,,
\quad
\L_2^- := \Big\{  (k, l) \in  \Lambda_2   \, , \, l = e_i - e_j  \Big\}\,.
$$
We first consider the most difficult  case $ \L_2^- $.
Setting
$ {\cal R}_{k,i,j} (\g) :=  {\cal R}_{k, e_i - e_j}(\g) $
we show that
\be\label{eccola}
\Big| \bigcup_{(k,l)  \in \Lambda_2^{-}} {\cal R}_{k,l}(\g) \Big| =
\Big| \bigcup_{(k,i,j) \in \, {\mathtt I}} {\cal R}_{k,i,j}(\g) \Big| \lessdot \g^{2/3} \rho^{n-1}
\ee
where
\be\label{unionr}
{\mathtt I} := \Big\{ (k,i,j) \in \Z^n \times (\Z \setminus {\cal I})^2 \, : \,   (k,i,j) \neq  (0, i,  i ) \, , \
\pluto \cdot k + i - j = 0  \Big\} \, .
\ee
Note that the indices in $ {\mathtt I}  $ satisfy 
\be\label{nonvuoto}
||i|-|j|| \leq \kappa \, | k |\quad {\rm and}
\quad k\neq 0
\, .
\ee
Since the matrix $ A $ in \eqref{omegaxi}  is invertible, the bound
\eqref{muflone} implies, for $ \e $ small enough, that
\be\label{inveromi}
\o^\infty : {\cal O} \to \o^\infty({\cal O}) \ \ {\rm is \ \ invertible \ \ and } \ \  |(\o^\infty)^{-1}|^{\rm lip}\leq  2 \| A^{-1} \| \, .
\ee

\begin{lemma}\label{stimaRij}  For $ (k,i,j) \in {\mathtt I }$, $ \eta \in (0,1) $, we have
\be\label{Rkij}
|{\cal R}^\t_{k,i,j} ( \eta ) | \lessdot \frac{\eta \rho^{n-1}}{ 1 + |k|^{\t+1}} \, .
\ee
\end{lemma}

\begin{pf}
By \eqref{muflone}  and \eqref{omegaxi}
$$
\om^\infty (\xi) \cdot k + \Om^\infty_i (\xi) - \Om^\infty_j (\xi)
 =
\om^\infty (\xi) \cdot k +  \sqrt{i^2+\mm} - \sqrt{j^2+\mm} +  r_{i,j}( \xi )
$$
where
\be\label{errors}
 |r_{i,j}(\xi)| = O( \e  \g) \, , \  |r_{i,j}|^{\rm lip} = O( \e ) \, .
\ee
We introduce  the final frequencies $ \zeta :=  \om^\infty (\xi) $ as parameters (see \eqref{inveromi}), and we consider
$$
 f_{k,i,j}(\zeta ) := \zeta \cdot k +  \sqrt{i^2+\mm} - \sqrt{j^2+\mm}  + {\tilde r}_{i,j}(\zeta)
$$
where also $   {\tilde r}_{i,j} :=  r_{i,j} \circ (\o^\infty)^{-1}  $ satisfies \eqref{errors}.
In the direction $ \zeta = s k|k|^{-1} + w  $, $ w \cdot k = 0 $, the function
$ {\tilde f}_{k,i,j}(s) := f_{k,i,j} (s k |k|^{-1} + w) $ satisfies
$$
{\tilde f}_{k,i,j}(s_2) -  {\tilde f}_{k,i,j}(s_1) \stackrel{\eqref{errors}}\geq (s_2 - s_1) (|k| - C \e )
\geq  (s_2 - s_1) |k| / 2 \, .
$$
Since $|k|\geq 1$ (recall \eqref{nonvuoto}), by Fubini theorem,
$$
\Big|\Big\{ \zeta \in \om^\infty ({\cal O}) \, : \,  |f_{k,i,j}(\zeta )| \leq \frac{2 \eta}{1+|k|^\t} \Big\}\Big| \lessdot
\frac{\eta \rho^{n-1}}{ 1 + |k|^{\t+1}} \, .
$$
By \eqref{inveromi}
the bound \eqref{Rkij} follows.
\end{pf}

We split
\begin{equation}\label{atoya}
{\mathtt I} = {\mathtt I}_> \, \cup \, {\mathtt I}_< \quad {\rm where} \quad
{\mathtt I}_> :=
\Big\{ (k,i,j) \in {\mathtt I} \, : \, \min \{|i|, |j|\} > C_\sharp \g^{-1/3}(1 + |k|^{\t_0}) \Big\}
\end{equation}
where $ C_\sharp > C_\star $ in \eqref{freco} and $ \t_0 := n +1 $.
We set $ {\mathtt I}_< := {\mathtt I} \setminus  {\mathtt I}_> $.

\begin{lemma}\label{lem:inclu}
For all  $ (k, i, j) \in {\mathtt I}_> $ we have
\be \label{claimrq1}
{\cal R}_{k,i,j}^{\t_0} (\g^{2/3}) \subset {\cal R}_{k, i_0, j_0}^{\t_0} (2 \g^{2/3})
\ee
(see  \eqref{stimaRkl}),
$ i_0, j_0 \in \Z \setminus {\cal I} $ satisfy
\be \label{segni}
\mathtt s(i_0) = \mathtt s(i) \, , \ \mathtt s(j_0 ) = \mathtt s(j)  \, , \  | i_0 | - | j_0 | = |i|- |j|
\ee
and
\be\label{quantog}
\min\{|i_0|, |j_0| \} = \Big[ C_\sharp \g^{-1/3}(1+ |k|^{\t_0}) \Big] \, .
\ee
\end{lemma}

\begin{pf}
Since $ |j| \geq \g^{-1/3} C_\star $,
by \eqref{freco} and \eqref{omegaxi} we have the frequency asymptotic
\begin{equation}\label{fare730}
\Om^\infty_j(\xi)=|j|+\frac{\mm}{2|j|}+ \vec a \cdot \xi+a^\infty_{\mathtt s(j)}(\xi)+
O\left(\frac{\mm^2}{|j|^3}\right)+O\left(\e \frac{\gamma^{2/3}}{|j|}\right)\,.
\end{equation}
By \eqref{nonvuoto} we have
$||i|-|j||=||i_0|-|j_0||\leq C |k|$, $ |k| \geq 1 $.
If $ \xi \in {\cal O} \setminus {\cal R}_{k, i_0, j_0}^{\t_0} ( 2 \g^{2/3} ) $,
since $ |i|,|j| \geq \mu_0 := \min\{|i_0|, |j_0|\}  $
(recall \eqref{atoya} and \eqref{quantog}),  we have
\begin{eqnarray*}
|\om^\infty (\xi) \cdot k + \Om_i^\infty (\xi) - \Om_j^\infty (\xi) |
 & \geq &
| \om^\infty (\xi) \cdot k + \Om_{i_0}^\infty (\xi) - \Om_{j_0}^\infty (\xi) | \\
& & - | \Om_i^\infty (\xi) -  \Om_{i_0}^\infty (\xi) - \Om_j^\infty (\xi) + \Om_{j_0}^\infty (\xi)|   \\
 & \stackrel{\eqref{fare730}} \geq &
 \frac{ 4 \g^{2/3}}{1+|k|^{\t_0}}  - | |i| -| i_0 |  - |j| + | j_0| |  \\
 && -| a^\infty_{\mathtt s(i)}-a^\infty_{\mathtt s(i_0)} - a^\infty_{\mathtt s(j)}+a^\infty_{\mathtt s(j_0)}  |
 \\
&  & - C \e \frac{\g^{2/3}}{ \mu_0 } - C \frac{\mm^2}{\m_0^3}
 - \frac{\mm}{2} \frac{||i|-|j||}{|i| \, |j|}
- \frac{\mm}{2} \frac{||i_0|-|j_0||}{|i_0| \, |j_0|}
\\
& \stackrel{\eqref{segni}} \geq &  \frac{ 4 \g^{2/3}}{1+|k|^{\t_0}}  -
C \e \frac{\g^{2/3}}{ \mu_0 } -  C  \frac{|k|}{ \mu_0^2} \stackrel{\eqref{quantog} }  \geq
\frac{2 \g^{2/3}}{1+|k|^{\t_0}}
\end{eqnarray*}
taking $ C_\sharp $ in \eqref{quantog}
 large enough. Therefore $ \xi \in  {\cal O} \setminus {\cal R}_{k,i,j}^{ \t_0} ( \g^{2/3}) $ proving
\eqref{claimrq1}.
\end{pf}


As a corollary we deduce:

\begin{lemma}\label{cor1} 
$
\Big| \bigcup_{ (k,i,j) \in {\mathtt I}_>} {\cal R}_{k,i,j}^\t ( \g )\Big| \lessdot
\g^{2/3} \rho^{n-1} \, .$
\end{lemma}

\begin{pf}
Since $ 0<\g\leq 1  $ and
$ \t \geq \t_0  $ (see \eqref{antiochia}), we have (see  \eqref{stimaRkl})
$ {\cal R}_{k,i,j}^\t ( \g)  \subset  {\cal R}_{k,i,j}^{ \t_0} (\g^{2/3}) $.
Then Lemma \ref{lem:inclu} and  \eqref{Rkij} imply that, 
for each  $k\in \Z^n$, $ p \in \Z $ fixed
$$
\Big| \bigcup_{ (k,i,j) \in {\mathtt I}_> , \,  |i|-|j| = p }{\cal R}_{k,i,j}^\t ( \gamma) \Big| \lessdot
\frac{\g^{2/3}  \rho^{n-1}}{ 1 + |k|^{\t_0 + 1 }}
\, .
$$
Therefore
$$
\Big| \bigcup_{ (k,i,j) \in {\mathtt I}_>} {\cal R}_{k,i,j}^\t ( \g)\Big|
\lessdot \sum_{k, |p| \leq C|k|}
 \frac{\g^{2/3}  \rho^{n-1}}{ 1 + |k|^{\t_0 +1}} \lessdot
\sum_{k} \frac{  \g^{2/3} \rho^{n-1}}{ 1 + |k|^{\t_0}}
$$
proving the lemma.
\end{pf}

\begin{lemma}\label{lem:SIC}
$ \Big| \bigcup_{ (k,i,j) \in {\mathtt I}_< } {\cal R}_{k,i,j}^\t (  \g)\Big| \lessdot \g^{2/3} \rho^{n-1} $.
\end{lemma}

\begin{pf}
For all $ (k,i,j) \in {\mathtt I}_< $ such that $ {\cal R}_{k,i,j}^\t ( \g) \neq \emptyset $ we have
(see  \eqref{unionr})
$$
 \min \{|i|, |j|\} < C_\sharp \g^{-1/3}(1+|k|^{\t_0}) \, , \ j-i=  k \cdot\mathtt j
 \  \Longrightarrow \
 \max\{|i|, |j|\}  < C' \g^{-1/3}(1+|k|^{\t_0}) \, .
$$
Therefore, using also Lemma \ref{stimaRij} and \eqref{nonvuoto}
$$
\Big| \bigcup_{ (k,i,j) \in {\mathtt I}_< } {\cal R}_{k,i,j}^\t ( \g)\Big| \lessdot  \sum_k
 \sum_{|i| \leq
C'\g^{-1/3} (1+ |k|^{\t_0}) \atop  j=i+ k \cdot\mathtt j }  \frac{\g \rho^{n-1}}{1+ |k|^{\tau+1}}
 \lessdot  \sum_k    \frac{\g^{ 2/3} \rho^{n-1}}{1+ |k|^{\tau - \t_0 +1}}
$$
which, by \eqref{antiochia}, gives the lemma. 
\end{pf}

Lemmata \ref{cor1}, \ref{lem:SIC} imply \eqref{eccola}.
 This concludes the case $(k,l)\in \L_2^-.$
Let consider the other cases. The analogue of Lemma \ref{stimaRij} is
\begin{lemma}\label{capriccetto}
For $ (k,l) \in \L_0\cup \L_1\cup \L_2^+$,
$ \eta \in (0, \sqrt{\mm} / 2 ) $, we have
\be\label{treporcellini}
|{\cal R}_{kl} ( \eta ) | \lessdot \frac{\eta \rho^{n-1}}{ 1 + |k|^{\t}} \, .
\ee
\end{lemma}

\begin{pf}
We consider only the case $(k,l)\in\L_2^+$, $ l = e_i+e_j$.
By \eqref{muflone}  and \eqref{omegaxi}
$$
f_{k,i,j}(\xi):= \om^\infty (\xi) \cdot k + \Om^\infty_i (\xi) + \Om^\infty_j (\xi)
=  \om^\infty (\xi) \cdot k +  \sqrt{i^2+\mm} + \sqrt{j^2+\mm} + 2 \vec a \cdot \xi + r_{i,j}( \xi )
$$
where $ |r_{i,j}(\xi)| = O( \e  \g) \, , \  |r_{i,j}|^{\rm lip} = O( \e ) $.
Changing variables $ \zeta := \om_\infty (\xi) $ we find
\be\label{duccio}
f_{k,i,j}(\zeta):= \zeta \cdot k +  \sqrt{i^2+\mm} + \sqrt{j^2+\mm} + 2 \vec a \cdot A^{-1} (\zeta - \bar \om ) +
{\tilde r}_{i,j}( \zeta)
\ee
where also
\be\label{errorsbis}
{\tilde r}_{i,j}( \zeta)  =  O( \e  \g) \, , \  | {\tilde r}_{i,j}|^{\rm lip} = O( \e ) \, .
\ee
If $ k = \vec  a = 0 $ then the function in \eqref{duccio} is bigger than
$ \sqrt{\mm} $ and ${\cal R}_{0l}(\eta)=\emptyset $, for $ 0 \leq \eta \leq \sqrt{\mm} /2 $.
Otherwise, by \eqref{pota},  the vector
\be\label{giotto}
\tilde a:= A^T k +2 \vec a = A^T \big(k+2 (A^{-1})^T \vec a \big)\quad
{\rm satisfies}\quad
|\tilde a|\geq c=c(A, \vec  a) > 0 \, , \ \forall k \neq 0 \, .
\ee
The function
$ {\tilde f}_{k,i,j}(s) := f_{k,i,j} (s \tilde a |\tilde a|^{-1} + w) $, $ \tilde a \cdot w = 0 $, satisfies
$ {\tilde f}_{k,i,j}(s_2) -  {\tilde f}_{k,i,j}(s_1) \geq $ $ (s_2 - s_1) (|\tilde a| - C \e )
\geq  $ $ (s_2 - s_1) |\tilde a| / 2 $ by \eqref{errorsbis}.
Then \eqref{treporcellini} follows by \eqref{giotto} and Fubini theorem.
\end{pf}

By Lemma \ref{capriccetto}, \eqref{stimaRkl}, \eqref{defRkl}, \eqref{defO0} and  standard arguments (as above)
\be\label{sopra}
\Big| \bigcup_{(k,l) \in \L_0 \cup \L_1 \cup \L_2^+} {\cal R}_{kl}(\g) \Big| \lessdot \g \rho^{n-1}, \
\Big| \bigcup_{(k,p) \in \Z^{n+1} \setminus \{0\} } \tilde {\cal R}_{kp}(\g^{2/3}) \Big| \, ,
\  |  \mathcal O \setminus \mathcal O_0
| \lessdot \g^{2/3} \rho^{n-1} \, .
\ee
Finally \eqref{santadorotea}, \eqref{eccola}, \eqref{sopra}
 imply \eqref{consolatrixafflictorum}.

\section{Application to DNLW}\label{sec:DNLW} 

For
$ \vec \jmath = (j_1, \ldots, j_d)\in\mathbb{Z}^d $,
$ {\vec \s} = (\s_1,\ldots,\s_d) \in \{ \pm \}^d $ we denote
$ {\vec \s} \cdot \vec \jmath := 
 \s_1 j_1  + \ldots  + \s_d j_d  $, and,
given $ (u_j,\bar u_j)_{j \in \Z } = $ $ (u_j^+,u_j^-)_{j \in \Z } $,
we define
the monomial
$ u_{\vec \jmath}^{\vec \s} := u_{j_1}^{\s_1}\cdots u_{j_d}^{\s_d} $ (of degree $ d $).

\subsection{The partial Birkhoff normal form}\setcounter{equation}{0}
\label{sec:BNF}
We now consider the Hamiltonian \eqref{tuc} when $ F (s) = s^4 /4  $
since terms of order five or more will not make any difference,
see remark \ref{carmelo}.

After a
rescaling of the variables (and of the Hamiltonian) it
becomes 
\begin{eqnarray}\label{Ham}
H
&=&
\sum_{j\in \Z}\lambda_j u_j^+ u_j^-
+
\sum_{{\vec \jmath} \in \Z^4, {\vec \s} \in \{\pm\}^4,
 {\vec \s}\cdot {\vec \jmath}=0} u_{\vec \jmath}^{\vec \s}
=:N+G
\\
&=&
\sum_{j\in\Z} \lambda_j u_j \bar u_j +
\sum_{|\a|+|\b|=4,\,\pi(\a,\b)=0} G_{\a,\b} u^\a \bar u^\b\,,
\quad
G_{\a,\b}:=\frac{(|\a|+|\b|)!}{\a! \b!}=\frac{4!}{\a!\b!}
\,,\ \ \ \ \ \ \
\nonumber
\end{eqnarray}
where $ (u^+ , u^- ) = (u, \bar u )  \in \ell^{a,p} \times  \ell^{a,p} $
for some  $a>0$, $p>1/2$, and the momentum  is (see \eqref{momento2})
$$
\pi(\a,\b) = \sum_{j\in \Z} j(\alpha_j-\beta_j) \, .
$$
 Note that $0 \leq G_{\a,\b}\leq 4! $ (recall  $ \a! = \Pi_{i \in \Z} \a_i ! $)

\begin{lemma}\label{finzioni}
For all $ R > 0 $, $ N_0$ satisfying \eqref{caracalla},
 the  Hamiltonian $ G $ defined in \eqref{Ham} belongs to
 ${\mathcal Q}_{R}^T(N_0,\CCquarti,\quattrocc)$ and
\begin{equation}\label{ragno}
\|G\|^T_{R, N_0,\CCquarti,\quattrocc}=
 \| X_G\|_{R} \lessdot \, R^2\, .
\end{equation}
\end{lemma}
\begin{pf}
The Hamiltonian vector field $ X_G := ( -\ii \partial_{\bar u} G,
\ii \partial_{u} G) $ has components
$$
\ii \s \partial_{u^\s_l} G= \ii \s
\sum_{|\a|+|\b|=3, \pi(\a,\b)=-\s l}
G_{\a,\b}^{l,\s} u^\a \bar u^\b\, , \quad \s=\pm \, ,  \ l \in\Z \, ,
$$
where
$$
G_{\a,\b}^{l,+}=(\a_l+1)G_{\a+e_l,\b}\,,\qquad
G_{\a,\b}^{l,-}=(\b_l+1)G_{\a,\b+e_l}\,.
$$
Note that $0\leq G_{\a,\b}^{l,\s}\leq 5! $
By Definitions
\ref{MNV}, \ref{Hregular} and \eqref{normaEsr}
$$
\|X_G\|_R=\frac{1}{R}
\sup_{\|u\|_{a,p},\|\bar u \|_{a,p}< R }
\bigg(
\sum_{l\in\Z,\s=\pm} e^{2a |l|} \langle l \rangle^{2p}
\Big(
\sum_{|\a|+|\b|=3\,, \pi(\a,\b)=-\s l}
G_{\a,\b}^{l,\s} |u^\a| |\bar u^\b|
\Big)^2
\bigg)^{1/2}\,.
$$
For each component
\begin{eqnarray*}
\sum_{|\a|+|\b|=3\,,\pi(\a,\b)=-\s l}
G_{\a,\b}^{l,\s} |u^\a| |\bar u^\b|
& \lessdot & \sum_{\s_1 j_1 + \s_2 j_2 + \s_3 j_3=-\s l}
 |u_{j_1}^{\s_1}| |u_{j_2}^{\s_2}| |u_{j_3}^{\s_3}|
 \\
 &\lessdot & \big(\tilde u* \tilde u* \tilde u\big)_{-\s l}
\end{eqnarray*}
where  $ \tilde u := (\tilde u_l)_{l\in\Z} $,
$ \tilde u_j:=|u_j|+|\bar u_j|$,
and  $*$ denotes the convolution of sequences.
 Note that $ \|\tilde u\|_{a,p}\leq \|u\|_{a,p}+\|\bar u\|_{a,p} $.
Since $\ell^{a,p}$ is an Hilbert algebra, 
 $\|\tilde u* \tilde u* \tilde u\|_{a,p}\lessdot
  \|\tilde u\|_{a,p}^3 $, and
\begin{eqnarray}
 \|X_G\|_R & \lessdot &
R^{-1} \sup_{\|u\|_{a,p},\|\bar u \|_{a,p} < R }
\bigg(
\sum_{l\in\Z,\s=\pm} e^{2a |l|} \langle l \rangle^{2p}
\big|(\tilde u* \tilde u* \tilde u)_{-\s l}\big|^2
\bigg)^{1/2} \label{boundG}
\\
&\lessdot & R^{-1} \,
\sup_{\|u\|_{a,p},\|\bar u \|_{a,p} < R}
\|\tilde u* \tilde u* \tilde u\|_{a,p}
\lessdot R^{-1} \,
\sup_{\|u\|_{a,p},\|\bar u \|_{a,p} < R} \|\tilde u\|_{a,p}^3
 \lessdot R^2
\,. \nonumber
\end{eqnarray}
Moreover $ G \in {\mathcal H}^{\rm null}_R $, namely
$ G $ Poisson commutes with the momentum
$ {\cal M} := \sum_{j\in \Z} j u_j \bar u_j $, 
because (see \eqref{Poissonbraket})
\begin{equation}\label{momentumMu}
\{ {\cal M} , u_{\vec \jmath}^{\vec \s}\} = - \ii {\vec \s}\cdot {\vec \jmath} \, u_{\vec \jmath}^{\vec \s} \, .
\end{equation}
We now prove that, for all $ N \geq N_0 $, the projection
$\Pi_{N,\CCquarti,\quattrocc} G\in
{\mathcal T}_{R}(N,\CCquarti,\quattrocc)$. Hence
\eqref{ragno} follows by \eqref{boundG} (see Definition
\ref{topbis_aa}). By Definition  \ref{BL} (with $ g
\rightsquigarrow G $, no $(x,y)$-variables and $z=u$, $\bar z=\bar
u$), in particular \eqref{evaristo}, \eqref{carriego},
we get 
\begin{eqnarray}
\Pi_{N,\CCquarti,\quattrocc} G
&=&
\sum_{|m|,|n|>  \CCNquarti, \s,\s'=\pm }
G_{m,n}^{\s,\s'}(w_L) u_m^\s \bar u_n^{\s'}\qquad {\rm with}
\nonumber
\\
G_{m,n}^{\s,\s'}(w_L)
&=&
\sum_{\sum_{j\in\Z} |j|(\a_j+\b_j)<\quattrocc N^L,
\atop \pi(\a,\b)=-\s m-\s' n}
G_{\a,\b,m,n}^{\s,\s'} u^\a \bar u^\b
\qquad {\rm and}
\nonumber
\\
G^{+,+}_{\a,\b,m,n}
&=&
\frac{1}{2-\d_{mn}}G_{\a+e_m+e_n,\b}
= \frac{1}{2-\d_{mn}} \frac{4!}{(1+\d_{mn})!}=12=
G^{-,-}_{\a,\b,m,n}
\nonumber
\\
G^{+,-}_{\a,\b,m,n}
&=&
G_{\a+e_m,\b+e_n}=24
=G^{-,+}_{\a,\b,m,n}\, .
\nonumber
\end{eqnarray}
These coefficients trivially satisfy \eqref{marco}
(with $f \rightsquigarrow G$), so
$\Pi_{N,\CCquarti,\quattrocc} G\in
{\mathcal T}_R (N,\CCquarti,\quattrocc)$.
\end{pf}

We now perform a Birkhoff semi--normal form on the \textsl{tangential sites}
\begin{equation}\label{S}
{\cal I} := \{\pluto_1,\dots, \pluto_n \} \subset \Z \, , \quad
\pluto_1 < \dots < \pluto_n \, ,
\end{equation}
recall \eqref{cC}. Let  $ {\cal I}^c:= \Z\setminus {\cal I}. $

Set
\begin{equation}\label{giulio}
\overline G:=\frac12\sum_{i\, {\rm or}\, j\in {\cal I}}
\overline G_{ij} u_i^+ u_i^- u_j^+u_j^- \, ,
\quad
\overline G_{ij} :=12 (2-\d_{ij}) \, ,
\quad
\hat G:=\sum_{\vec \jmath\in\Z^4,\,\vec \s\in\{+,-\}^4,
\atop
 {\vec \s}\cdot {\vec \jmath}=0,\, {\vec \jmath}\in ({\cal I}^c)^4}
u^{\vec \s}_{\vec \jmath}\, .
\end{equation}
By 
\eqref{ragno} and noting that $ \overline G $, $\hat G$
are projections of $ G $, for $ R > 0 $,
$N_0$ satisfying \eqref{caracalla}, we have
\begin{equation}\label{stellamaris}
\|\overline G\|^T_{R, N_0,\CCquarti,\quattrocc}\,,\
\|\hat G\|^T_{R, N_0,\CCquarti,\quattrocc}
\lessdot R^2 \, .
\end{equation}

\begin{proposition}\label{BNF}
{\bf (Birkhoff normal form)}
For any 
$ {\cal I} \subset \Z $ and $ \mm > 0 $, there exists $ R_0 >  0 $ and
a real analytic, symplectic change of variables
$$
\Gamma 
: B_{R/2}\times B_{R/2}
\subset \ell^{a,p}\times \ell^{a,p} \ \to\
B_{R}\times B_{R} \subset \ell^{a,p}\times \ell^{a,p} \, , \quad  0 < R < R_0  \, ,
$$
that takes the Hamiltonian $ H = N + G $ in \eqref{Ham} into
\be\label{Birkhoff}
H_{\rm Birkhoff} := H \circ \Gamma =  N+\overline{G}+\hat G + K
\ee
where
$\overline{G},\hat G$ are defined in \eqref{giulio} and
\begin{equation}\label{Kappa}
K :=\sum_{\vec \jmath\in\Z^{2d},\,\vec \s\in\{+,-\}^{2d},\atop
 \,d\geq 3,\, {\vec \s}\cdot {\vec \jmath}=0}
K_{{\vec \jmath},{\vec \s}}u^{\vec \s}_{\vec \jmath}
\end{equation}
satisfies, for $ N_0' := N_0' ( \mm,{\cal I},L,b) $ large enough,
\begin{equation}\label{stellamarisbis}
\|K\|^T_{R/2, N_0',\CCterzi, \trecc}
\lessdot  R^4\, .
\end{equation}
\end{proposition}

The rest of this subsection is devoted to the proof of Proposition \ref{BNF}.
We start following the strategy of \cite{Po3}.  By \eqref{Poissonbraket} the Poisson bracket
\be\label{PoissN}
\{ N , u_{\vec \jmath}^{\vec \s}\} = - \ii {\vec \s}\cdot \l_{\vec \jmath} \, u_{\vec \jmath}^{\vec \s}
\ee
where  $ \l_{ \vec{\jmath}} := (  \l_{j_1}, \ldots,  \l_{j_d} ) $ and $ \l_j := \l_j(\mm) := \sqrt{j^2+\mm} $.

\smallskip

The following lemma extends  Lemma 4 of \cite{Po3}.

\begin{lemma}\label{poeschel}
{\bf (Small divisors)}
Let $ \vec \jmath\in \Z^4,\, \vec \s \in \{\pm\}^4 $ be such that  $ {\vec \s}\cdot {\vec \jmath}=0 $ and (up to permutation of the indexes)
\begin{eqnarray} 
&  &  \; \vec \jmath = 0 \, , \sum_{i=1}^4 \s_i\neq 0 
\,,   \label{pampers3}\\ & &\nonumber \\
& {\rm or} &
\vec \jmath = (0,0,q,q) \, , \ q \neq 0 \, , 
\s_1 = \s_2   \, ,  \label{pampers2}  \\& &\nonumber\\
&  {\rm or} &
{\vec \jmath} = (p, p, -p, -p)\, ,  \ p \neq 0  \, ,
\s_1 = \s_2 \,,   \label{pampers4}   \\& &\nonumber \\
&{\rm or} & \
{\vec \jmath}\neq(p, p,q, q)\,  . \label{pampers1}
\end{eqnarray}
Then,  there exists an absolute constant $c_* >  0 $, such that,
for every $\mm\in(0,\infty) $,   
\begin{equation}\label{pecorino}
|{\vec \s}\cdot \l_{\vec \jmath}(\mm)|\geq \frac{c_* \mm}{(n_0^2+\mm)^{3/2}} > 0 \quad
\mbox{where}
\quad
n_0:= 
\min \{ \langle j_1 \rangle,  \langle j_2  \rangle,   \langle j_3  \rangle,   \langle j_4  \rangle \} \, .
\end{equation}
\end{lemma}

\begin{pf}
In the Appendix.
\end{pf}

The map  $ \Gamma := \Phi^1_F $ is obtained as the time-$1$ flow generated
by the Hamiltonian
\begin{equation}\label{brodino}
F:=- \sum_{  \vec\jmath\cdot \vec\s=0 \, , \vec\s \cdot \l_{\vec\jmath}\neq 0 \\ \atop
\ {\rm and}\ {\vec \jmath}\notin ({\cal I}^c)^4}
\frac{\ii}{{\vec \s}\cdot \l_{\vec \jmath}} u^{\vec \s}_{\vec \jmath}
\end{equation}
We notice that the condition $\vec\jmath\cdot \vec\s=0 \, , \vec\s \cdot \l_{\vec\jmath}\neq 0$ is equivalent to  requiring
that
$  \vec\jmath\cdot \vec\s=0 $ and $ {\vec \jmath},{\vec \s}$ satisfy  \eqref{pampers3}-\eqref{pampers1}.
 By Lemma \ref{poeschel}
there is a constant $ \bar c >  0  $ (depending only on $\mm$ and ${\cal I}$) such that
\begin{equation}\label{chiaraeluna}
\vec\jmath\cdot \vec\s=0 \, , \vec\s \cdot \l_{\vec\jmath}\neq 0 \ {\rm and}\
{\vec \jmath}\notin ({\cal I}^c)^4
\quad
\Longrightarrow
\quad
|{\vec \s}\cdot \l_{\vec \jmath}|\geq \bar c > 0 \, .
\end{equation}
We have proved that  the moduli of the small divisors in \eqref{brodino} are uniformly bounded away from zero. 
Hence $ F $ is well defined and, arguing as in Lemma \ref{finzioni}, we get
\begin{equation}\label{ameriga}
\|X_F\|_{R}\lessdot  R^2\, .
\end{equation}
Moreover $ F\in {\mathcal H}_{R}^{\rm null}  $ because
in   \eqref{brodino} the sum is restricted to $\vec\s\cdot\vec \jmath=0 $ (see also \eqref{momentumMu}).

\begin{lemma}\label{bototo}
$ F $ in \eqref{brodino} solves the homological equation
\begin{equation}\label{claudia}
\{N,F\} + G = {\rm ad}_F(N) + G = \overline G + \hat G
\end{equation}
where $ \overline G $, $ \hat G $ are defined in \eqref{giulio}.
\end{lemma}

\begin{pf}
We claim that the only $ \vec \jmath \in \Z^4 $, $ \vec \s \in \{\pm\}^4 $ with $  \vec \jmath \cdot  \vec \s = 0 $ which do not satisfy \eqref{pampers3}-\eqref{pampers1} have the form
\be\label{zerisecchi}
j_1=j_2 \, ,  \ j_3=j_4 \, , \  \s_1 = -\s_2 \, , \
\s_3=-\s_4 \ {\rm (or  \ permutations \ of\ the\ indexes)} \, .
\ee
Indeed:

If $ \vec \jmath = 0 $, $ \sum_i \s_i=0 $: the $\s_i$ are pairwise equal and  \eqref{zerisecchi}  holds.

If  $ \vec \jmath = (0,0,q,q) $, $ q \neq 0 $,  and $ \s_1 = - \s_2 $: by  $  \vec \jmath \cdot  \vec \s = 0 $  we have also $ \s_3 = - \s_4 $ and \eqref{zerisecchi}  holds.

If $ {\vec \jmath} = (p, p, -p, -p)$,  $ p \neq 0 $ and  $ \s_1 =-\s_2   $:    by  $  \vec \jmath \cdot  \vec \s = 0 $  we have also $ \s_3 = - \s_4 $ and    \eqref{zerisecchi}  holds.

If $  j_1=j_2 $, $ j_3 = j_4 $,  $ j_1 , j_3 \neq 0 $, $ j_1 \neq - j_3 $:
\\[1mm]
{\sc Case 1:} $   j_1\neq  j_3 $.  Then
$0={\vec \s}\cdot {\vec \jmath}= (\s_1+\s_2)j_1+ (\s_3+\s_4)j_3$  implies
$\s_1=-\s_2 $,
$\s_3=-\s_4$.
\\[1mm]
{\sc Case 2:}
$ j_1 = j_3 $ and so $j_1=j_2= j_3=j_4\neq 0$. Hence
$0= (\s_1+\s_2+\s_3+\s_4)j_1$ and \eqref{zerisecchi} follows.

\smallskip

By  \eqref{brodino} and \eqref{PoissN}
all the monomials in $\{N,F\} $ cancel the monomials of $ G $ in \eqref{Ham}
except for those in $ \hat G $ (see \eqref{giulio}) and those
of the  form $ |u_p |^2 | u_q |^2 $, $ p $ or $ q  \in { \cal I } $,
which contribute to $ \overline G $.  The expression in \eqref{giulio} of 
$\overline G $ follows by counting the multiplicities.
\end{pf}

The Hamiltonian $ F  \in {\mathcal H}_{R}^{\rm null}  $ in \eqref{brodino} is  quasi-T\"oplitz:

\begin{lemma}\label{isaia}
Let $ R > 0 $. If
$N_0 := N_0 (\mm,{\cal I},L,b) $ is large enough, then
$ F $ defined in \eqref{brodino} belongs to
 ${\mathcal Q}_{R}^T(N_0,\CCquarti, \quattrocc)$ and
\begin{equation}\label{libellula}
\| F \|_{R, N_0, \CCquarti, \quattrocc}^T \lessdot R^2 \, .
\end{equation}
\end{lemma}

\begin{pf}
We have to show that $ F  \in {\mathcal H}_{R}^{\rm null}  $ verifies Definition \ref{topbis_aa}.
For all  $N\geq N_0 $, we compute, by \eqref{brodino} and
 Definition \ref{BL} (in particular \eqref{evaristo}),  the projection
\begin{equation}\label{geremia}
\Pi_{N,\CCquarti, \quattrocc} F
=
\sum_{|n|,|m|>CN/4\,,\atop
\sigma,\sigma'=\pm \,,|\sigma m+\sigma' n|<\quattrocc N^L}
 F^{\s,\s'}_{m,n}(w^L) u_m^\s  u_n^{\s'}
\end{equation}
where
\begin{eqnarray}
 F^{\s,\s'}_{m,n}(w^L)
&:=&
-12 \ii \sum_{|i|+|j|<\quattrocc N^L,\ i\,{\rm or}\, j\in {\cal I}\,,\atop
{\s_i i+ \s_j j +\s m +\s' n=0,\ i\neq j \, {\rm if}\, m=n}}
\frac{u_i^{\s_i} u_j^{\s_j}}{\s_i\l_i+\s_j\l_j+\s\l_m+\s'\l_n} \qquad
\label{camomilla}
\\
&=&
\sum_{\sum_j |j|( \a_j +  \b_j)<\quattrocc N^L,\ \sum_{j\in {\cal I}}(\a_j+ \b_j)>0\,,
\atop  \s m+\s' n=-\pi(\a,\b),\  |\a|+|\b|=2, \ \a\neq \b \, {\rm if}\, m=n}
 F^{\s,\s'}_{\a,\b,m,n} u^{\a} \bar u^{\b}
 \label{ezechiele}
\end{eqnarray}
and
\be\label{lambdab}
 F^{\s,\s'}_{\a,\b,m,n}
:=
-\frac{24\ii}{\a ! \b!} \frac{1}{\l_{\a,\b}+\s \l_m+\s'\l_n} \, ,
\quad \  \l_{\a,\b}
:=
\sum_h \l_h  (\a_h -\b_h) \,.
\ee

Notice that in   \eqref{camomilla} the restriction   $i\neq j $ if $ m=n$
is equivalent to requiring $$\{(i,j,m,n),(\s_i,\s_j,\s,\s')\}\neq\{(i,i,m,m),(\s_i,-\s_i,\s,-\s)\}\,,$$ see Formula \eqref{brodino} and \eqref{zerisecchi}.
Indeed
 if $m=n$ , $|i|+|j|<\quattrocc N^L$ and $|m|> CN /4 $ then, by momentum conservation, we have a contribution to \eqref{camomilla} only if $\s=-\s'$ and hence
$|i|=|j|$.

We define the  T\"oplitz approximation  
\begin{equation}\label{Ft}
\Ft := \sum \tilde F^{\s,\s'}_{m,n}(w^L) u_m^\s  u_n^{\s'}
\quad
{\rm with}
\quad
\tilde F^{\s,\s'}_{m,n}(w^L):=\sum \tilde  F^{\s,\s'}_{\a,\b,m,n} u^{\a} \bar u^{\b}
\end{equation}
where the indexes in the two sums have the same restrictions
as in \eqref{geremia}, \eqref{ezechiele}, respectively, and the coefficients are
\begin{equation}\label{sonno}
\Ft^{\s,-\s}_{\a,\b,m,n} := -\frac{24\ii}{\a ! \b!}
\frac{1}{\l_{\a,\b}+\s |m|-\s |n|}\,, \qquad \
\Ft^{\s,\s}_{\a,\b,m,n} :=0\,.
\end{equation}
The coefficients in \eqref{sonno} are well defined for $ N \geq N_0 $ large enough, because
\begin{eqnarray}
| \l_{\a,\b}+\s |m|-\s |n| | & \geq &  | \l_{\a,\b}+\s \l_m -\s \l_n | - | \l_m - |m| | - | \l_n - |n| | \nonumber \\
& \stackrel{\eqref{chiaraeluna}, \eqref{primoas}} \geq & \bar c -  \frac{\mm}{2}\left( \frac{1}{|m|}+\frac{1}{|n|}\right) \geq
 \bar c -  \quattrofrattoCC\frac{\mm}{N_0} \geq
 \frac{\bar c}{2} \, , \label{mimmo}
\end{eqnarray}
($\bar c$ defined in \eqref{chiaraeluna})
having used the elementary inequality
\be\label{primoas}
|\sqrt{ n^2 + \mm } - |n | | \leq 1 / (2|n|)  \, .
\ee
Then \eqref{Ft}, \eqref{sonno}, \eqref{mimmo} imply, arguing  as in the proof of Lemma \ref{finzioni}, that 
\begin{equation}\label{ameriga2}
\|X_\Ft\|_{R} 
\lessdot  R^2 \,.
\end{equation}
For proving that 
$ \Ft \in {\mathcal T}_{R}(N_0,\CCquarti, \quattrocc) $
we have to show
\eqref{marco} (with $f\rightsquigarrow \Ft$), namely
\begin{equation}\label{bagnetto}
\Ft^{\s,\s'}_{\a,\b,m,n} = \Ft^{\s,\s'}_{\a,\b}(\mathtt s(m),\s m+\s' n)
\end{equation}
with
$$
\Ft^{\s,-\s}_{\a,\b}(s,h) := -\frac{24\ii}{\a ! \b!}
\frac{1}{\l_{\a,\b}+sh} \,,\qquad \Ft^{\s,\s}_{\a,\b}(s,h)=0
\, , \quad s=\pm \, , \ h\in\Z \, .
$$
Recalling \eqref{sonno}, 
this is obvious when $\s'=\s$.
When  $\s'=-\s$ 
we first note that $ \mathtt s(m)=\mathtt s(n)$.
Indeed the restriction on the first sum in \eqref{Ft} is
(recall \eqref{geremia})
$ |m|,|n|>\CCNquarti $, $ |\s m-\s n|<\quattrocc N^L $,
which implies
$\mathtt s(m)=\mathtt s(n)$ by \eqref{caracalla}.
Then
$$
\s |m|-\s |n| = \s \mathtt s(m)m-\s \mathtt s(n)n=\mathtt s(m)(\s m-\s n)
$$
and \eqref{bagnetto} follows.
We have proved that $\Ft \in {\mathcal T}_{R}(N_0,\CCquarti, \quattrocc).$

The  T\"oplitz defect, defined by \eqref{limi_aa}, is  
\begin{equation}\label{Fr}
\Fr:=\sum \hat F^{\s,\s'}_{m,n}(w^L) u_m^\s  u_n^{\s'}
\quad
{\rm with}
\quad
\hat F^{\s,\s'}_{m,n}(w^L):=\sum \hat  F^{\s,\s'}_{\a,\b,m,n} u^{\a} \bar u^{\b}
\end{equation}
where the indexes in the two sums have the same restrictions
as in \eqref{geremia}-\eqref{ezechiele}, and 
\begin{eqnarray}
\Fr^{\s,\s}_{\a,\b,m,n} &=& - \frac{24\ii}{\a ! \b!}
\frac{N}{\l_{\a,\b}+\s \l_m+\s \l_n} \label{malachia}
\\
\Fr^{\s,-\s}_{\a,\b,m,n} &=& -N\frac{24\ii}{\a ! \b!}
\left( \frac{1}{\l_{\a,\b}+\s \l_m -\s \l_n} -
\frac{1}{\l_{\a,\b}+\s |m|-\s |n|} \right) \nonumber \\
& = & \frac{24\ii}{\a ! \b!}
\frac{  N \s(\l_m -|m| - \l_n + |n|)}{(\l_{\a,\b}+\s \l_m -\s \l_n) ( \l_{\a,\b}+\s |m|-\s |n|     )  }
\label{daniele}
\end{eqnarray}
We now proof that the
coefficients
in  \eqref{malachia}-\eqref{daniele}
are bounded by a constant independent of $N$.

The coefficients 
in  \eqref{malachia} are bounded because
$$
|\l_{\a,\b}| \leq
\sum_h \l_h (|\a_h|+|\b_h|)\leq
\sum_h |h| (|\a_h|+|\b_h|) +\sqrt{\mm} \sum_h (|\a_h|+|\b_h|)
\leq \quattrocc N^L + 2\sqrt{\mm}
$$
by \eqref{lambdab}-\eqref{ezechiele} (note that $\l_h\leq |h|+\sqrt\mm$)
and
$$
| \l_{\a,\b}+\s \l_m+\s \l_n| \geq | \l_m + \l_n | - |\l_{\a,\b}|
 \geq \CCNmezzi - \quattrocc N^L - 2\sqrt{\mm} \geq  \CCNquarti
$$
for $ N \geq N_0 $ large enough.

 The  coefficients in 
 \eqref{daniele} are bounded by  \eqref{chiaraeluna}, \eqref{mimmo}, and
$$
N | \l_m -|m| - \l_n + |n| |
\stackrel{\eqref{primoas}} \leq N \frac{\mm}{2} \Big( \frac{1}{|m|} +  \frac{1}{|m|} \Big) \leq
 \quattrofrattoCC\mm  \, .
$$
 Hence arguing  as in the proof of Lemma \ref{finzioni} we get
\be\label{ameriga3}
 \|X_\Fr\|_{R} \lessdot R^2 \,.
\ee
In conclusion, \eqref{ameriga}, \eqref{ameriga2}, \eqref{ameriga3}  imply \eqref{libellula}
(recall \eqref{unobisbis_aa}).
\end{pf}

\smallskip

\noindent
{\sc Proof of Proposition \ref{BNF} completed.}
We have
\begin{eqnarray*}
e^{{\rm ad}_F}H
&=&
e^{{\rm ad}_F} N+e^{{\rm ad}_F}G
=N+\{N,F\}+\sum_{i\geq 2}\frac{1}{i!}{\rm ad}_F^i(N)+G+  \sum_{i\geq 1}\frac{1}{i!}{\rm ad}_F^i(G)
\\
&\stackrel{\eqref{claudia}}=&
N+ \overline G+\hat G+
\sum_{i\geq 1}\frac{1}{(i+1)!}{\rm ad}_F^i\big({\rm ad}_F(N)\big)+  \sum_{i\geq 1}\frac{1}{i!}{\rm ad}_F^i(G)
\\
&=& N+ \overline G+\hat G+K
\end{eqnarray*}
where, using again \eqref{claudia},
\begin{equation}\label{Kappabis}
K :=
\sum_{i\geq 1}\frac{1}{(i+1)!}{\rm ad}_F^i(\overline G+\hat G-G)+
\sum_{i\geq 1}\frac{1}{i!}{\rm ad}_F^i G=:K_1+K_2\,.
\end{equation}
\\[1mm]
{\sc Proof of \eqref{Kappa}.}
We claim that in the expansion of $ K $ in \eqref{Kappabis}  there are only
monomials $  u_{\vec \jmath}^{\vec \s} $ with
$\vec \jmath\in\Z^{2d},\,\vec \s\in\{+,-\}^{2d} $, $ d \geq 3 $.
Indeed  $F,G,\overline G$, $\hat G$
contain only monomials of degree four and, for any monomial
$\mathfrak m$,
${\rm ad}_F (\mathfrak m)$ contains only monomials
of degree equal to the ${\rm deg}(\mathfrak m) + 2 $.
The restriction $\vec\s\cdot\vec \jmath=0$
follows by the Jacobi identity \eqref{Jacobi},
since $F,G,\overline G, \hat G$ preserve momentum, i.e. Poisson commute with $ M $.
\\[1mm]
{\sc Proof of \eqref{stellamarisbis}.} 
We apply  Proposition \ref{main} with
(no $(x,y)$ variables and)
$$
f\rightsquigarrow F\,, \quad
g\rightsquigarrow
\begin{cases}
 \overline G+\hat G-G \quad {\rm for} \ K_1 \, ,  \cr
G  \quad \, \qquad \qquad {\rm for} \ K_2 \, ,
\end{cases}
r \rightsquigarrow R \,,\quad
r' \rightsquigarrow R/2\,,\quad
\d\rightsquigarrow 1/2\,,
$$
$$
\theta \rightsquigarrow \CCquarti\,, \quad
\theta' \rightsquigarrow \CCterzi\,, \quad
\mu \rightsquigarrow \quattrocc\,, \quad
\mu' \rightsquigarrow \trecc\,,
$$
 $ N_0$ defined in Lemma \ref{isaia}
and $N_0'\geq N_0$
satisfying \eqref{stoppa} and
\begin{equation}\label{Giobbebis}
\kappa  (N_0')^{b-L} \ln N_0'  \leq \cc \,, \qquad\quad
(\CC  + \kappa) (N_0')^{L-1} \ln N_0' \leq 1 / 2 \, .
\end{equation}
Note that \eqref{Giobbe} follows by \eqref{Giobbebis}.
By \eqref{libellula}, the assumption  \eqref{piccof} is verified for
every $ 0<R<R_0,$ with $R_0$
small enough.
Then Proposition \ref{main} applies and
 \eqref{stellamarisbis} follows by
 \eqref{gPhif12} (with $h\rightsquigarrow 1$),
 \eqref{ragno},
 \eqref{libellula} and \eqref{stellamaris}.
 \rule{2mm}{2mm}

\subsection{Action--angle variables}\label{sec:aa}

We introduce action-angle variables on the  tangential sites
$ {\cal I } := \{\pluto_1,\dots, \pluto_n \} $  (see \eqref{S}) via the
analytic and symplectic map
\be\label{variableAA}
\Phi (x,y,z,\bar z;\xi) := (u,\bar u)
\ee
defined by
\be\label{actionangle}
u_{\pluto_l } := \sqrt{\xi_l +y_l} \, e^{{\rm i} x_l}, \, {\bar u}_{\pluto_l } := \sqrt{\xi_l +y_l} \, e^{- {\rm i} x_l}, \,
 l =1, \dots, n \, , \ \ u_j := z_j \, ,  \  \bar u_j := \bar z_j \, ,
\, j \in \Z \setminus {\cal I} \, .
\ee
Let
\be\label{paraO}
\mathcal O_\rro := \Big\{ \xi \in \R^n \, : \,
\frac{\rro}{2}\leq \xi_l \leq \rro \, , \ l = 1, \ldots, n \Big\} \, .
\ee
\begin{lemma} {\bf (Domains)}
Let $r,R,\rho > 0 $ satisfy
\begin{equation}\label{condro}
16 r^2 < \rro\,,\quad  \
 \rro  = C_* R^2  \quad {\rm with} \quad C_*^{-1} := 48 n \kappa^{2p} e^{2 ( s+ a\kappa)} \, .
\end{equation}
Then, for all $\xi \in  \mathcal O_\rro \cup \mathcal  O_{2 \rro} $,
the map
\be\label{inclusio}
\Phi(\,\cdot\,;\xi) : D(s,2r) \to {\cal D}(R/2) := B_{R/2} \times B_{R/2}  \subset \ell^{a,p} \times \ell^{a,p}
\ee
is well defined
and analytic  $($$D(s,2r)$ is defined in \eqref{Dsr} and $ \kappa $ in
\eqref{caracalla}$)$.
\end{lemma}

\begin{pf}
Note first that for  $ (x,y,z,\bar z)\in D(s,2r) $ we have (see \eqref{seconda}) that
$ |y_l | < 4 r^2 \stackrel{\eqref{condro}} < \rho / 4 < \xi_l $,
$ \forall \xi \in {\cal O}_\rho \cup  {\cal O}_{2\rho} $. Then
the map $ y_l \mapsto \sqrt{\xi_l + y_l} $ is
well defined and analytic. Moreover, for $ \xi_l  \leq 2 \rho $,
$ |\pluto_l | \leq \kappa $, $ x \in \T^n_s $, $\| z \|_{a,p} < 2 r $, we get
\begin{eqnarray*}
\| u(x,y,z,\bar z; \xi) \|_{a,p}^2 & \stackrel{\eqref{variableAA}} = &
\sum_{l =1}^n ( \xi_l + y_l) |e^{2{\rm i} x_l}| |\pluto_l |^{2p}
e^{2a|\pluto_l |} + \!\!
  \sum_{j\in \Z\setminus \cal I}\!\! |z_j|^2 \langle j \rangle^{2p} e^{2 a|j|} \\
& \leq &  n  \Big( 2\rro+ \frac{\rro}{4} \Big) e^{2 s}  \kappa^{2p}  e^{2 a\kappa}   + 4 r^2 \stackrel{\eqref{condro}} <
R^2 / 4 \,
\end{eqnarray*}
proving \eqref{inclusio} (the bound for  $ \bar u $ is the same).
\end{pf}

Given a function $ F : \, {\cal D}(R/2)\to \C $,
the previous Lemma shows that the composite map  $ F \circ \Phi :
D(s, 2r) \to \C $. The main result of  this section is Proposition
\ref{qtop0}:
if $ F $ is  quasi-T\"oplitz in the variables $ (u, \bar u) $ then
the composite $ F \circ \Phi $ is  quasi-T\"oplitz in the
variables $ (x,y,z, \bar z ) $  (see Definition \ref{topbis_aa}).

We write
\be\label{defFab} F = \sum_{\a,\b} F_{\a,\b} \mathfrak m_{\a,\b} \, ,
\quad
\mathfrak m_{\a,\b}  :=
(u^{(1)})^{\alpha^{(1)}}(\bar u^{(1)})^{\beta^{(1)}}(u^{(2)})^{\alpha^{(2)}}(\bar u^{(2)})^{\beta^{(2)}} \, ,
\ee
where
$$
u= (u^{(1)}, u^{(2)}) \, , \quad  u^{(1)}:= \{u_j\}_{j\in \cal I}
\, , \ u^{(2)}:= \{u_j\}_{j\in \Z\setminus \cal I}  \, , \quad {\rm similarly \ for}  \ \bar u \, ,
$$
and
\be\label{defa1b1}
 (\a,\b) =  (\a^{(1)}  + \a^{(2)} ,\b^{(1)} + \b^{(2)}) \, , \ \
(\a^{(1)},\b^{(1)}):= \{\a_j,\b_j\}_{j\in \cal I} \, , \ \
(\a^{(2)},\b^{(2)}):= \{\a_j,\b_j\}_{j\in \Z\setminus \cal I} \, .
\ee
We define
\begin{equation}\label{ananas}
{\cal H}^d_R := \Big\{  F\in \mathcal H_R\  :  \ F
= \sum_{|\a^{(2)}+\b^{(2)}| \geq d} F_{\a,\b}  u^\a\bar u^\b \Big\}\, .
\end{equation}

\begin{proposition}{\bf (Quasi--T\"oplitz)} \label{qtop0}
Let $ N_0,\teta,\mu,\mu' $ satisfying \eqref{caracalla} and
\begin{equation}\label{fuso}
(\mu'-\mu)N_0^L> N_0^b\,,\qquad  N_0
2^{-\frac{N_0^b}{2\kappa}+1}<1 \, .
\end{equation}
If
$ F \in \mathcal Q_{R/2}^T (N_0,\theta,\mu') \cap {\cal H}_{R/2}^{d} $ with $ d = 0, 1 $,
then $ f := F \circ \Phi\in \mathcal Q_{s,r}^T(N_0,\theta,\mu)$
and
\begin{equation}\label{topstim}
 \|f\|_{s,r, N_0,\theta,\mu,\mathcal O_\rro}^T
 \lessdot (  8r/R )^{d-2}\|F\|_{R/2, N_0,\theta,\mu'}^T \, .
 \end{equation}
\end{proposition}

The rest of this section is devoted to the proof of Proposition  \ref{qtop0}.
Introducing the action-angle variables \eqref{actionangle}
 in \eqref{defFab}, and using the Taylor expansion
\be\label{fractional}
(1 + t )^\g = \sum_{h \geq 0} \binom{\g}{h} t^h \, , \quad   \binom{\g}{0} := 1 \, ,
\ \ \binom{\g}{h} := \frac{\g (\g-1) \ldots (\g - h + 1)}{h!} \, , \ h \geq 1 \, ,
\ee
we get
\be\label{fnewc}
 f := F\circ \Phi  = \sum_{k,i,\a^{(2)},\b^{(2)}}f_{k,i,\a^{(2)},\b^{(2)}} e^{ \ii k
\cdot x} y^i z^{\a^{(2)}} \bar z^{\b^{(2)}} \ee
with   Taylor--Fourier coefficients
\be\label{coeffit} f_{k,i,\a^{(2)},\b^{(2)}} :=
\sum_{\a^{(1)}-\b^{(1)}=k} F_{\a,\b}
\prod_{l=1}^n\xi_l^{\frac{\a^{(1)}_l+\b^{(1)}_l}{2} -i_l}
\binom{\frac{\a^{(1)}_l +\b^{(1)}_l}{2}}{i_l} \, .
\ee
We need an upper bound on the binomial coefficients.
\begin{lemma}\label{tecno}
 For $ |t| < 1 / 2 $ we have
\begin{equation}\label{embeh}
(i) \, \sum_{h\geq 0}| t|^h \Big|\binom{\frac k2}{ h} \Big| \leq 2^{k}\,, \  \forall k \geq 0 \, ,
\quad
(ii) \,
\sum_{h\geq 1}| t|^h \Big|\binom{\frac k2}{ h} \Big| \leq 3^{k} |t| \, , \   \forall k \geq 1  \,.
\end{equation}
\end{lemma}
\begin{pf}
By \eqref{fractional} and the definition of majorant (see \eqref{formalpower1A}) we have
\be\label{caten}
\sum_{h\geq 0} \Big|\binom{\frac k2}{ h} \Big| t^h  =
M (1+t)^{\frac k2}
\stackrel{\eqref{sumproA}}
\prec (M (1+t)^{\frac12})^k = \Big( \sum_{h \geq 0} \Big|\binom{\frac 12}{ h}\Big|  t^h \Big)^k
\prec \Big( \sum_{h \geq 0} t^h \Big)^k
\ee
because 
$ \Big| \binom{\frac 12}{h} \Big| \leq 1$ by \eqref{fractional}.
For $ |t| < 1/2 $ the bound \eqref{caten} implies  \eqref{embeh}-($i$). Ne
$$
\sum_{h\geq 1}| t|^h \Big|\binom{\frac k2}{ h} \Big|
\leq
|t| \sum_{h\geq 0}| t|^h \Big|\binom{\frac k2}{ h+1} \Big|
\stackrel{\eqref{fractional}} =
|t|\sum_{h\geq 0}| t|^h \Big|\binom{\frac k2}{ h} \Big| \frac{|\frac k2 - h|}{h+1}
\leq
k|t|\sum_{h\geq 0}| t|^h \Big|\binom{\frac k2}{ h} \Big| \stackrel{\eqref{embeh}\!-(i)} \leq \!\!\!\!\!\!k 2^k |t|
$$
which implies \eqref{embeh}-($ii$) for $ k \geq 1 $. 
\end{pf}

\begin{lemma}\label{artaserse}
{\bf ($ M$-regularity)} If $ F \in {\cal H}_{R/2}^d $ then $ f := F
\circ \Phi\in \mathcal H_{s,2r}$  and
\begin{equation}\label{coord}
\|X_f\|_{s,2r,\mathcal O_\rro\cup \mathcal O_{2\rro}}
\lessdot  (8r/R)^{d-2} \|X_F\|_{R/2} \, .
\end{equation}
Moreover if $F$ preserves momentum then so does $F\circ \Phi$.
\end{lemma}

\begin{pf}
We first bound the majorant  norm  
\begin{equation}\label{mattonella}
\| f\|_{s,2r,\mathcal O_\rro\cup \mathcal O_{2\rro}}
\stackrel{ \eqref{fnewc},  \eqref{ananas}} {:=}
\!\!\! \!\!\! \!\!  \sup_{\xi\in \mathcal O_\rro\cup \mathcal
O_{2\rro}} \sup_{(y,z, \bar z) \in D(2r)} \sum_{k,i,|\a^{(2)}+\b^{(2)}|\geq d}\!\!\! \!\!
|f_{k,i,\a^{(2)}, \b^{(2)}}| e^{|k|s} |y^i|
|z^{\a^{(2)}}| |{\bar z}^{\b^{(2)}}| \, .
\end{equation}
Fix $\a^{(2)},\b^{(2)} $.
Since for all  $ \xi \in  \mathcal O_{\rro} \cup \mathcal O_{2\rro} $, $ y \in B_{(2r)^2}  $, we have
$ |y_l / \xi_l | < 1 / 2 $ by \eqref{condro}, we have
\begin{eqnarray}
&&
\sum_{k} e^{|k|s}\sum_i| f_{k,i,\a^{(2)},\b^{(2)}}||y|^i
\label{achemenide}
\\
&\stackrel{\eqref{coeffit}} \leq&
\sum_{\a^{(1)},\b^{(1)}} e^{ s( |\a^{(1)}|+|\b^{(1)}|)}|F_{\a,\b} |
 \xi^{\frac{\a^{(1)}+\b^{(1)}}{2} } \prod_{l=1}^n\sum_{i_l \geq 0} \left |
 \frac{y_l}{\xi_l}\right |^{i_l}
 \Big|\binom{\frac{\a^{(1)}_l +\b^{(1)}_l}{2}}{i_l} \Big|
 \label{achemenidebis}
\\
 &\stackrel{\eqref{embeh}} \leq&
\sum_{\a^{(1)},\b^{(1)}} e^{ s( |\a^{(1)}|+|\b^{(1)}|)}|F_{\a,\b}
|\xi^{\frac{\a^{(1)} +\b^{(1)}}{2}}\prod_{l=1}^n 2^{ \a^{(1)}_l
+\b^{(1)}_l }  \label{achemenideter} \\
& \leq &
\sum_{\a^{(1)},\b^{(1)}} e^{s( |\a^{(1)}|+|\b^{(1)}|)}|F_{\a,\b}
|(2\rro)^{\frac{|\a^{(1)}| +|\b^{(1)}|}{2}}
2^{|\a^{(1)}|+|\b^{(1)}|} =
\sum_{\a^{(1)},\b^{(1)}} (2 e^s \sqrt{2\rro})^{|\a^{(1)}|+|\b^{(1)}|}
|F_{\a,\b}|\,. \nonumber
\end{eqnarray}
Then, substituting in \eqref{mattonella},
\begin{eqnarray}\label{horus}
\|f\|_{s,2r,\mathcal O_\rro\cup \mathcal O_{2\rro}}
&\leq& \sup_{\|
z\|_{a,p}, \| \bar z \|_{a,p} < 2r} G(z,\bar z)\qquad {\rm where}
\\
G(z,\bar z)
&:=&
\sum_{|\a^{(2)}+\b^{(2)}|\geq d} (2 e^s \sqrt{2\rro})^{|\a^{(1)}|+|\b^{(1)}|}|F_{\a,\b}|
|z^{\a^{(2)}}| |\bar z^{\b^{(2)}}| \,. \label{Gdefa}
\end{eqnarray}
By \eqref{condro},  for all $ \| z\|_{a,p}, \| \bar z\|_{a,p} < 2r $,  the vector $(u^*,\bar u^*)$ defined by
\begin{equation}\label{cambise}
  u_j^*=\bar u_j^*:=
2e^s \sqrt{2\rro} \,,\ \ j\in {\cal I}\,,\quad
u_j^*:= (R/(8r)) |z_j|\,,\
\bar u_j^*:=(R/(8r)) |\bar z_j|\,, \ \ {j\in \Z\setminus \cal I}
\end{equation}
belongs to  $ B_{R/2} \times   B_{R/2} $.
Then, by \eqref{Gdefa}, recalling \eqref{formalpower1A}, Definition \ref{majsc}
(and since  $ R / (8r) > 1 $ by \eqref{condro}),
$$
G (z,\bar z) \leq (8r/R)^d
(M F) (u^*,\bar u^*) \leq (8r/R)^d \| F \|_{R/2}\,, \quad
\forall\,\| z\|_{a,p}, \| \bar z\|_{a,p} < 2r\,.
$$
Hence by \eqref{horus}
\begin{equation}\label{iside}
\|f\|_{s,2r,\mathcal O_\rro\cup \mathcal O_{2\rro}}
\leq
(8r/R)^d
\| F \|_{R/2}\,.
\end{equation}
This shows that $ f $ is $ M $-regular.
Similarly we get
\be\label{compoz}
\|  \partial_{z} f \|_{s, 2r,\mathcal O_\rro\cup \mathcal O_{2\rro}}
\leq
\| \partial_{u^{(2)}} F \|_{R/2}  (8r/R)^{d-1}\, , \ \ {\rm same \ for \ } \partial_{\bar z} \, .
\ee
Moreover,  by the chain rule, and \eqref{iside}
\begin{eqnarray*}
\|  \partial_{x_i} f\|_{s, 2r,\mathcal O_\rro\cup \mathcal O_{2\rro}}
&\leq&
( \| \partial_{u^{(1)}_i} F
\|_{R/2} +\|\partial_{ \bar
u^{(1)}_i}F\|_{R/2})\sqrt{2\rro+\rro/4}e^s (8r/R)^d
\\
 \|  \partial_{y_i} f\|_{s, 2r,\mathcal O_\rro\cup \mathcal O_{2\rro}}
 &\leq&
 ( \| \partial_{u^{(1)}_i} F \|_{R/2} +\|\partial_{ \bar u^{(1)}_i}F\|_{R/2}) \frac{e^s}{\sqrt{\rro/2-\rro/4}} (8r/R)^d \, .
 \end{eqnarray*}
Then \eqref{coord} follows by \eqref{condro} (recalling \eqref{normaEsr}).
\end{pf}

\begin{definition}
For a monomial
$ \mathfrak m_{\a,\b}  :=
(u^{(1)})^{\alpha^{(1)}}(\bar u^{(1)})^{\beta^{(1)}}(u^{(2)})^{\alpha^{(2)}}(\bar u^{(2)})^{\beta^{(2)}} $
$($as in \eqref{defFab}$)$
we set
\be\label{tange}
 {\mathfrak p} ( \mathfrak m_{\a,\b} ):= \sum_{l = 1}^n
 \langle \pluto_l \rangle (\alpha^{(1)}_{\pluto_l} +  \beta^{(1)}_{\pluto_l} ) \,  , \quad \langle j \rangle :=
  \max \{1, |j| \} \, .
\ee
For any  $ F $ as in \eqref{defFab}, $ K \in \N $, we define the projection
\be\label{protang}
\Pi_{\mathfrak p \geq  K} F :=
\sum_{{\mathfrak p} ( \mathfrak m_{\a,\b} ) \geq K} F_{\a,\b}  \mathfrak
m_{\a,\b} \, , \quad \Pi_{\mathfrak p <  K} := I - \Pi_{\mathfrak p \geq  K} \, .
\ee
\end{definition}

\begin{lemma}\label{cometichiami}
Let $ F \in {\cal H}_{R/2} $. Then
\be\label{stimaXP}
\| X_{(\Pi_{\mathfrak p \geq  K}F)\circ
\Phi}\|_{s,r,\mathcal O_\rro}
\leq
2^{- \frac{K}{2\kappa }+1} \|X_{F\circ \Phi}\|_{s,2r,\mathcal O_{2\rro}}\,.
\ee
\end{lemma}

\begin{pf}
For each monomial $ \mathfrak m_{\a,\b} $ as in \eqref{defFab}
with $ \mathfrak p ( \mathfrak m_{\a,\b} ) \geq  K $  we have
$$
| \alpha^{(1)}+\beta^{(1)} | \stackrel{\eqref{defa1b1}} = \sum_{l=1}^n
 \alpha^{(1)}_{\pluto_l } + \beta^{(1)}_{\pluto_l }
\stackrel{\eqref{caracalla}} \geq
\kappa ^{-1} \sum_{l=1}^n  \langle \pluto_l \rangle
(  \alpha^{(1)}_{\pluto_l } + \beta^{(1)}_{\pluto_l } )
\stackrel{\eqref{tange}} =
\kappa ^{-1}   \mathfrak p ( \mathfrak m_{\a,\b} ) \geq
\kappa ^{-1} K
$$
and then,  $ \forall \xi \in \mathcal O_\rro $, $ y \in B_{r^2} $,
\begin{eqnarray}
|(\mathfrak m_{\a,\b} \circ \Phi ) (x,y,z, \bar z; \xi)| & \stackrel{\eqref{actionangle}} = &
| (\xi + y)^{\frac{\alpha^{(1)} + \beta^{(1)}}{2}}
 e^{\ii (\alpha^{(1)}-\beta^{(1)})\cdot x} z^{\alpha^{(2)}} \bar
z^{\beta^{(2)}}|  \\
& = &
2^{- \frac{|\alpha^{(1)}+\beta^{(1)}|}{2}}|(2\xi+2y)^{\frac{\alpha^{(1)}+\beta^{(1)}}{2}}
e^{\ii (\alpha^{(1)}-\beta^{(1)})\cdot x} z^{\alpha^{(2)}} \bar
z^{\beta^{(2)}}| \nonumber \\
& \leq &
2^{- \frac{K}{2 \kappa }} |(\mathfrak m_{\a,\b} \circ \Phi ) (x,2y,z,\bar z;2\xi)| \, . \nonumber
\end{eqnarray}
The bound \eqref{stimaXP} for the Hamiltonian vector field
follows applying the above rescaling argument to each component,
and noting that the derivatives with respect to $y$  in the
vector field  
 decrease the degree in $ \xi $ by one. 
\end{pf}

Let $ N_0, \teta, \mu,\mu'$ be as in Proposition \ref{qtop0}.
For $ N \geq N_0 $  and $ F \in {\cal H}_{R/2} $ we set
\begin{equation}\label{adamo}
 f^*:= \Pi_{N,\theta,\mu}\Big( (F-\Pi_{N,\theta,\mu' }F)\circ
\Phi \Big) \, .
\end{equation}
Note that $\Pi_{N,\theta,\mu'}$ is the projection on the
bilinear functions in the variables $u,\bar u$,
while $\Pi_{N,\theta,\mu}$
in the variables $x,y,z,\bar z $.

\begin{lemma}\label{ciccio} We have
\begin{equation}\label{cojo2}
\| X_{ f^* }\|_{s,r,\mathcal O_{\rro}}
 \leq
2^{- \frac{N^b}{2 \kappa} + 1 }  \|X_{F\circ\Phi}\|_{s,2r,\mathcal O_{2\rro}} \, .
\end{equation}

\end{lemma}
\begin{pf}
We first claim that
if $ F =  \mathfrak m_{\a,\b} $  is  a  monomial as in \eqref{defFab} with
$ \mathfrak p( \mathfrak m_{\a,\b} ) < N^b $ then  
$ f^* = 0 $.

{\sc Case $ 1 $}: $  \mathfrak m_{\a,\b} $ is $( N,\theta,\mu')$--bilinear, see Definition \ref{BL}. Then
$ \Pi_{N,\theta,\mu' }  \mathfrak m_{\a,\b} =  \mathfrak m_{\a,\b} $ and $ f^*= 0 $, see \eqref{adamo}.

{\sc Case $ 2 $}:  $  \mathfrak m_{\a,\b} $ is {\em not }$(N,\theta,\mu')$--bilinear.
Then $\Pi_{N,\theta,\mu' }  \mathfrak m_{\a,\b} = 0 $ and $ f^*=\Pi_{N,\theta,\mu}(  \mathfrak m_{\a,\b} \circ \Phi) $,
see \eqref{adamo}.
We claim that $ \mathfrak m_{\a,\b} \circ \Phi $
is not $(N,\theta,\mu)$--bilinear, and so $  f^* = \Pi_{N,\theta,\mu}(  \mathfrak m_{\a,\b} \circ \Phi) = 0 $.
Indeed,
\be\label{mabco}
\mathfrak m_{\a,\b} \circ \Phi = (\xi + y)^{\frac{\alpha^{(1)} + \beta^{(1)}}{2}}
 e^{\ii (\alpha^{(1)}-\beta^{(1)})\cdot x} z^{\alpha^{(2)}} \bar
z^{\beta^{(2)}}
\ee
is  $(N,\theta,\mu)$--bilinear if and only if  (see  Definitions \ref{BL} and \ref{LM})
$$
z^{\alpha^{(2)}} {\bar z}^{\beta^{(2)}}=
z^{\tilde \alpha^{(2)}} {\bar z }^{\tilde \beta^{(2)}} z_m^\sigma
z_n^{\sigma '} \, ,
$$
\begin{equation}\label{rompic}
\sum_{j\in \Z\setminus \cal I} |j| (\tilde \alpha^{(2)}_j+ \tilde \beta^{(2)}_j)<
\mu N^L \,,\quad  |m|,|n|> \theta N\,,
\quad
|\a^{(1)}-\b^{(1)}|< N^b\,.
\end{equation}
We deduce the contradiction that 
$ \mathfrak m_{\a,\b} =
(u^{(1)})^{\alpha^{(1)}}(\bar u^{(1)})^{\beta^{(1)}}
(u^{(2)})^{\tilde \alpha^{(2)}}(\bar u^{(2)})^{\tilde \beta^{(2)}} u_m^\sigma u_n^{\sigma '} $
is  $(N,\theta,\mu')$-bilinear 
because (recall that we suppose $ \mathfrak p(  \mathfrak m_{\a,\b} ) < N^b $)
$$
\sum_{l =1}^n |\pluto_l |(\alpha^{(1)}_{\pluto_l}
+ \beta^{(1)}_{\pluto_l})+ \sum_{j\in \Z\setminus
\cal I} |j| (\tilde \alpha^{(2)}_j +  \tilde \beta^{(2)}_j)\!\! \stackrel{\eqref{tange}, \eqref{rompic}}{<} \!\!
\mathfrak p(  \mathfrak m_{\a,\b} )+ \mu N^L < N^b + \mu N^L
\stackrel{\eqref{fuso}} < \mu' N^L \, .
$$
For the general case, we divide $ F= \Pi_{\mathfrak p <  N^b} F + \Pi_{\mathfrak p \geq   N^b}F
$. By the above claim
$$
 f^*=\Pi_{N,\theta,\mu}\Big(
\big((Id-\Pi_{N,\theta,\mu'})\Pi_{\mathfrak p \geq   N^b} F\big)\circ
\Phi\Big)=\Pi_{N,\theta,\mu}\Big(
\big(\Pi_{\mathfrak p \geq   N^b}(Id-\Pi_{N,\theta,\mu'}) F\big)\circ
\Phi\Big) \, .
$$
Finally, \eqref{cojo2} follows by \eqref{proiezaa} and applying Lemma \ref{cometichiami}
to $ \big(\Pi_{\mathfrak p \geq   N^b}(Id-\Pi_{N,\theta,\mu'}) F\big)\circ
\Phi $.
\end{pf}
\begin{lemma}\label{top:aa}
Let 
 $ F\in {\cal T}_{R/2}(N,\theta,\mu')$ 
with $\Pi_{\mathfrak p \geq   N^b}F=0 $. Then $ F \circ \Phi(\cdot;\xi)\in {\cal T}_{s,
2r}(N,\theta,\mu') $, $\forall\, \xi \in\mathcal O_{\rro}\cup \mathcal O_{2\rro}$ .
\end{lemma}
\begin{pf}
Recalling Definition \ref{matteo_aa} we have
$$
F = \sum_{|m|,|n|> \teta N, \s,\s'=\pm}
F^{\s,\s'}(\mathtt s(m),\s m+\s' n) u_m^\s u_n^{\s'} \ \ {\rm with} \ \ F^{\s,\s'}(\varsigma , h)
\in {\cal L}_{R/2}(N,\mu', h ) \, .
$$
Composing with the map $ \Phi $ in \eqref{actionangle}, since $ m, n \notin {\cal I } $,  we get
$$
F\circ \Phi = \sum_{\s,\s'=\pm\,, |m|,|n|> \teta N}
F^{\s,\s'}(\mathtt s(m),\s m+\s' n)\circ \Phi \, z_m^\s z_n^{\s'} \, .
$$
Each coefficient
$ F^{\s,\s'}(\mathtt s(m),\s m+\s' n)\circ \Phi$  depends on $n,m,\s,\s'$ only
through $\mathtt s(m),\s m+\s' n,\s,\s' $. Hence, in order to conclude that  $ F\circ \Phi
 \in {\cal T}_{s, 2r}(N,\theta,\mu')
$
 it remains only to prove
 that  $ F^{\s,\s'}(\mathtt s(m),\s m+\s' n) \circ\Phi\in {\cal L}_{s,2r}(N,\mu',\s m+\s' n) $,
 see Definition \ref{LM}.
Each monomial $ \mathfrak m_{\a,\b} $ of $F^{\s,\s'}( \mathtt s(m) , \s m+\s' n ) \in {\cal L}_{R/2}(N,\mu',\s m+\s' n) $
satisfies
$$
\sum_{l=1}^n (\a_{\pluto_l}+\b_{\pluto_l} )|\pluto_l |
+ \sum_{j\in \Z\setminus\mathcal I}(\a_j+\b_j)|j| < \mu' N^L  \ \quad
{\rm and } \quad \ \mathfrak p(\mathfrak m_{\a,\b} ) < N^b
$$
by the hypothesis $ \Pi_{\mathfrak p \geq   N^b}F =0 $. 
Hence 
 $ \mathfrak m_{\a,\b} \circ \Phi$ (see \eqref{mabco}) is $(N,\mu')$-low momentum,
 in particular $  | \a^{(1)} - \b^{(1)} | \leq
\mathfrak p(\mathfrak m_{\a,\b} ) < N^b $.
\end{pf}

\smallskip

\begin{pfn}{\sc of Proposition \ref{qtop0}}.
Since $  F \in {\cal Q}_{R/2}^T(N_0,\theta,\mu' ) $ (see Definition \ref{topbis_aa}), for all $ N\geq N_0 $,
there is a T\"oplitz approximation  $ \tilde F \in {\cal
T}_{R/2}(N,\theta,\mu' ) $  of $ F $, namely
\be\label{Pmu1}
\Pi_{N, \teta, \mu'} F= \tilde F+N^{-1}\hat F \,\quad{\rm  with }
\quad   \|X_F \|_{R/2}, \|X_{\tilde F}\|_{R/2},\|X_{\hat F}\|_{R/2} <
2\|F\|_{R/2, N_0 ,\theta,\mu'}^T \, .
\ee
In order to prove that $ f := F \circ \Phi \in  {\cal Q}_{s,r}^T(N_0,\theta,\mu)  $  we define
its candidate T\"oplitz approximation
\be\label{tildefap}
\tilde f:=
\Pi_{N,\theta,\mu}( (\Pi_{ {\mathfrak p} < N^b} \tilde F)
\circ \Phi) \, ,
\ee
see \eqref{protang}.
Lemma \ref{top:aa} applied to
$  \Pi_{ {\mathfrak p} < N^b} \tilde F \in  {\cal T}_{R/2}(N,\theta,\mu' )  $ 
implies that
$ (\Pi_{ {\mathfrak p} < N^b} \tilde F) \circ \Phi  \in {\cal T}_{s,2r}(N,\theta,\mu') $ and
then, applying the projection $ \Pi_{N,\teta, \mu}$ we get
$ \tilde f \in {\cal T}_{s,2r}(N,\theta,\mu ) \subset  {\cal T}_{s,r}(N,\theta,\mu ) $.
Moreover,  by \eqref{tildefap} and applying Lemma \ref{artaserse} to
$ \Pi_{ {\mathfrak p} < N^b} \tilde F $ (note that $  \Pi_{ {\mathfrak p} < N^b} \tilde F $ is either zero or it
is in $  {\cal H}^{d}_{R/2} $ with $ d \geq 2 $ because it is bilinear),
we get 
\begin{eqnarray}\label{hofame}
\|X_{\tilde f}\|_{s,r,\mathcal O_\rro} \stackrel{\eqref{proiezaa}} \leq  \|X_{ (\Pi_{ {\mathfrak p} < N^b} \tilde F)
\circ \Phi) }\|_{s,r,\mathcal O_\rro}
& \stackrel{\eqref{coord}} \lessdot & (8r/R)^{d-2}
\|X_{\Pi_{ {\mathfrak p} < N^b} \tilde F} \|_{R/2} \nonumber \\
&  \stackrel{\eqref{proiezaa}, \eqref{Pmu1}} \lessdot  & (8r/R)^{d-2} \| F \|_{R/2, N_0, \teta, \mu' }^T \, .
\end{eqnarray}
Moreover the T\"oplitz defect is
\begin{eqnarray*}
\hat f
&:=&
N(\Pi_{N,\theta,\mu} f-\tilde f)
\stackrel{\eqref{tildefap}} = N \, \Pi_{ N,\theta,\mu }\big(( F - \Pi_{{\mathfrak p}< N^b} \tilde F)\circ \Phi\big)
\\
&=&
 N\Pi_{N,\theta,\mu}\big((F-\tilde F)\circ \Phi\big) + N \Pi_{ N,\theta,\mu }\big(( \tilde F - \Pi_{{\mathfrak p} < N^b} \tilde F)\circ \Phi\big)
\\
& \stackrel{\eqref{Pmu1}, \eqref{protang}} =&
\Pi_{ N,\theta,\mu }(\hat F\circ \Phi)
+ N\Pi_{N,\theta,\mu} \Big( \big(F-\Pi_{N,\theta,\mu'}F\big)\circ \Phi\Big)
+ N\Pi_{ N,\theta,\mu
} \big((\Pi_{{\mathfrak p} \geq N^b} \tilde F)\circ \Phi\big)
\\
&\stackrel{\eqref{adamo}}=&
\Pi_{ N,\theta,\mu }(\hat F\circ \Phi)
+ Nf^*
+ N\Pi_{ N,\theta,\mu
} \big((\Pi_{{\mathfrak p} \geq N^b} \tilde F)\circ \Phi\big)
\,.
\end{eqnarray*}
Using \eqref{proiezaa}, Lemmata \ref{cometichiami} and
\ref{ciccio}
 imply that, since $N 2^{-\frac{N^b}{2\kappa}+1}\leq 1,$
  $\forall\, N\geq N_0$  by \eqref{fuso},
\begin{eqnarray}
\|X_{\hat f}\|_{s,r,\mathcal O_{\rro}}
&\leq&
\|X_{\hat F\circ \Phi}\|_{s,r,\mathcal O_{\rro}} + N
2^{-\frac{N^b}{2\kappa}+1}
(\|X_{F\circ \Phi}\|_{s,2r,\mathcal O_{2\rro}}
+\|X_{ \tilde F\circ \Phi}\|_{s,2r,\mathcal O_{2\rro}})
\nonumber
 \\
 &\lessdot &
  \|X_{\hat F\circ \Phi}\|_{s,2r,\mathcal O_{\rro}}+
 \|X_{ F\circ \Phi}\|_{s,2r,\mathcal O_{2\rro}}
 +\|X_{ \tilde F\circ \Phi}\|_{s,2r,\mathcal O_{2\rro}} \nonumber
  \\
&  \stackrel{\eqref{coord}}{\lessdot} &
(8r/R)^{d-2}
(\|X_{\hat F}\|_{R/2}+\|X_{F}\|_{R/2}+\|X_{ \tilde F}\|_{R/2})
\label{ultimopassa}
  \\
& \stackrel{\eqref{Pmu1}} \lessdot & (8r/R)^{d-2} \|F\|_{R/2,N_0,\theta,\mu'}^T
 \label{hofamebis}
 \end{eqnarray}
(to get  \eqref{ultimopassa}  we also note that $ F, \hat F , \tilde F \in {\cal H}^{d}_{R/2} $ with $ d = 0, 1 $, unless are zero).

The bound \eqref{topstim} follows by \eqref{coord}, \eqref{hofame}, \eqref{hofamebis}.
\end{pfn}

\noindent
We conclude this subsection with a lemma,
similar to Lemma \ref{artaserse},
used in Lemma \ref{dario} (see \eqref{coma2}). 

\begin{lemma}\label{artasersebis}
 Let $ F \in {\cal H}_{R/2} $, $ f := F
\circ \Phi$  and $\tilde f(x,y):=f(x,y,0,0)-f(x,0,0,0).$
Then, assuming \eqref{condro},
\begin{equation}\label{coordbis}
\|X_{\tilde f}\|_{s,2r,\mathcal O_\rro\cup \mathcal O_{2\rro}}
\lessdot
 \|X_F\|_{R/2} \, .
\end{equation}
Moreover if $F$ preserves momentum then so does $\tilde f$.
\end{lemma}
\begin{pf}
We proceed as in Lemma \ref{artaserse}.
The main difference is that here there are no $(z,\bar z)$-variables
and the sum in \eqref{achemenide} runs over $i\neq 0$.
Then in  the product in \eqref{achemenidebis}
(at least) one of the sums  is on $i_l\geq 1$.
Therefore we can use the second estimate
in \eqref{embeh}
gaining a factor\footnote{Actually we have the constant 3 instead of 2 in \eqref{achemenideter}
and $ 3 e^s $ instead of $ 2 e^s $ in \eqref{horus} and \eqref{cambise}.}
$8r^2/\rro$ (since $|y_l|/|\xi_l|\leq 8r^2/\rro$
by \eqref{paraO}).
Continuing as in the proof of Lemma \ref{artaserse}
we get (recall \eqref{coord} with $d=0$)
$$
\|X_{\tilde f}\|_{s,2r,\mathcal O_\rro\cup \mathcal O_{2\rro}}
\lessdot (r^2/\rro) (r/R)^{-2}  \|X_F\|_{R/2} \stackrel{ \eqref{condro}} \lessdot  \|X_F\|_{R/2}
$$
proving \eqref{coordbis}. 
\end{pf}

\subsection{Proof of Theorem \ref{thm:DNLW}} \label{sec:pr}

We now introduce the action-angle variables \eqref{actionangle} (via the map \eqref{variableAA})
in the Birkhoff normal form Hamiltonian \eqref{Birkhoff}.
Hence we obtain the parameter dependent family of Hamiltonians
\be\label{def:H'}
H' := H_{\rm Birkhoff} \circ \Phi = {\cal N} + P
\ee
where (up to a constant), by \eqref{giulio},
\be\label{NP}
{\cal N } := \om (\xi) \cdot y + \Om (\xi) z \bar z \, , \ \
P:=  \frac12 Ay \cdot y + By \cdot z \bar z + \hat G (z, \bar z) + K' (x,y,z, \bar z; \xi) \, ,
\ee
\be\label{frequenzeshi}
\om (\xi ) :=  \bar \om + A \xi \, , \
\bar \om := (\l_{\pluto_1}, \ldots, \l_{\pluto_n}) \, , \quad
\Om (\xi ) :=  \bar \Om + B \xi \, , \
\bar \Om := (\l_{j})_{j     \in \Z \setminus {\cal I}} \, ,
\ee
\be \label{K'def}
A = (A_{l h})_{1 \leq l,h \leq n} \, , \ A_{lh} := 12 ( 2 - \d_{lh} ) \, , \
B = (B_{jl})_{j \in \Z \setminus {\cal I}, 1 \leq l \leq n} \, , \ B_{jl} := 24 \, ,
\ K' := K \circ \Phi \, .
\ee
The parameters $ \xi $ stay in the set $ \mathcal O_\rro $ defined in \eqref{paraO} with
$ \rho = C_* R^2 $ as  in \eqref{condro}.
As in \eqref{P0Pbar} we decompose the perturbation 
\begin{equation}\label{serse}
P  = P_{00} + \bar P  \quad {\rm where} \quad
P_{00}(x; \xi) := K'(x,0,0,0;\xi) \, ,
\quad \bar P := P - P_{00} \, .
\end{equation}
\begin{lemma}\label{dario}
Let $s,r>0$ as in \eqref{condro} and $N$
large enough (w.r.t. $\mm,{\cal I},L,b$). Then
\begin{equation}\label{milziade}
\| X_{P_{00}} \|_{s,r}  \lessdot R^6 r^{-2} \,,\qquad
\| \bar P \|_{s,r,N,\CCterzi, \duecc}^T
 \lessdot   r^2 + R^5 r^{-1}
\end{equation}
and, for $\l>0,$
\begin{equation}\label{milziadebis}
| X_{P_{00}} |_{s,r}^\l  \lessdot (1+\l/\rro) R^6 r^{-2} \,,\qquad
| X_{\bar P} |_{s,r}^\l
 \lessdot (1+\l/\rro) (r^2 + R^5 r^{-1})\,,
\end{equation}
for $ \xi $ belonging to
\begin{equation}\label{temistocle}
    \mathcal{O}(\rro):=
   \Big\{ \xi\in \R^n\  :\  \frac23 \rro \leq \xi_l\leq\frac34 \rro\,,
   \ \  l=1,\ldots , n   \Big\} \subset {\cal O}_\rho \,.
\end{equation}
\end{lemma}

\begin{pf}
By the definition \eqref{serse} we have
\begin{eqnarray}
\| X_{P_{00}} \|_{s,r}
 \stackrel{\eqref{caligola}} \leq   \|  X_{K'} \|_{s,r}
\stackrel{\eqref{maggio}}
\leq
  \|  K' \|_{s,r,N,\CCterzi,\duecc}^T
  & \stackrel{\eqref{K'def}} = &
\|  K \circ \Phi \|_{s,r,N,\CCterzi,\duecc}^T \nonumber \\
& \stackrel{\eqref{topstim}}  \lessdot &
\Big(\frac{r}{R} \Big)^{-2}\|  K  \|_{R/2, N , \CCterzi,\duecc}^T
\label{p00}
\end{eqnarray}
(applying \eqref{topstim} with $ d \rightsquigarrow 0 $,
$ N_0 \rightsquigarrow N $, $ \teta  \rightsquigarrow \CCterzi $,
$ \mu  \rightsquigarrow  \duecc $,
 $ \mu'  \rightsquigarrow  \trecc $) and
taking $ N $ large enough so that \eqref{fuso} holds.
Take also $ N \geq N_0' $ defined in Proposition \ref{BNF}. Then
by \eqref{p00} we get
$$
\| X_{P_{00}} \|_{s,r} \stackrel{\eqref{inscatola}} \lessdot
\Big(\frac{r}{R} \Big)^{-2}\|  K  \|_{R/2, N_0' , \CCterzi,\duecc}^T
\stackrel{\eqref{stellamarisbis}}
\lessdot
\Big(\frac{r}{R} \Big)^{-2} R^4
\lessdot \frac{R^6}{r^2}
$$
proving the first estimate in \eqref{milziade}.
Let us prove the second bound.
By \eqref{serse} and \eqref{NP} we write
\be\label{barP}
\bar P = \frac12 Ay \cdot y + By \cdot z \bar z +
 \hat G (z, \bar z) + K_1 + K_2
\ee
where
$$
K_1 :=  K' (x,y,z, \bar z; \xi) - K' (x,y,0,0; \xi) \, ,
\quad  K_2 :=  K' (x,y,0,0; \xi) -  K' (x,0,0, 0; \xi) \, .
$$
Using \eqref{stellamaris} (note that $r<R$ by \eqref{condro})
for $ N \geq N_0 $ large enough to fulfill \eqref{caracalla},
 we have by \eqref{inscatola}
\be\label{coma1}
\Big\| \frac12 Ay \cdot y + By \cdot z \bar z + \hat G (z, \bar z)
\Big\|_{s,r,N,\CCterzi,\duecc}^T \lessdot r^2
\, .
\ee
By  \eqref{topstim} (with $ d \rightsquigarrow 1 $,
$ N_0 \rightsquigarrow N $,
$ \mu  \rightsquigarrow  \duecc$,
 $ \mu'  \rightsquigarrow  \trecc $),
for $ N \geq N_0 (\mm,{\cal I},L,b)$ large enough,  we get
\be\label{coma3}
\| K_1 \|_{s,r,N,\CCterzi,\duecc}^T
\lessdot
\Big(\frac{r}{R} \Big)^{-1} R^4  \lessdot \frac{R^5}{r} \, .
\ee
Moreover, since $ K_2 $ does not depend on $ (z, \bar z) $, we have
\be\label{coma2}
\| K_2 \|_{s,r,N,\CCterzi,\duecc}^T
\stackrel{\eqref{clearly}} =
\| X_{K_2} \|_{s,r}
 \stackrel{\eqref{coordbis}} \lessdot
\|X_{K}\|_{R/2}  \stackrel{\eqref{maggio}}  \lessdot
\| K\|_{R/2, N_0', \CCterzi,\trecc}^T
 \stackrel{\eqref{stellamarisbis}} \lessdot R^4 \, .
\ee
In conclusion,  \eqref{barP}, \eqref{coma1}, \eqref{coma3}, \eqref{coma2}  imply the second estimate in \eqref{milziade}:
$$
\| \bar P \|_{s,r,N,\CCterzi,\duecc}^T
 \lessdot  r^2 +  \frac{R^5}{r} + R^4   \stackrel{\eqref{condro}}
 \lessdot  r^2 + \frac{R^5}{r}  \, .
$$
Let us prove the estimates  \eqref{milziadebis}
for the Lipschitz norm  defined in \eqref{filadelfia}
 (which involves only 
the  sup-norm of the vector fields).
First 
$$
|X_{P_{00}} |_{s,r}  \stackrel{\eqref{XMXA}} \leq  \|X_{P_{00}} \|_{s,r}
  \stackrel{\eqref{milziade}}  \lessdot R^6 r^{-2} \, , \ \
|X_{\bar P} |_{s,r}  \stackrel{\eqref{XMXA}} \leq  \| X_{\bar P} \|_{s,r}
\stackrel{\eqref{maggio}} \leq  \| \bar P \|_{s,r,N,\CCterzi,\duecc}^T
 \stackrel{\eqref{milziade}} \lessdot   r^2 + R^5 r^{-1}
 $$
Next, 
since the vector fields $ X_{P_{00}}  $, $ X_{\bar P}  $ are \textsl{analytic}
in the parameters $\xi\in\mathcal{O}_\rro $,
Cauchy estimates in the domain
$\mathcal{O}(\rro)\subset \mathcal{O}_\rro $ (see \eqref{temistocle}) 
imply
$$
| {X}_{P_{00}} |^{{\rm lip}}_{s,r, {\cal O}(\rro)} \lessdot
\rho^{-1} | {X}_{P_{00}} |_{s,r, {\cal O}_\rho} \lessdot  R^6
r^{-2} , \quad | {X}_{\bar P} |^{{\rm lip}}_{s,r, {\cal O}(\rro)}
\lessdot \rho^{-1} | {X}_{\bar P} |_{s,r, {\cal O}_\rho}
 \lessdot   r^2 + R^5 r^{-1}
$$
and \eqref{milziadebis} are proved.
\end{pf}

\smallskip

All the assumptions of Theorems
 \ref{thm:IBKAM}-\ref{thm:measure} are fulfilled by $ H' $ in \eqref{def:H'} with parameters
 $ \xi \in {\cal O}(\rro) $ defined in \eqref{temistocle}.
 Note that the sets ${\cal O}=[\rho/2,\rho]^n $
 defined in Theorem \ref{thm:measure}
 and ${\cal O}(\rro)$ defined in \eqref{temistocle}
 are diffeomorphic through $\xi_i\mapsto (7\rho +2\xi_i)/12$.
The hypothesis (A1)-(A2) follow from \eqref{frequenzeshi}, \eqref{K'def}
with
$$
a(\xi) = 
24 \sum_{l=1,\ldots, n} \xi_l \, , 
\quad \
{\rm and} \ \quad M_0 = 24 + \| A \|  \, .
$$
Then (A3)-(A4) and the quantitative bound \eqref{KAMconditionT} follow by
\eqref{milziade}-\eqref{milziadebis},
choosing 
\begin{equation}\label{leonida}
s = 1,\  r = R^{1 + \frac34}\,,
\
\rro = C_* R^2 \  {\rm as \ in \ }\eqref{condro},
\
N \ {\rm as \ in \ Lemma \ }\ref{dario},\
  \theta = \CCterzi, \
   \mu = \duecc,
   \
\gamma = R^{3 + \frac15}   \, ,
\end{equation}
and taking $ R $ \textsl{small enough}. Hence Theorem  \ref{thm:IBKAM} applies.

Let us verify that also the assumptions  of Theorem \ref{thm:measure} are fulfilled.
Indeed \eqref{omegaxi} follows by 
\eqref{frequenzeshi}, \eqref{K'def} with $ \vec  a = 24 (1, \ldots, 1 ) \in \R^n $.
The matrix $ A $ defined in \eqref{K'def} is invertible and
$$
A^{-1} = (A^{-1}_{lh})_{1 \leq l,h \leq n} \, , \quad
A^{-1}_{lh} = \frac{1}{12} \Big( \frac{2}{2n - 1} - \d_{lh} \Big) \, .    
$$
Finally the non-degeneracy assumption \eqref{pota} is satisfied because $ A = A^T $ and
$$
2 A^{-1} \vec  a = \frac{4}{2n-1} (1, \ldots, 1) \notin \Z^n \setminus {0} \, .
$$
We deduce that the Cantor  set of parameters
$ \mathcal{O}_\infty \subset \mathcal{O} $ in \eqref{Cantorinf} has asymptotically full density
because
$$
\frac{|\mathcal{O} \setminus\mathcal{O}_\infty|}{| \mathcal{O}|}
\stackrel{\eqref{consolatrixafflictorum}} \lessdot
\rho^{-1} \g^{2/3}\stackrel{\eqref{leonida}}\lessdot R^{-2}
R^{\frac23 (3+\frac15)}=R^{\frac{2}{15}} \to 0 \, . 
$$

The proof of Theorem \ref{thm:DNLW} is now completed.

\begin{remark}\label{carmelo}
The  terms $\sum_{k\geq 5}f_k s^k$ in \eqref{funzi}
 contribute to the Hamiltonian \eqref{Ham}  with monomials of order 6 or more and \eqref{Birkhoff} holds
 (with a possibly different $ K $ satisfying \eqref{stellamarisbis}).
On the contrary, the term $ f_4 s^4 $  in  \eqref{funzi} would add monomials of order $5$
to the Hamiltonian in \eqref{Ham}. Hence \eqref{stellamarisbis} holds  with $R^3$
instead of $R^4$. This estimate is not sufficient. These $5$-th order terms should
be removed by a  Birkhoff normal form.
For simplicity,  we did not pursue this point.
\end{remark}

\section{Appendix}\setcounter{equation}{0}\label{sec:appendix}

{\sc Proof of Lemma \ref{Cauchyx}.}
We need some notation: we write
$ E = \oplus_{j=1}^4 E_j $, $ E_1 := (\C^n,|\ |_\infty) $, $ E_2 := (\C^n,|\ |_1) $, $ E_3 :=E_4:=  \ell^{a,p}_{\cal I}$
so that a vector $v=(x,y,z,\bar z)\in E$ can be expressed by its
four components $v^{(j)}\in E_j$,
$v^{(1)}:=x,$ $v^{(2)}:=y,$ $v^{(3)}:=z,v^{(4)}:=\bar z,$
 and the norm \eqref{normaEsr}
is \be\label{norma4} \| v\|_{E,s,r}:= \sum_{j=1}^4\frac{
|v^{(j)}|_{E_j}}{\rho_j}\,,\quad{\rm where}\;  \rho_1=s\,,\quad
\rho_2= r^2\,,\quad \rho_3=\rho_4= r\,. \ee
We are now ready to prove \eqref{nitore}. By definition
\begin{eqnarray*}
 \| d X(v) \|_{\mathcal L((E,s,r); (E,s',r'))} \! \! \! \!  & :=  & \! \! \! \!  \sup_{\|Y\|_{E, s,r}\leq 1}   \|  d X(v)[Y] \|_{E,s',r'}
\! \!   \stackrel{\eqref{norma4}}  = \! \!   \sup_{\|Y\|_{E,s,r}\leq 1}  \sum_{i=1}^4 \ \frac{|  d X^{(i)}(v)[Y] |_{E_i}}{\rho_i'} \\
 & = & \! \! \! \!    \sup_{\|Y\|_{E,s,r}\leq 1}\sum_{i=1}^4
 \frac{|\sum_{j=1}^4 d_{v^{(j)}}X^{(i)}(v)Y^{(j)}|_{E_i}}{\rho'_i} \\
& \leq & \! \!  \! \!  \sup_{\|Y\|_{E,s,r}\leq 1}\sum_{i,j=1}^4 \frac{ |d_{v^{(j)}}X^{(i)}(v)Y^{(j)}|_{E_i}}{\rho'_i} \\
& \leq & \! \! \! \!   \sup_{\|Y\|_{E,s,r}\leq 1}\sum_{i,j=1}^4 \frac{1}{\rho'_i}
 \|d_{v^{(j)}}X^{(i)}(v)\|_{{\cal L}(E_j,E_i)} |Y^{(j)}|_{E_j} \\
& \leq & \! \!  \! \!   \sup_{\|Y\|_{E,s,r}\leq 1} \sup_{\tilde v \in D(s,r)}  \sum_{i,j=1}^4
\frac{1}{\rho'_i} \frac{|X^{(i)}(\tilde v)|_{E_i}}{(\rho_j-\rho'_j)} |Y^{(j)}|_{E_j}
\end{eqnarray*}
by the Cauchy estimates in Banach spaces. Then
 \begin{eqnarray*}
\| d X(v) \|_{\mathcal L((E,s,r); (E,s',r'))}  \!\!\! & \leq& \!\!\!
 \sup_{\tilde v \in D(s,r)}\sum_{i=1}^4 \frac{\rho_i}{\rho'_i} \frac{|X^{(i)}(\tilde v)|_{E_i}}{\rho_i}\sup_{\|Y\|_{E,s,r}\leq 1}\sum_{j=1}^4
 \Big( 1-\frac{\rho'_j}{\rho_j} \Big)^{-1}\frac{|Y^{(j)}|_{E_j}}{\rho_j} \\
 & \stackrel{\eqref{norma4}} \leq & \!\! \! \max_{i=1,\ldots,4} \frac{\rho_i}{\rho'_i} \,
\max_{j=1,\ldots, 4 } \Big( 1-\frac{\rho'_j}{\rho_j} \Big)^{-1}
\sup_{\tilde v \in D(s,r)} \| X (\tilde v) \|_{E,s,r} \leq 4
\delta^{-1} | X |_{s,r}
\end{eqnarray*}
by \eqref{tiatira}, \eqref{diffusivumsui}. This proves \eqref{nitore}.
\rule{2mm}{2mm}

\medskip

\textsc{Proof of Lemma \ref{poeschel}.}
We first extend
Lemma 4  of \cite{Po3} proving that:
\begin{lemma}\label{Poschelzero}
If $0 \leq i\leq j\leq k\leq l$ with $i\pm j\pm k\pm l=0$ for {\sc some}
combination of plus and minus signs and $(i,j,k,l)\neq(p,p,q,q)$ for $ p, q \in \N $,
then, there exists an absolute constant $ c >  0 $, such that
\be\label{dalbasso}
|\pm \l_i(\mm) \pm \l_j(\mm) \pm \l_k(\mm) \pm \l_l(\mm)|\geq c\mm(i^2+ \mm)^{-3/2}
\ee
for {\sc all} possible combinations  of plus and minus signs
\end{lemma}

\begin{pf}
When $ i >  0 $ it is a reformulation of the statement of
Lemma 4 of \cite{Po3}.  Let us prove it also for $ i = 0 $.
Then $ j \pm k\pm l=0$
for some combination of plus and minus signs.  
Since $ (i,j,k,l) \neq (0,0,q,q) $,
the only possibility
is $l=j+k$ with $ j \geq 1$ (otherwise $i=j=0$ and $k=l$). We have to study
$$
\d (\mm) := \pm \l_0 (\mm)  \pm \l_j (\mm)  \pm \l_k (\mm) \pm \l_l (\mm)
$$
for all possible combinations  of plus and minus signs. To this end, we distinguish
them according to their number of plus and minus signs.
To shorten notation we let, for example, $\d_{++-+}=\l_0+\l_j-\l_k+\l_l$, similarly
for the other combinations.
The only interesting cases are when there are one or two minus signs.
The case when there are no  (or four) minus signs is trivial. When there are $ 3$ minus signs we reduce to the case with one minus sign by a global sign change.
\\[1mm]
{\sl One minus sign.}
Since $\d_{++-+},\d_{+-++},\d_{-+++}\geq \d_{+++-}:=\d$
we study only  the last case.
We have 
$$
\d(0)=j+k-l=0\,,\qquad \d'(\mm)=\frac12 \left( \frac{1}{\l_0}+
\frac{1}{\l_j}+\frac{1}{\l_k}-\frac{1}{\l_l} \right)\geq \frac{1}{2\l_0}
=\frac{1}{2\sqrt \mm}\,.
$$
Therefore $\d(\mm)\geq \sqrt \mm\geq c\mm (1+\mm)^{-3/2} $ for an absolute constant $ c >  0 $.
\\[1mm]
{\sl Two minus signs.}
Now we have $\d_{-+-+},\d_{--++}\geq \d_{+--+}$
and all other cases reduce to these ones by inverting signs.
So we consider only $\d=\d_{+--+}.$
Since the function $f(t):=\sqrt{t^2+\mm}$ is monotone increasing and convex for $t\geq 0$,
we have the estimate
\begin{equation}\label{fave}
\l_l - \l_k \geq \l_{l-p}-\l_{k-p} \, , \quad \forall \, 0\leq p\leq k \,.
\end{equation}
Hence $\l_l-\l_k\geq \l_{j+1}-\l_1$
and $\l_{j+1}-\l_j\geq \l_2-\l_1$ (using $j=l-k\geq 1$).
Therefore
$$
\d=\l_0-\l_j-\l_k+\l_l\geq \l_0-\l_j-\l_1+\l_{j+1}
\geq
\l_2-2\l_1+\l_0\geq \mm (4+\mm)^{-3/2}\,.
$$
The last inequality follows since $f''(t)=\mm (t^2+\mm)^{-3/2} $
is decreasing and $ \l_2-2\l_1+\l_0=f(2)-2f(1)+f(0) = f'' (\xi)
\geq f''(2) $
for some $ \xi \in (0,2) $.
\end{pf}

\noindent
We complete the proof of Lemma \ref{poeschel}.
We first consider the trivial cases \eqref{pampers3}-\eqref{pampers4}.
\\[1mm]
{\sc case \eqref{pampers3}.} Since $ \sum_i \s_i\neq 0 $ is even, \eqref{pecorino} follows by
$$
|\s \cdot \l_{\vec \jmath} | =| \sum_i \s_i \l_{\vec 0}| \geq 2 \l_{\vec 0} = 2 \sqrt{\mm}  \geq
\mm (1+\mm)^{-3/2}   \, .
$$
{\sc case \eqref{pampers2}.} By $ \vec \s \cdot \vec \jmath = (\s_3 + \s_4 ) q  = 0 $, $ q \neq 0 $, we deduce $ \s_3 = - \s_4 $. Hence \eqref{pecorino} follows by
$$
|\s \cdot \l_{\vec \jmath} | = | (\s_1 + \s_2 ) \l_0 | = 2 \sqrt{\mm}  \geq  \mm (1+\mm)^{-3/2}
$$
{\sc case \eqref{pampers4}.}  Since $\vec \jmath= (p,p,-p,-p)$ and $\s_1=\s_2$ then  $\vec\s\cdot \vec \jmath=0$ implies
$\s_3=\s_4=\s_2$ and
$$
|\s \cdot \l_{\vec \jmath} | = | 4 \l_p | = 4 \sqrt{p^2+\mm}  \geq  \mm (p^2+\mm)^{-3/2}  \, .
$$
{\sc case \eqref{pampers1}.} Set 
$ |j_1| =: i $ , $ |j_2| =: j $, $ |j_3| =:k $, $ |j_4|=: l  $. After reordering we can assume $ 0 \leq i\leq j\leq k\leq l.$
Since, by assumption, $ \vec \s \cdot \vec \jmath = 0 $, the following combination of
plus and minus signs gives
$ \mathtt s(j_1) \s_1 i  + \mathtt s(j_2) \s_2 j + \mathtt s(j_3) \s_3 k +  \mathtt s(j_4) \s_4  l = 0 $.
Hence Lemma \ref{Poschelzero}
implies (\ref{pecorino}) for every ${\vec \jmath}$ except when
$|j_1 |= |j_2| $ and $|j_3|=|j_4|$
(in this case $i=j$ and $k=l$ and Lemma  \ref{Poschelzero} does not apply).
We now prove that \eqref{pecorino} holds also in these cases.
We have that
$ {\vec \s}\cdot \l_{\vec \jmath} = (\sigma_1+\sigma_2)\l_{j_1}+ (\sigma_3+\sigma_4) \l_{j_3} $
 where $\s_a+\s_b=0,\pm 2$ so that \eqref{pecorino} holds trivially unless $\sigma_1+\sigma_2= -(\sigma_3+\sigma_4)$.
 We consider this last case. If $\sigma_1+\sigma_2=-(\sigma_3+\sigma_4)=0$ then  the equality
 $ \vec\s\cdot\vec\jmath= \s_1(j_1-j_2)+\s_3(j_3-j_4)=0 $
 implies that $ j_1, \ldots, j_4 $ are pairwise equal,
  contrary to our hypothesis.
 If $\sigma_1+\sigma_2=\pm 2$ and $i:=|j_1|<k:=|j_3|$ then
 $$
|{\vec \s}\cdot \l_{\vec \jmath}|\geq 2\l_{j_3}-2\l_{j_1} =2\l_k-2\l_i
\stackrel{\eqref{fave}} \geq 2\l_{k-i}-2\l_0
\stackrel{(k>i)}{\geq}
2\l_1-2\l_0\geq
1/\sqrt{1+\mm}
$$
giving (\ref{pecorino}).
If $|j_1|=|j_2|=|j_3|=|j_4|$ and  $\sigma_1+\sigma_2=-(\sigma_3+\sigma_4)=\pm 2$ then the relation
$ \vec\s\cdot\vec\jmath= \s_1 (j_1+j_2-j_3-j_4)=0$ implies  that the $j_1, \ldots , j_4 $ are pairwise equal, contrary to the hypothesis.
\rule{2mm}{2mm}

\noindent
{\it Massimiliano Berti, Michela Procesi}, Dipartimento di Matematica e Applicazioni ``R. Caccioppoli",
Universit\`a degli Studi di Napoli Federico II,  Via Cintia, Monte S. Angelo,
I-80126, Napoli, Italy,  {\tt m.berti@unina.it}, {\tt michela.procesi@unina.it}
\\[1mm]
{\it Luca Biasco}, Dipartimento di Matematica, Universit\`a di Roma 3, Largo San Leonardo Murialdo, I-00146, Roma, Italy,
{\tt biasco@mat.uniroma3.it}.
\\[1mm]
\noindent
This research was supported by the European Research Council under FP7 ``New Connections between dynamical systems
and Hamiltonian PDEs with small divisors phenomena" and partially by
the PRIN2009 grant ``Critical Point Theory and Perturbative Methods for Nonlinear Differential Equations".

\begin{thebibliography}{10}

\bibitem{BBo} Berti M., Bolle P. {\it Quasi-periodic solutions 
 with Sobolev regularity of NLS on $ \T^d $ with a multiplicative potential}, to appear on Eur. Jour. Math.

\bibitem{BB10}  Berti M., Biasco L., {\it Branching of Cantor manifolds of elliptic tori
and applications to  PDEs}, Comm. Math. Phys, 305, 3,  741-796, 2011.

\bibitem{B1}  Bourgain J., {\it  Construction of quasi-periodic solutions
for Hamiltonian perturbations of linear equations and applications
to nonlinear PDE}, Internat. Math. Res. Notices, no. 11, 1994.

\bibitem{B3}  Bourgain J., {\it Quasi-periodic solutions of Hamiltonian
perturbations of $2D$ linear Schr\"odinger equations},
Annals of Math. 148, 363-439, 1998.

\bibitem{B2} Bourgain J., {\it Periodic solutions of nonlinear wave
equations},  Harmonic analysis and partial differential equations,
69--97, Chicago Lectures in Math., Univ. Chicago Press, 1999.

\bibitem{B5}  Bourgain J.,  {\it Green's function estimates for lattice Schr\"odinger
operators and applications}, Annals of Mathematics Studies 158,
Princeton University Press, Princeton, 2005.


\bibitem{CY} Chierchia L., You J., {\it KAM tori for 1D
nonlinear wave equations with periodic boundary conditions}, Comm. Math. Phys. 211, 497-525, 2000.

\bibitem{C} Craig W.,  {\it Probl\`emes de petits diviseurs dans
les \'equations aux d\'eriv\'ees partielles},
Panoramas et Synth\`eses, 9,
Soci\'et\'e Math\'ematique de France, Paris, 2000.

\bibitem{CW} Craig W., Wayne C. E., {\it Newton's method and periodic solutions
of nonlinear wave equation}, Comm. Pure  Appl. Math. 46, 1409-1498, 1993.


\bibitem{DS} Delort J.M., Szeftel J., {\it Long-time existence for small data of nonlinear
Klein-Gordon equations on tori and spheres}, Int. Math. Res. Notices, 37, 1897-1966, 2004.

\bibitem{De} Delort J.M., {\it A quasi-linear Birkhoff normal forms method. Application to the quasi-linear
Klein-Gordon equation on $ S^1 $}, Ast\'erisque, Volume 2009/2010, expos\'es 1012-1026.





\bibitem{EK} Eliasson L.H., Kuksin S.,
{\it KAM for non-linear Schr\"odinger equation},  Annals of Math., 172, 371-435,
2010.

\bibitem{GY} Geng J., You J.,
{\it A KAM theorem for Hamiltonian partial differential equations in higher dimensional spaces},
Comm. Math. Phys. 262, 343-372, 2006.

\bibitem{GYX}
Geng J., Xu X., You J.,
{\it An infinite dimensional KAM theorem and its application to the two dimensional cubic Schr\"odinger equation},
Adv. Math. 226, 6, 5361-5402, 2011.

\bibitem{GG}
G\'{e}rard P., Grellier S.,
{\it Effective integrable dynamics for some nonlinear wave equation},
preprint 2011, \texttt{http://hal.archives-ouvertes.fr/hal-00635686/fr/}

\bibitem{GT}  Grebert B.,  Thomann L., {\it KAM for the
quantum harmonic oscillator}, Comm. Math. Phys. 307, 2, 383-427, 2011.


\bibitem{KaP} Kappeler T., P\"{o}schel J.,
{\it KAM and KdV}, Springer, 2003.

\bibitem{Ku}  Kuksin S.,
{\it Hamiltonian perturbations of infinite-dimensional linear
systems with imaginary spectrum},
Funktsional. Anal. i Prilozhen. 21, no. 3, 22--37, 95, 1987.

\bibitem{K2} Kuksin S., {\it A KAM theorem for equations of the Korteweg-de Vries type},
Rev. Math-Math Phys., 10, 3, 1-64, 1998.

\bibitem{k1}  Kuksin S.,
{\it Analysis of Hamiltonian PDEs},
Oxford Lecture Series in Mathematics and its Applications, 19.
Oxford University Press, 2000.

\bibitem{KP} Kuksin S., P\"oschel J., {\it
Invariant Cantor manifolds of quasi-periodic oscillations
for a nonlinear Schr\"{o}dinger equation},
Ann. of Math, 2,  143,  no. 1, 149-179, 1996.


\bibitem{LY}  Liu J., Yuan X.,
{\it A KAM Theorem for Hamiltonian Partial Differential
Equations with Unbounded Perturbations}, Comm. Math. Phys,
  307 (3), 629-673, 2011.


\bibitem{Po}  P\"oschel J., {\it On elliptic lower dimensional tori in Hamiltonian systems}, Math. Z., 202, 559-608, 1989.


\bibitem{Po2} P\"oschel J.,
{\it A KAM-Theorem for some nonlinear PDEs},  Ann. Sc. Norm. Pisa, 
23, 119-148, 1996.

\bibitem{Po3} P\"oschel J., {\it Quasi-periodic solutions for
a nonlinear wave equation}, Comment. Math. Helv.,  71,
no. 2, 269-296, 1996.

\bibitem{PT} P\"oschel J., Trubowitz E.,  {\it Inverse spectral theory}, Academic Press, Orlando,  1987.

\bibitem{PX} Procesi M., Xu X., {\it Quasi-T\"oplitz Functions in KAM Theorem}, preprint 2011.

\bibitem{W2} Wang W. M.,
{\it Supercritical nonlinear Schr\"odinger equations I: quasi-periodic solutions},  2011.

\bibitem{W1}  Wayne E., {\it Periodic and quasi-periodic solutions
of nonlinear wave equations via KAM theory}, Comm. Math. Phys. 127, 479-528, 1990.


\end{thebibliography}
\end{document}